\documentclass[9pt]{extarticle}
\usepackage[utf8]{inputenc}
\usepackage{graphicx}
\usepackage{cite}
\usepackage{comment}
\usepackage{appendix}
\usepackage{color}
\usepackage{amsfonts,amsmath,amssymb,amsthm}
\usepackage{float}
\newtheorem{theorem}{Theorem}
\usepackage{lipsum}
\usepackage{caption}
\usepackage{subcaption}
\usepackage{multicol}
\usepackage[left=1.5 cm,right=1.2 cm,top=1.5cm,bottom=1.2cm]{geometry}
\numberwithin{equation}{section}
\title{{\Large Predator-prey dynamics pertaining to structuralizing predator species into three stages coupled with maturation delay owing to juvenile hunting}}
\author{Debasish Bhattacharjee$^1$, Tapasvini Roy$^2$, Santanu Acharjee$^3$, Tarini Kumar Dutta$^4$ \\$^{1,2,3}$Department of Mathematics, Gauhati University, Assam, India\\
$^4$Department of Mathematics, Assam Don Bosco University, Assam, India\\
e-mails:$^1$debabh2@gmail.com, $^2$tapasviniroy@gmail.com, $^3$sacharjee326@gmail.com, $^4$tkdutta2001@yahoo.co.in}
\date{}
\begin{document}
\maketitle
\begin{abstract}
    The predator-prey dynamic appertaining to two species is explored, wherein the predator species is structured into different stages. As evidenced from natural documentation, the immature predators possess the potential to predate albeit not as competently as the adults. Nevertheless, this potentiality is not acquired immediately after their incipience of life, hence, the immature stage is branched off into the infant stage, the stage with extensive reliance on the adults, and the juvenile stage, the stage with the potential to predate but not to procreate. In this paper, this inaugural concept is coupled with injuries in the juvenile stage as the repercussion of their incompetency in predating, thereby ensuing a delay in their maturation. With the incentive to investigate the ascendancy of these refinements over the whole system, stability analyses along with various bifurcation analyses around the equilibrium points of the system are corroborated. In addition to Hopf, transcritical, and saddle node bifurcations, the existence of Bogdanov-Takens point, cusp point, Bautin bifurcation point, ‘bloom’ phenomenon, twice occurring Hopf bifurcation, and bi-stability phenomenon make the paper appreciably more rich and efficacious. 
\end{abstract}
\textbf{2020 AMS Classifications:} 34C23,92D25,92D40.\\
\textbf{Keywords: } Population dynamics; dangerous prey; limit cycle.
\section{Introduction}
After the inaugural work of Lotka and Volterra \cite{p1,p2}, a wealth of theoretical and pragmatical research has been deliberated upon the dynamical synergy among the interacting species regarding predation, competition, parasitism, and other such analogous symbiosis, owing to its ubiquitous influence on other copious aspects of the natural world. For representing these population dynamics of the interactions, mathematical models are contrived, based upon the experiments and surveillance, so as to prognose the distribution of population and structure of communities, and also how miniscule alteration in biological or environmental parameters effectuates the entire bio-system. 
\par The standard Lotka-Volterra mathematical model has been found to be lacking in several aspects. With the exertion of many researchers, the models epitomizing the predator-prey relationship have become a lot more realistic and precise by assimilating instances  such as, time-delay, functional response, stage-structure, and others into the models. Stage structured mathematical models paint more realistic and precise scenarios as compared to that of unstructured ones; after the assimilation, the models are known to show substantial changes in their behaviors. Owing to the improvisation on the assumption that each and every individual of the predator species has the same ability for predation and survivability, the stage-structured predator-prey models have received much attention \cite{ppr1,ppr2,ppr3,ppr4,ppr5,ppr6}. The rudimentary works of constructing the mathematical models with stage structures are done by Gurney \cite{p3}, Nisbet et al \cite{p4,p5}, and their works have been further progressed by Aiello et al \cite{p6,p7}.
\par The life-history of predators is always divided into two stages: the immature stage  and the matured stage. The immature predators are predominantly assumed to not being able to hunt or reproduce and are exhaustively dependent on the adults for their sustainability. However, there are several documentations of predators hunting as immatures, albeit not as efficiently as the matured ones. Caimans\cite{caiman}, pumas\cite{puma}, snakes\cite{snake}, snails\cite{snail}, sharks\cite{shark}, \textit{Clinus superciliosus}\cite{clinus}, spiders\cite{spider}, and raptors\cite{raptor} are only a handful among the myriad of such predators that hunt during their developing period of lives. At the same time, it is to be noted that the immature predators do not instantaneously start hunting from the incipience of their life. Puma \textit{(Puma concolor)} kittens are dependent on their mothers  upto two years, before dispersing away as juveniles or sub-adults \cite{puma}, similarly snakes, spiders and the likes begin their life as an egg. The same fact of there being a period of time before the juveniles of all other predator species start hunting, can be unearthed. Thereupon, while considering predation by immature predators, it is quite logical to further partition the immature stage into two stages: the infancy stage and the juvenile stage. Eliciting from all these, we have reckoned with the idea of structuring the predator species into three stages: a stage with its individuals unmitigatedly reliant on the adults (infant stage), a stage whose individuals are equipped for only predation (juvenile stage), and a stage with its individuals equipped with the ability to both predate and proliferate (adult stage); which (to the best of our knowledge) has not yet been done.\\
\par Predatism, being a highly complex and dynamic interaction, is not always a simple case scenario of predator successfully killing their prey. Prey exhausts all of their morphological and behavioural traits that enhance their survivability during their encounter with the predators.  There are countless instances where failure on the part of predators can cost more than just a missed opportunity to feed, specifically yielding injury and sometimes even death. Caribou and other adult ungulates, pose an injury risk to coyotes, despite this, coyotes single out these large prey owing to high energetic gains \cite{a1}. The same scenario is witnessed between the wolves and their prey, elk and bison \cite{a2}. Mountain lions may be pierced by porcupine quills or left blind in one eye from the blow of deer hoof \cite{a3}. Jaguars are frequently seen with fresh wounds owing to the defensive attack by peccary or giant anteater \textit{(Myrmecophaga tridactyla)} \cite{a4}. Praying Mantids vomited and showed signs of poisoning after feeding on $O.fasciatus$ \cite{a5}.
\par The predators that include dangerous prey in their diet evolve their hunting strategy accordingly so as to mitigate the risk of injury \cite{muk}. Interestingly, it is found that the younger the predator is, the greater would be the prospect of catastrophe of predation, owing to their relatively less body size and inexperience when compared to that of the adult individuals of their population \cite{b1}. The ingenuousness of the juveniles towards the specific approach orientations required while handling dangerous prey leads to them being injured, and sometimes even killed. Thereon, the predatory behavior and its effectiveness varies with ontogeny. 
\par Ants, considered as dangerous prey for spiderlings because of their mandibles and the secretion of noxious substances by stinging, lacerates the spiderlings, which has serious implications on autonomy. Spiderlings change their attacking tactic through experience to that of the adults, thereby capturing their prey faster and more efficiently \cite{ex1,ex2}.  Centipedes, capable of delivering extremely damaging and painful strikes using their forcipules, are known to be fed upon by snakes, which hunt as juveniles \cite{snake}. Juvenile migratory birds, being inadequate, require stop-over sites to replenish their fuel reserves. The unfamiliar terrain as well as handling and capturing unfamiliar prey poses a huge risk to them \cite{ex3}. Serious injuries associated with prey capture have been witnessed among the likes of big cats such as lions, tigers, jaguars, cougars, etc, especially in the juveniles of the populations \cite{b1}. The injuries suffered by the juvenile predators  as a consequence of their incompetency may lead to them being unwilling/unable to hunt while nursing the wound or some other foraging costs similar to those arising from the risk of predation, and with similar consequences, which inevitably would lead to a delay in their transition into the adult stage.\\
\par In \cite{ppr1,ppr2,ppr3,ppr4}, the authors have studied the mathematical models with general stage structure in the predator species, where the immature predators are plagued with mortal perilousness while preying on dangerous prey. But in this paper, rather than death, the concentration is solely kept on injury as being the negative repercussion of juvenile hunting (as the frequency of injuries is much higher than that of death), which is yet to be probed into mathematically. We have endeavoured  towards assembling a system of equations where only the adult predators  are fully competent in predating and can reproduce;  the incompetency  of the juvenile predator, while predating, in handling of the dangerous prey causes a delay in their transition rate; during the infancy period of life, the predators are completely dependent on the adults.  The structuralization of predator species and the inclusion of injury induced delay import the significance of each stage, on the grounds of their individually disparate impact on the whole ecosphere. 
\par In this paper, the construction of the equations of bio-system is done in section 2 and its enclosedness is shown in section 3.  The discussion of the equilibrium points in conjunction with their local and global stabilities is done in section 4; the same is done for Hopf, saddle node, and transcritical bifurcations in section 5. Numerical simulation is done in section 6, and finally the conclusion is given in section 7. 
\section{Construction of the equations of bio-system}
To begin with, we have considered the prey species to have a logistic growth rate, a biologically feasible scenario owing to the constraints on natural resources, along with the existence of the predator species in the same habitat as that of the former. The functional response, which symbolises the interaction between both the species is taken to be Holling-type 1 functional response. 
Incorporating stage structure of the predator species into the mathematical model in the generally used manner \cite{khajan2, p-ex1, p-ex2}, and taking $X(T)$ as the biomass density of the prey and $Y_i(T)$, and $Y_m(T)$ as the immature and matured predators’ biomass densities respectively, we have, 
\begin{equation}
\begin{split} \label{e1}
    \frac{dX}{dT}&=r X(1-\frac{X}{K})-AXY_m,\\
\frac{dY_i}{dT}&=U_1 AXY_m-B_1Y_i-D_1Y_i,\\
\frac{dY_m}{dT}&=B_1 Y_i-D_2 Y_m,\\
&X(0)>0, Y_i(0)>0, Y_m(0)>0.
\end{split}
\end{equation}
Here, matured predators can kill their prey hence, their rate of predation is taken as $A$, also they possess the reproductive ability which can be translated as the transformation of consumed prey into immature predators, $U_1$ denotes this rate of transformation. $D_1,\;\text{and } D_2$ are the death rates of immature and matured predators, respectively, and $B_1$ is the transition rate of immature predators into the matured stage, which is taken as constant when the nutrient required by the immature predators is not dependent on the prey. Even among the specialist predator, diet is known to vary during ontogeny, and the immatures can comfortably survive for a time being on alternate prey with low energetic returns.  Conjointly, the immature predator is completely dependent on their adults. \\ 
\par Eliciting from all the theoretical results as given in the introduction, we have partitioned the immature stage of predator into infancy and juvenile stages. The infancy stage is similar to that of the immature stage of \ref{e1}, as the individuals of this stage are exhaustively reliant on the adults. While the juvenile predators can hunt, they are yet to obtain the ability to reproduce and hence, their predation do not correspond to an increase in their total numbers. Thus, the mathematical model \ref{e1}  is remodelled into the following system of equations:
\begin{equation}
    \begin{split} \label{e2}
    \frac{dX}{dT}&=rX(1-\frac{X}{K})-A_1XY_2-A_2XY_3,\\
\frac{dY_1}{dT}&=uA_2XY_3-BY_1-D_1Y_1,\\
\frac{dY_2}{dT}&=BY_1-CY_2-D_2Y_2,\\
\frac{dY_3}{dT}&=CY_2-D_3Y_3,\\
&X(0)>0, Y_1(0)>0,\, Y_2(0)>0,\, Y_3(0)>0.
\end{split}
\end{equation}
Here, $Y_1=$ infant predators, $Y_2=$ juvenile predators, $Y_3=$ adult predators, $u=$ transformation rate of consumed prey into infant predators,
 $A_1$= predation rate of juvenile predators, $A_2=$ predation rate of adult predator. $D_1=$ natural death rates of infant predators, $D_2=$ natural death rates of juvenile predators, $D_3=$ natural death rates of adult predators, $B=$ transition rate from infancy to juvenile, $C=$ transition rate from juvenile to adulthood.
 \par Due to lack of efficiency in predating by the juveniles, consequential to lack of robustness and experience, they are faced with injuries. The cost of being injured leads to them taking a reprieve from predation or from hunting their required amount of nutritive prey, resulting in causing a delay in their transition rate. This delay has been incorporated into the model \ref{e2}, by remoulding $C$ into $(C-A_3X)$, where $A_3$ is the rate of injury caused by the prey. \\
Therefore, our equations of bio-system is:
\begin{equation}
    \begin{split} \label{eq1}
        \frac{dX}{dT}&=rX(1-\frac{X}{K})-A_1XY_2-A_2XY_3,\\
        \frac{dY_1}{dT}&=uA_2XY_3-BY_1-D_1Y_1,\\
        \frac{dY_2}{dT}&=BY_1-(C-A_3X)Y_2-D_2Y_2,\\
        \frac{dY_3}{dT}&=(C-A_3X)Y_2-D_3Y_3,\\
        &\text{With initial conditions: } X(0)>0,\hspace{3pt}Y_1(0)>0,\hspace{3pt}Y_2(0)>0,\hspace{3pt}Y_3(0)>0.
    \end{split}
\end{equation}
To make the model simpler for working out various calculations, the following steps are taken:\\
$t=rT,\hspace{3pt}x=X/K,\hspace{3pt}y_1=Y_1/K,\hspace{3pt}y_2=Y_2/K,\hspace{3pt}y_3=Y_3/K$ along with $a_1=A_1\frac{K}{r},\hspace{3pt}a_2=A_2\frac{K}{r},\hspace{3pt}b=B/r,\hspace{3pt}d_1=D_1/r,\hspace{3pt}c=C/r,\hspace{3pt}a_3=A_3\frac{K}{r},\hspace{3pt}d_2=D_2/r,\hspace{3pt}d_3=D_3/r$ and then using them in \ref{eq1} we have,
\begin{equation}\label{ma-eq}
    \begin{split}
        \frac{dx}{dt}&=x(1-x)-a_1xy_2-a_2xy_3,\\
        \frac{dy_1}{dt}&=ua_2xy_3-by_1-d_1y_1,\\
        \frac{dy_2}{dt}&=by_1-(c-a_3x)y_2-d_2y_2,\\
        \frac{dy_3}{dt}&=(c-a_3x)y_2-d_3y_3,\\
        &\text{With initial conditions: } x(0)>0,\hspace{3pt}y_1(0)>0,\hspace{3pt}y_2(0)>0,\hspace{3pt}y_3(0)>0.
   \end{split}
\end{equation}
\section{Boundedness of the equations of bio-system}
From the first equation of (\ref{ma-eq}), we can observe $\Dot{x}\leq x(1-x)$, i.e.,
$\text{Hence }\lim_{t\to\infty}\, sup \,x(t)\leq 1$.\\
Next, we consider a variable $\beta$ in the following manner,
\begin{align*}
    \beta=&ux+y_1+y_2+y_3,\\
    &\text{differentiating $\beta$ with respect to time t and then adding $\zeta \beta$, where $\zeta=min\,\{d_1,d_2,d_3\}$, we have, }\\
    \frac{d\beta}{dt}+\zeta \beta=& ux(1-x)-ua_1xy_2+\zeta ux-(d_1-\zeta)y_1-(d_2-\zeta)y_2-(d_3-\zeta)y_3 \leq ux(1-x+\zeta)\\
    i.e.\,\frac{d\beta}{dt}+\zeta\beta\leq & u(\frac{1+\zeta}{2})^2\\
    i.e.,\beta\leq & u\frac{(1+\zeta)^2}{4\zeta}(1-e^{-\zeta t})+\beta(x(0),y_1(0),y_2(0),y_3(0)) e^{-\zeta t}\\
    \therefore \lim_{t \to \infty}\, sup \, \beta(t)& \leq \frac{u(1+\zeta)^2}{4\zeta}\\
    & \,i.e.\,\lim_{t \to \infty}\, sup \,(y_1(t)+y_2(t)+y_3(t))\leq\frac{u(1+\zeta)^2}{4\zeta}-u
    \leq\frac{u(1-\zeta)^2}{4\zeta}.
\end{align*}
Therefore, we arrive at the following theorem.
\begin{theorem}
Every solution of the equations of bio-system (\ref{ma-eq}) originating from the region $R_+^4$ is confined in\\
$\{(x,y_1,y_2,y_3):0\leq x\leq1,0\leq (y_1+y_2+y_3)\leq \frac{u(1-\zeta)^2}{4\zeta}\}$ where $\zeta=min\{d_1,d_2,d_3\}$.
\end{theorem}
\section{Equilibrium points and their stabilities}
Equilibrium points are the stationary points experienced by the system of equations where the augmentation and reduction velocities of the species vanish. Our given equations of bio-system (\ref{ma-eq}) is encountered with four types of such situations. First of them is the extinction equilibrium point $E_1(0,0,0,0)$, and it always exists. In the absenteeism of predators, the prey can attain the highest ecological load of its environment, resulting in the axial equilibrium $E_2(1,0,0,0)$, the prey-only equilibrium point, which too always exists, also with the predator being a specialist type, other axial equilibriums are not possible. The most advantageous and desirable case scenario in the ecological system is the co-existence of all the species.  The equations of bio-system (\ref{ma-eq}) manifests two compresent equilibrium points $E_3(x^*,y_1^*,y_2^*,y_3^*)$ and $E_4(x^{**},y_1^{**},y_2^{**},y_3^{**})$. Stability analysis is requisite for the alteration of behavior of the equilibrium points into a desirable one. With this inducement, the theorems in this section have been deliberated upon.
\subsection{Extinction equilibrium point $E_1(0,0,0,0)$} 
The equilibrium point $E_1(0,0,0,0)$ is a saddle point because:
\begin{align}\label{jaco}
J= \left(
\begin{array}{cccc}
 -2 x-a_1 y_2-a_2 y_3+1 & 0 & -a_1 x & -a_2 x \\
 a_2 u y_3 & -b-d_1 & 0 & a_2 u x \\
 a_3 y_2 & b & -c-d_2+a_3 x & 0 \\
 -a_3 y_2 & 0 & c-a_3 x & -d_3 \\
\end{array}
\right)\end{align}
is the Jacobian matrix of the equations of bio-system (\ref{ma-eq}),
the characteristic polynomial of which at $(0,0,0,0)$ is\\
$(1 - \omega) (\omega+b + d_1) (c d_3 + d_2 d_3 + c \omega + d_2 \omega + d_3 \omega + \omega^2)=0, \; \; 
i.e., \; \; \omega=1,-b-d_1,-c-d_2,-d_3$.\\
$\therefore$ All the eigen values are negative except one.
\subsection{Prey-only equilibrium point $E_2(1,0,0,0)$}
\begin{theorem} \label{thm-lo-pf}
The equilibrium point $E_2(1,0,0,0)$ can be either a center, stable point, or an unstable point. All the conditions required for those related to stability are given within the following proof.
\end{theorem}
\textbf{Proof:} Putting $(x,y_1,y_2,y_3)= (1,0,0,0)$ at the Jacobian matrix (\ref{jaco}), we arrive at the required matrix $J_1$, the determinant and trace of whose are:\\ 
$Det\,[J_1]= d_1(c+d_2) d_3 + b (c d_3 + d_2 d_3 -a_2 c u)-a_3 (d_1 d_3+ b (d_3-a_2 u))$, and\\
$Trace\,[J_1]=a_3-(1 + b + d_1 + c + d_2 +d_3) $.\\
For (1,0,0,0) to be a saddle point, $Det\,[J_1]<0$ is required,
$$i.e., \hspace{3pt}  either.\hspace{5pt}  a_3>c+d_2,a_2<\frac{(b+d_1)(a_3-c-d_2)d_3}{b(a_3-c)u}\hspace{5pt} or\hspace{4pt}a_3<c,a_2>\frac{(b+d_1)(a_3-c-d_2)d_3}{b(a_3-c)u}.$$
For (1,0,0,0) to be a unstable point, $Det\,[J_1]>0$, and $Trace\,[J_1]>0$ are required,
$$i.e.,\hspace{3pt} a_3>1+b+c+d_1+d_2+d_3,\text{ and }a_2>\frac{(b+d_1)(a_3-c-d_2)d_3}{bu(a_3-c)}.$$
For (1,0,0,0) to be a centre, $Det\,[J_1]>0$, and $Trace\,[J_1]=0$ are required,
$$i.e.,\hspace{3pt}  a_2>\frac{(b+d_1)(1+b+d_1+d_3)d_3}{b(1+b+d_1+d_2+d_3)u},\text{ and }a_3=1+b+c+d_1+d_2+d_3.$$
For (1,0,0,0) to be a stable point, $Det\,[J_1]>0$, and $Trace\,[J_1]<0$ are required, resulting the following cases:
\begin{align*}
    Case1:&\hspace{5pt}a_3<c,\text{ and }a_2<\frac{(b+d_1)(a_3-c-d_2)d_3}{bu(a_3-c)}, \\
    Case2:&\hspace{5pt}a_3=c,\\
    Case3:&\hspace{5pt}c<a_3<1+c,\text{ along with, either }\left (a_2>\frac{(b+d_1)(a_3-c-d_2)d_3}{bu(a_3-c)},a_3>c+d_2\right)\text{ or }\left(a_3\leq c+d_2\right),\\
    Case4:&\hspace{5pt} a_3>1+c,a_2>\frac{(b+d_1)(a_3-c-d_2)d_3}{bu(a_3-c)},\text{ and } a_3<1+c+b+d_1+d_2+d_3.
\end{align*}
\begin{figure}[H]
    \centering
    \begin{subfigure}[b]{0.44\textwidth}
         \centering
    \includegraphics[width=\textwidth]{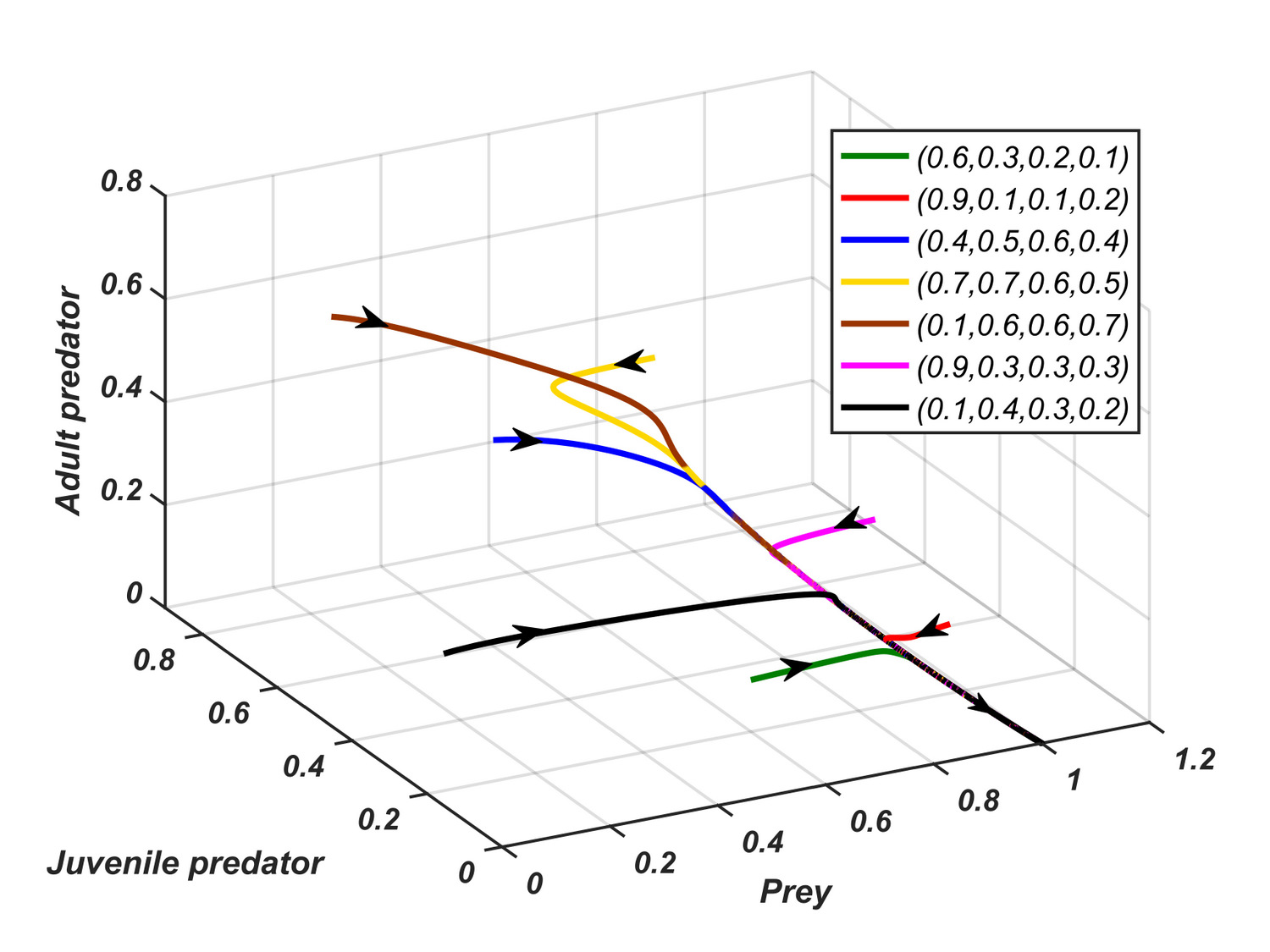}
    \caption{}
    \label{lo-axi-1}
    \end{subfigure}\hfill
    \begin{subfigure}[b]{0.44\textwidth}
         \centering
    \includegraphics[width=\textwidth]{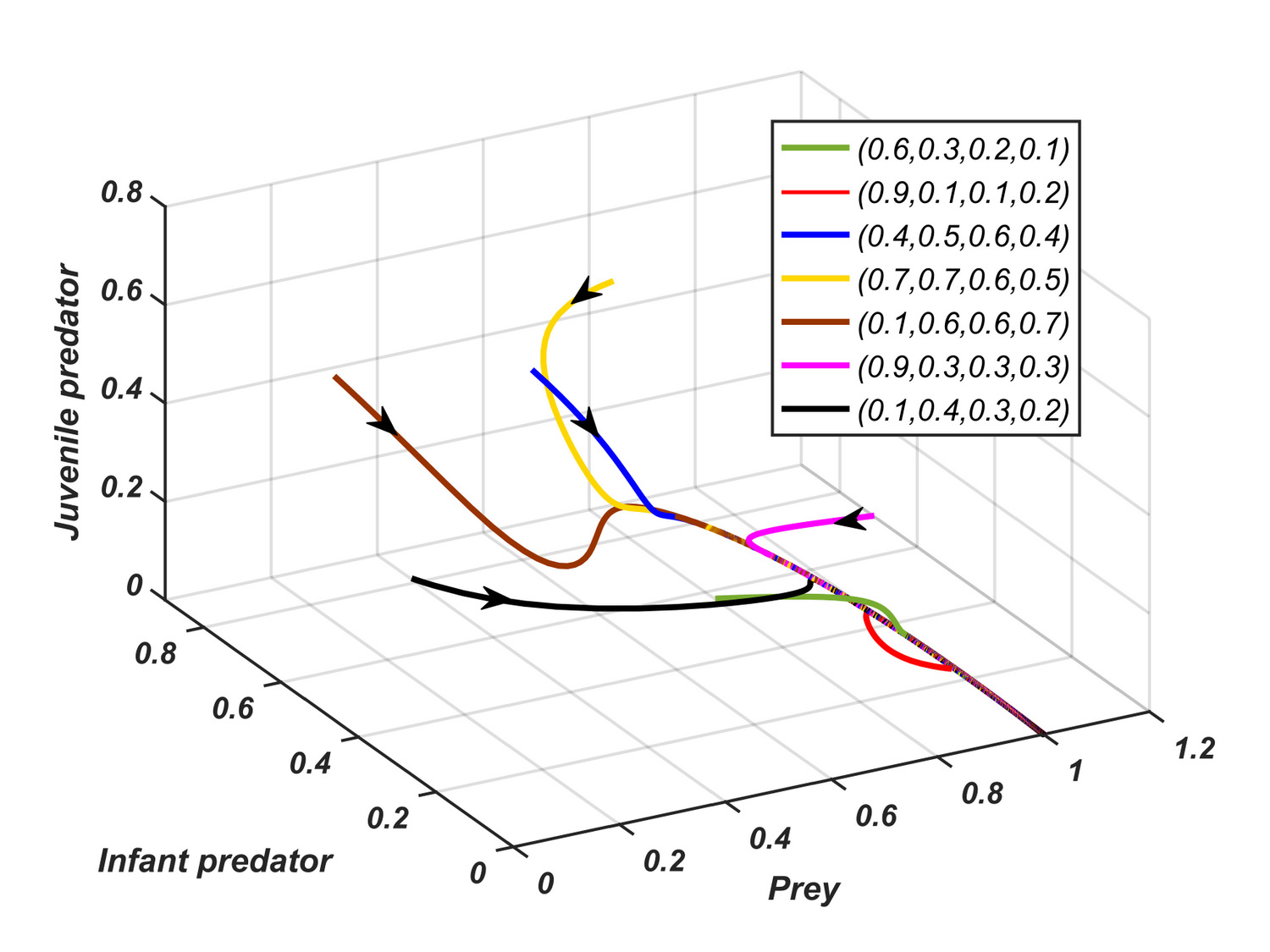}
    \caption{}
    \label{lo-axi-2}
    \end{subfigure}
    \caption{Local and global stability of prey-only equilibrium point $(1,0,0,0)$: (a)depicts extinction of juvenile and adult predators, (b)depicts extinction of infant and juvenile predators.}
    \label{lo-axi}
\end{figure}
\begin{theorem} \label{thm-glo-pf}The prey-only equilibrium point $E_2(1,0,0,0)$ is globally asymptotically stable when  $\;(1-x)(a_1y_2+x-1)+a_2y_3<0$. \end{theorem}
\textbf{Proof:} Considering the Lyapunov function $L=\tau_1(x-x_1-x_1ln\frac{x}{x_1})+\tau_2 y_1+\tau_3 y_2+\tau_4 y_3$ so as to investigate the global stability of prey-only equilibrium point $E_1(1,0,0,0)$, where $\tau_i (i=1,2,3,4)$ are positive constants, which are to be taken accordingly.\\
Differentiating the Lyapunov function with respect to time t, we have,
\begin{align*}
    \frac{dL}{dt}=&\tau_1 \frac{x-x_1}{x}\dot{x}+\tau_2 \dot{y_1}+\tau_3 \dot{y_2}+\tau_4 \dot{y_3}\\
    &\text{Now, putting $x_1=1$ and $\tau_1=u\tau_2, \tau_2=\tau_3=\tau_4$, we have,}\\
   \frac{dL}{dt} =&\tau_2(-(1-x)^2-a_1xy_2+a_1y_2+a_2y_3-d_1y_1-d_2y_2-d_3y_3)\\
\leq &\tau_2\{(1-x)(a_1y_2-1+x)+a_2y_3\}.
\end{align*}
Therefore, the sufficient condition for the equilibrium point to be globally asymptotically stable would be:
$$(1-x)(a_1y_2+x-1)+a_2y_3<0.$$
\subsection{Compresent equilibrium points $E_3(x^*,y_1^*,y_2^*,y_3^*)$ and $E_4(x^{**},y_1^{**},y_2^{**},y_3^{**})$}
The two compresent equilibrium points manifested by the equations of bio-system (\ref{ma-eq}) are:
\begin{align*}
&E_3(x^*,y_1^*,y_2^*,y_3^*)=
\left(\frac{\delta_3 + \mu}{ 2 a_2 a_3 b u}, -\frac{ (\delta_4- 2 a_2 b d_2 u + \mu) \mu_1}{ 4 a_2^2 a_3 b^2 u \delta_1+ \delta_2)}, \frac{\mu_1}{ 2 a_2 a_3 \delta_1+ \delta_2)},-\frac{ (\delta_4+\mu) \mu_1}{ 4 a_2^2 a_3 b d_3 u (\delta_1+\delta_2)}\right),\\
&\\
&E_4(x^{**},y_1^{**},y_2^{**},y_3^{**})=
\left(\frac{\delta_3 - \mu}{ 2 a_2 a_3 b u}, \frac{ (-\delta_4+ 2 a2 b d_2 u + \mu) \mu_2}{ 4 a_2^2 a_3 b^2 u (\delta_1+\delta_2)},\frac{\mu_2}{ 2 a_2 a_3 (\delta_1+\delta_2)},\frac{ (-\delta_4+ \mu) \mu_2 }{ 4 a_2^2 a_3 b d_3 u (\delta_1 + \delta_2)}\right),
\end{align*}
where, $\mu= \sqrt{a_3^2 (b+d_1)^2 d_3^2-2 a_2 a_3 b (b+d_1) (c+2 d_2) d_3 u+a_2^2 b^2 c^2 u^2}$,\\
$\mu_1=a_2^2 b (a_3 - c) c u – a_1 d_3 (a_3 (b + d_1) d_3 + \mu) + a_2 (-a_3^2 (b + d_1) d_3 - c (a_1 b d_3 u + \mu) + a_3 (c d_1 d_3 + 2 d_1 d_2 d_3 + b d_3 (c + 2 (d_2 + a_1 u)) +\mu))$,\\
$\mu_2= a_2^2 b (a_3 - c) c u + a_1 d_3 (-a_3 (b + d_1) d_3 +\mu) + a_2 (-a_3^2 (b + d_1) d_3 + a_3 (c d_1 d_3 + 2 d_1 d_2 d_3 + b d_3 (c + 2 (d_2 + a_1 u)) -\mu) + c (-a_1 b d_3 u +\mu)),$\\
$\delta_1=(a_1 d_3 (-a_3 (b + d_1) + a_1 b u), \; \delta_2=a_2 (a_3 (b + d_1) d_2 + a_1 b c u),\; \delta_3=a_3 (b + d_1) d_3 + a_2 b c u, \;\text{and } \delta_4=a_3 (b + d_1) d_3 – a_2 b c u $.\\ \\
\textbf{Existence}\vspace{0.16cm}\\
For both the compresent equilibrium points $E_3$, and $E_4$ to exist, the following necessary and sufficient conditions should be satisfied:
\begin{align*}
& \hspace{1cm}d_2<\frac{(a_3 (b+d_1) d_3-a_2 b c u)^2}{4 a_2 a_3 b (b+d_1) d_3 u} \text{ along with}\\
    Case1:\hspace{5pt}&a_3=c, \text{ and } a_2>\frac{a_3 b d_3+a_3d_1d_3}{2a_3bu-bcu},\\
Case2:\hspace{5pt}&a_3>c,\text{ and } a_2bcu>a_3(b+d_1)d_3,\\
Case3: \hspace{5pt} &a_3<c,\;2a_3>c,\;  a_2>\frac{a_3  d_3(b+d_1)}{2a_3bu-bcu},  \text{ and } \frac{(a_3-c)(d_1d_3+b(d_3-a_2u))}{(b+d_1)d_3}<d_2.
\end{align*}
\begin{theorem} \label{sta-co}
The compresent equilibrium points $E_3$ and $E_4$ are locally asymptotically stable iff
$\epsilon_1>0,\;\epsilon_4>0,\;\epsilon_1\epsilon_2-\epsilon_3>0, \;\text{ and } \; \epsilon_1\epsilon_2\epsilon_3-\epsilon_3^2-\epsilon_1^2\epsilon_4>0.$ \hspace{0.2cm}
The values of $\epsilon_1,\, \epsilon_2, \, \epsilon_3, \, \epsilon_4$ are explained within the proof.
\end{theorem}
\textbf{Proof:} The characteristic equation of the Jacobian (\ref{jaco}) is:
\begin{align*}
    a_2 b u x &[a_3 x (-1+\omega+2 x)-c (-1+\omega+2 x+a_1 y_2)+a_1 (d_3+\omega) y_3]+(-b-d_1-\omega) [a_2 a_3 (d_2+\omega) x y_2\\
    &+(-d_3-\omega) (a_1 a_3 x y_2+(c+d_2+\omega-a_3 x) (-1+\omega+2 x+a_1 y_2+a_2 y_3))]=0.
\end{align*}
This is rewritten as:
\begin{align} \label{ch-eq}
\omega^4+\epsilon_1\omega^3+\epsilon_2\omega^2+\epsilon_3\omega+\epsilon_4=0.
\end{align}
Using the Routh-Hurwitz criteria, and putting either $(x,y_1,y_2,y_3)=(x^*,y_1^*,y_2^*,y_3^*)$, or $(x,y_1,y_2,y_3)=(x^{**},y_1^{**},y_2^{**},y_3^{**})$ according to the requirement, we would have $E_3$ or $E_4$ to be locally asymptotically stable iff 
$\epsilon_1>0,\hspace{7pt}\epsilon_4>0,\hspace{7pt}\epsilon_1\epsilon_2-\epsilon_3>0,\hspace{7pt} \text{and }\epsilon_1\epsilon_2\epsilon_3-\epsilon_3^2-\epsilon_1^2\epsilon_4>0$.
\begin{theorem} \label{thm-glo-co}
A sufficient condition for the equations of bio-system (\ref{ma-eq}) to be globally asymptotically stable around the compresent equilibrium point  $E_3(x^*,y_1^*,y_2^*,y_3^*),$ or $E_4(x^{**},y_1^{**},y_2^{**},y_3^{**})$ is: \hspace{0.2cm}
$(x-X^*)(1-x-a_1y_2)+a_2X^*y_3-a_2xy_3 \frac{Y_1^*}{y_1}<0,\;\; \; Y_1^*(d_1+b)<d_1y_1+d_2y_2, \; \; \;  c+d_2<a_3x+by_1/y_2, \; \; \text{and }\; Y_3^*(d_3-c \frac{y_2}{y_3}+a_3x \frac{y_2}{y_3})<d_3y_3$.
\end{theorem}
\textbf{Proof: } The positive definite Lyapunov function to be considered to investigate the global stability of the system of equations around the equilibrium point $E_3(x^*,y_1^*,y_2^*,y_3^*)$ is taken as:
$$L_1=\eta_1(x-x^*-x^*ln\frac{x}{x^*})+\eta_2(y_1-y_1^*-y_1^*ln \frac{y_1}{y_1^*}) +\eta_3(y_2-y_2^*-y_2^*ln \frac{y_2}{y_2^*})+\eta_4(y_3-y_3^*-y_3^*ln \frac{y_3}{y_3^*}).$$
Here $\eta_1,\eta_2,\eta_3,\text{ and }\eta_4$ are all positive constants that are to be determined.
Similarly, for the equilibrium point $E_4(x^{**},y_1^{**},y_2^{**},y_3^{**})$, the Lyapunov function would be\\
$$L_2=\eta_1(x-x^{**}-x^{**}ln\frac{x}{x^{**}})+\eta_2(y_1-y_1^{**}-y_1^{**}ln \frac{y_1}{y_1^{**}}) +\eta_3(y_2-y_2^{**}-y_2^{**}ln \frac{y_2}{y_2^{**}})+\eta_4(y_3-y_3^{**}-y_3^{**}ln \frac{y_3}{y_3^{**}}).$$
Because of their similarities, both the Lyapunov functions can be written as follows: 
 $$L^*=\eta_1(x-X^*-X^*ln\frac{x}{X^*})+\eta_2(y_1-Y_1^*-Y_1^*ln \frac{y_1}{Y_1^*}) +\eta_3(y_2-Y_2^*-Y_2^*ln \frac{y_2}{Y_2^*})+\eta_4(y_3-Y_3^*-Y_3^*ln \frac{Y_3}{Y_3^*}),$$
 where, $(X^*, Y^*, Z^*)$ represent either of the equilibrium points $E_3$, $E_4$.\\
Therefore, differentiating $L^*$ with respect to time t, 
\begin{align*}
    \frac{dL^*}{dt}=&\eta_1\frac{x-X^*}{x}\dot{x}+ \eta_2\frac{y_1-Y_1^*}{y_1}\dot{y_1}+ \eta_3\frac{y_2-Y_2^*}{y_2}\dot{y_2}+ \eta_4\frac{y_3-Y_3^*}{y_3}\dot{y_3}\\
&\text{Taking } \eta_1=u\eta_2, \text{ and } \eta_2=\eta_3=\eta_4, \text{ we get,}\\
\frac{dL^*}{dt}=&\eta_2[u(x-x^2-X^*+xX^*-a_1xy_2+a_1X^*y_2-a_2xy_3+a_2X^*y_3)+ua_2xy_3-ua_2xy_3\frac{Y_1^*}{y_1}-by_1\\
&+bY_1^*-d_1y_1+d_1Y_1^*+by_1-by_1\frac{Y_2^*}{y_2}-(c-a_3x)y_2+(c-a_3x)Y_2^*-d_2y_2+d_2Y_2^*+(c-a_3x)y_2\\
&-(c-a_3x)y_2\frac{Y_3^*}{y_3}-d_3y_3+d_3Y_3^*]\\
=&\eta_2[\{u(x-x^2-X^*+xX^*-a_1xy_2+a_1X^*y_2+a_2X^*y_3-a_2xy_3\frac{Y_1^*}{y_1})\}+\{bY_1^*-d_1y_1+d_1Y_1^*-d_2y_2\}\\
&+\{-by_1\frac{Y_2^*}{y_2}+(c-a_3x)Y_2^*+d_2Y_2^*\}+\{-(c-a_3x)y_2\frac{Y_3^*}{y_3}-d_3y_3+d_3Y_3^*\}].
\end{align*}
Therefore, the sufficient condition requisite for $\dot{L^*}$ to be negative would be
$$ (x-X^*)(1-x-a_1y_2)+a_2X^*y_3-a_2xy_3 \frac{Y_1^*}{y_1}<0, \; \; Y_1^*(d_1+b)<d_1y_1+d_2y_2, \; \; c+d_2<a_3x+by_1/y_2, \; \text{and } \;Y_3^*(d_3-c \frac{y_2}{y_3}+a_3x \frac{y_2}{y_3})<d_3y_3.$$
Putting $(X^*,Y_1^*,Y_2^*,Y_3^*)=(x^*,y_1^*,y_2^*,y_3^*)$ or $(X^*,Y_1^*,Y_2^*,Y_3^*)=(x^{**},y_1^{**},y_2^{**},y_3^{**})$ accordingly,
we would have \\the sufficient requirements for global stability of equilibrium points $E_3(x^*,y_1^*,y_2^*,y_3^*)$ and $E_4(x^{**},y_1^{**},y_2^{**},y_3^{**})$.
\begin{figure}[H]
    \centering
    \begin{subfigure}[b]{0.44\textwidth}
         \centering
    \includegraphics[width=\textwidth]{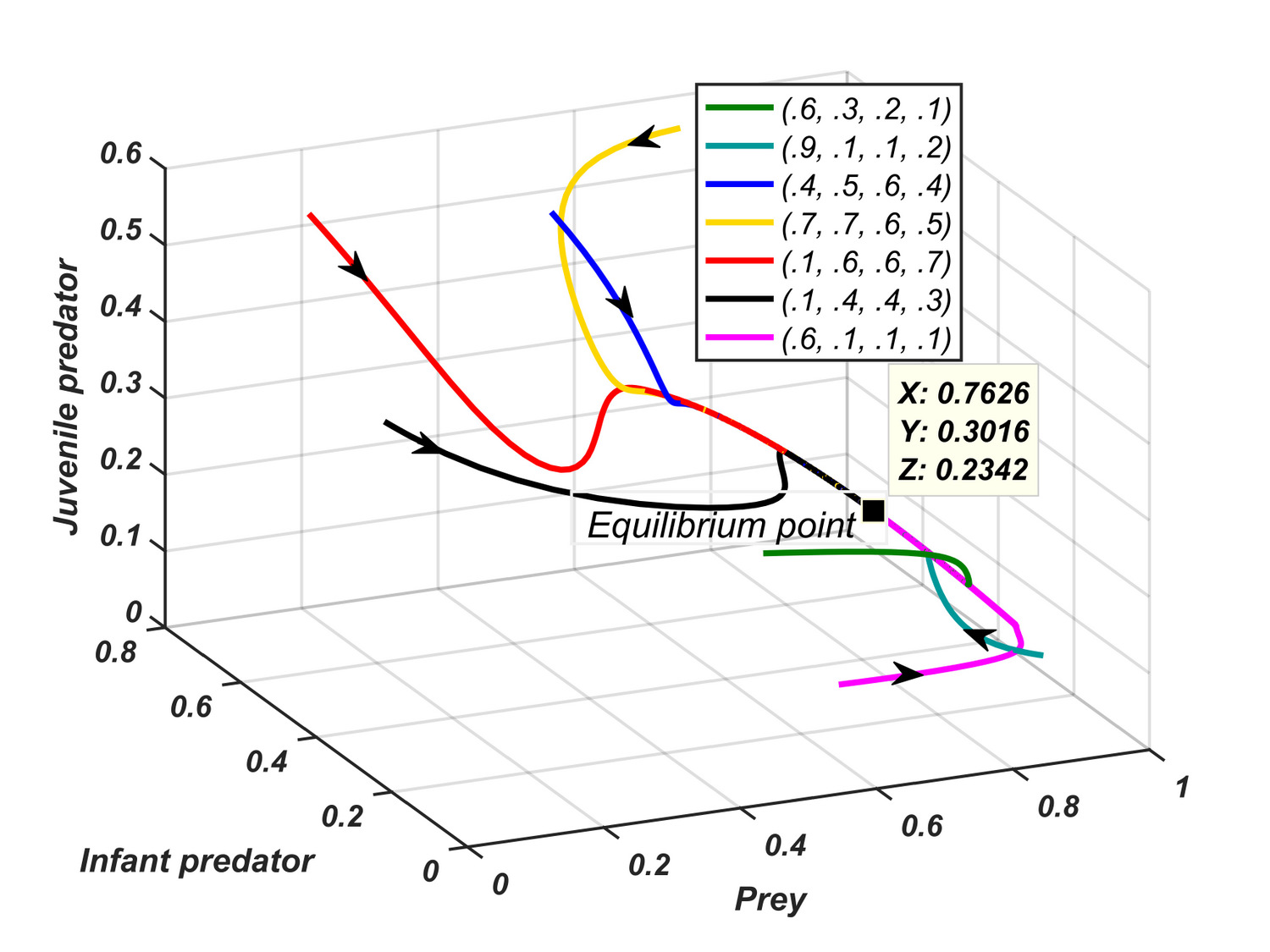}
    \caption{}
    \label{lo-xyz-1}
    \end{subfigure}\hfill
    \begin{subfigure}[b]{0.44\textwidth}
         \centering
    \includegraphics[width=\textwidth]{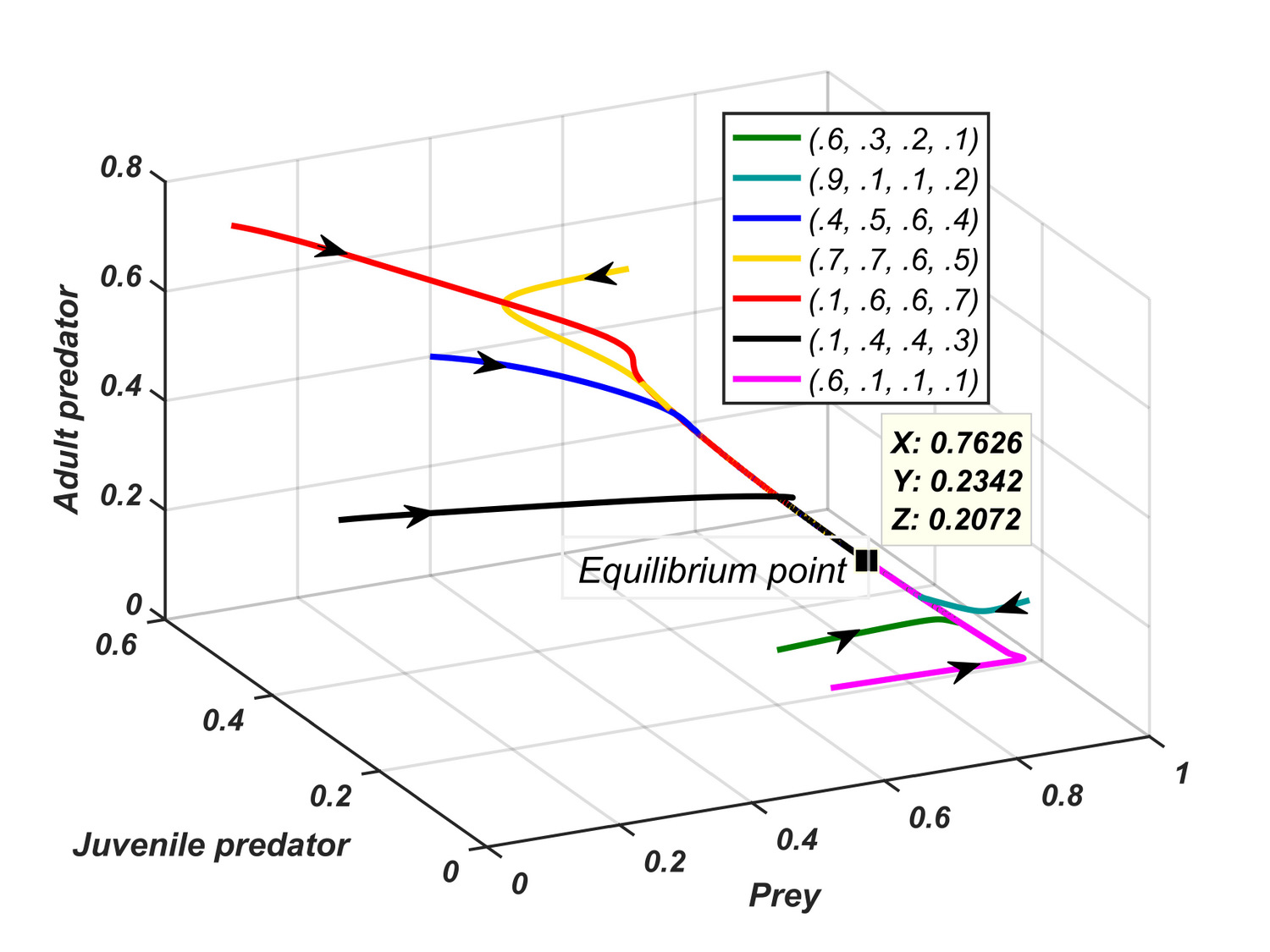}
    \caption{}
    \label{lo-xyz-2}
    \end{subfigure}
    \caption{Local and global stability of compresent equilibrium point: (a)depicts co-existence of prey, infant predators, and juvenile predators, (b)depicts co-existence of prey, juvenile predators, and adult predators. }
    \label{lo-xyz}
\end{figure}
\section{Bifurcations}
The focal point of this section is to discuss the existence of various transition points where the behavior of equations of bio-system (\ref{ma-eq}) modulates. The following transition points are investigated:
\subsection{Hopf bifurcation}
The transition point whose existence causes the system’s stability to switch and the birth or death of a periodic solution around an equilibrium point is the Hopf bifurcation point.\\ \\
First we strive towards the acquisition of the requisites for the manifestation of Hopf bifurcation around either of the compresent equilibrium points $E_3(x^*,y_1^*,y_2^*,y_3^*)$ and $E_4(x^{**},y_1^{**},y_2^{**},y_3^{**})$ with respect to any one of the parameters.\\
The characteristic equation of the variational matrix (\ref{jaco}) of the equations of bio-system (\ref{ma-eq}), as given in (\ref{ch-eq}) is
$$\omega^4+\epsilon_1 \omega^3+\epsilon_2 \omega^2+\epsilon_3 \omega+\epsilon_4=0. $$
If the point $c_h$ of parameter \textit{c} (transition rate of juvenile predator into adult) is presumed to be a bifurcation point, then we make the following assumptions:
$$(a)\; \; \epsilon_1>0,\epsilon_4>0,\epsilon_1 \epsilon_2-\epsilon_3>0, \hspace{0.5cm}
(b)\; \; \epsilon_1 \epsilon_2 \epsilon_3-\epsilon_3^2-\epsilon_1^2 \epsilon_4=0, \hspace{0.5cm}
(c)\; \;[\frac{d}{dc}(\epsilon_1 \epsilon_2 \epsilon_3-\epsilon_3^3-\epsilon_1^2 \epsilon_4)]_{c=c_h} \neq 0.$$
Taking the roots of the characteristic equation (\ref{jaco})  as: $-j,-k,-l,-m$ ; we would have
$$\epsilon_1=j+k+l+m,\hspace{0.7cm}
\epsilon_2=jk +lm+(j+k)(l+m),\hspace{0.7cm}
\epsilon_3=(j+k)lm+(l+m)jk, \hspace{0.7cm}
\epsilon_4=jklm.$$
$$\text{From (b): }\hspace{1cm}\epsilon_1 \epsilon_2 \epsilon_3-\epsilon_3^2-\epsilon_1^2\epsilon_4=0, \hspace{0.7cm}
i.e.,\; (j+k)(j+l)(j+m)(k+l)(k+m)(l+m)=0.$$
Therefore, at least one of the pair of roots should be additionally reciprocal to each other. Say $j+k=0$.  
$$\text{Hence,}\hspace{0.5cm}\epsilon_1 \epsilon_2-\epsilon_3=lm(l+m),\hspace{0.7cm}
\epsilon_1= l+m,\hspace{0.7cm}
\epsilon_4=-j^2lm. $$
Also, from $(a)$, we would have, $\hspace{0.4cm}lm(l+m)>0,\; \;l+m>0,\text{ and }-j^2lm>0$.\\
Therefore, no other pair of roots of the characteristic equation would be additionally reciprocal to each other, and $Re(l)>0,\; Re(m)>0$ and $j,k$ are the purely complex conjugate numbers, i.e., we would have two eigen values whose real part are negative and one pair of purely complex conjugates as the roots of the characteristic equation (\ref{ch-eq}). \\
Next, because of the property of continuity of roots, $\exists$ interval $(c_h-\varsigma,c_h+\varsigma)$ for some $\varsigma>0$, such that for $c$ belonging to this interval we would have,
$$j(c)=\psi_1(c)+\psi_2(c)i, \; \text{ and }\;
 k(c)=\psi_1(c)-\psi_2(c)i. $$
For the transversality condition, we need to prove $[\frac{d Re(j)}{dc}]_{c=c_h} \neq 0$ and $[\frac{d Re(k)}{dc}]_{c=c_h} \neq 0$.\\
For this we substitute $j(c)=\psi_1(c)+\psi_2(c)i$ in the characteristic equation (\ref{ch-eq}) and we get the following:
\begin{align}
\notag &(\psi_1+\psi_2i)^4+\epsilon_1 (\psi_1+\psi_2i)^3+\epsilon_2 (\psi_1+\psi_2i)^2+\epsilon_3(\psi_1+\psi_2i)+\epsilon_4=0\\
\notag  i.e.&\;
[4(\psi_1+\psi_2i)^3+3 \epsilon_1 (\psi_1+\psi_2i)^2+2 \epsilon_2 (\psi_1+\psi_2i)+\epsilon_3](\psi_1’+\psi_2’i)+[\epsilon_1’(\psi_1+\psi_2i)^3\\
\label{heq}&\hspace{0.5cm}+\epsilon_2'(\psi_1+\psi_2i)^2+ \epsilon_3’(\psi_1+\psi_2i)+\epsilon_4’]=0 \hspace{0.5cm}[\text{where }\epsilon'_i=\frac{d\epsilon_i}{dt}, \; \psi_i'=\frac{d\psi_i}{dt}].
\end{align}
Also using $ \epsilon_1 \epsilon_2 \epsilon_3-\epsilon_3^2-\epsilon_1^2\epsilon_4=0$ in the characteristic equation at $c=c_h$ we get,
$$\omega^4+ \epsilon_1\omega^3+(\frac{\epsilon_3}{\epsilon_1}+\frac{\epsilon_1 \epsilon_4}{\epsilon_3}) \omega^2+\epsilon_3 \omega+\epsilon_4=0, \hspace{0.5cm}
i.e.\; (\omega^2+\frac{\epsilon_3}{\epsilon_1})(\omega^2+\epsilon_1 \omega +\frac{\epsilon_1 \epsilon_4}{\epsilon_3})=0.$$
Utilizing $\psi_2= (\frac{\epsilon_3}{\epsilon_1})^{\frac{1}{2}} $ and $\psi_1=0$ at $c=c_h$ in \ref{heq}, we have
\begin{align*}
    &[-4i(\frac{\epsilon_3}{\epsilon_1})^{\frac{3}{2}}-3 \epsilon_1 (\frac{\epsilon_3}{\epsilon_1})+2 \epsilon_2 (\frac{\epsilon_3}{\epsilon_1})^{\frac{1}{2}}i+\epsilon_3](\psi_1’+\psi_2’i)+[-\epsilon_1’(\frac{\epsilon_3}{\epsilon_1})^{\frac{3}{2}}i- \epsilon_2’(\frac{\epsilon_3}{\epsilon_1})+\epsilon_3’(\frac{\epsilon_3}{\epsilon_1})^{\frac{1}{2}}i+\epsilon_4’]=0\\
&\text{Now, separating real and imaginary parts, we get the following:}\\
&\hspace{0.5cm} \text{\underline{Real part}}\\
&4 (\frac{\epsilon_3}{\epsilon_1})^{\frac{3}{2}} \psi_2’-3\epsilon_3 \psi_1’ -2 \epsilon_2 (\frac{\epsilon_3}{\epsilon_1})^{\frac{1}{2}} \psi_2’+\epsilon_3 \psi_1’-(\frac{\epsilon_3}{\epsilon_1})\epsilon_2’+\epsilon_4=0\\
\implies & \psi_2’=\frac{2 \epsilon_3 \psi_1’+ \epsilon_3 \epsilon_2’/\epsilon_1-\epsilon_4’}{2 (\frac{\epsilon_3}{\epsilon_1})^{\frac{1}{2}} (2 \frac{\epsilon_3}{\epsilon_1}-\epsilon_2)}.\\
\text{and}& \hspace{0.5cm} \text{\underline{Imaginary part}}\\
&-4(\frac{\epsilon_3}{\epsilon_1})^{\frac{3}{2}}\psi_1’ – 3 \epsilon_3 \psi_2’+2 \epsilon_2 (\frac{\epsilon_3}{\epsilon_1})^{\frac{1}{2}} \psi_1’+\epsilon_3 \psi_2’-\epsilon_1’ (\frac{\epsilon_3}{\epsilon_1})^{\frac{3}{2}}+(\frac{\epsilon_3}{\epsilon_1})^{\frac{1}{2}} \epsilon_3’=0\\
\implies& 2 (\epsilon_2-2 \frac{\epsilon_3}{\epsilon_1}) \psi_1 +\frac{-2 \epsilon_1 \epsilon_3 \psi_1’-\epsilon_3 \epsilon_2’+\epsilon_1 \epsilon_4’}{2 \frac{\epsilon_3}{\epsilon_1}-\epsilon_2}+(\epsilon_3’-\frac{\epsilon_1’ \epsilon_3}{\epsilon_1})=0 \hspace{10pt} \text{[putting the value of $\psi_2’$]}\\
\implies& \psi_1’=-\epsilon _1\frac{\epsilon_1’(\epsilon_2 \epsilon_3-2 \epsilon_1 \epsilon_4)+\epsilon_2’\epsilon_1 \epsilon_3+ \epsilon_3’(-2 \epsilon_3+\epsilon_1 \epsilon_2)-\epsilon_4’ \epsilon_1^2}{2[(\epsilon_3- \epsilon_2 \epsilon_2)^2+\epsilon_1^3 \epsilon_3]}.
\end{align*}
Therefore, transversality condition holds when
$$\epsilon_1’(\epsilon_2 \epsilon_3-2 \epsilon_1 \epsilon_4)+\epsilon_2’\epsilon_1 \epsilon_3+ \epsilon_3’(-2 \epsilon_3+\epsilon_1 \epsilon_2)-\epsilon_4’ \epsilon_1^2  \neq 0.$$
i.e. when our third assumption holds.
Hence, we arrive  at the following theorem:
\begin{theorem} \label{ho-thm}
The sufficient condition for the occurence of Hopf bifurcation in the given equations of bio-system (\ref{ma-eq}) with respect to any parameter, say $c$ is:
$$(a)\;\;\epsilon_1>0,\epsilon_4>0,\epsilon_1 \epsilon_2-\epsilon_3>0, \hspace{0.7cm}
(b)\;\; \epsilon_1 \epsilon_2 \epsilon_3-\epsilon_3^2-\epsilon_1^2 \epsilon_4=0,\text{ and}\hspace{0.5cm}
(c)\; \; [\frac{d}{dc}(\epsilon_1 \epsilon_2 \epsilon_3-\epsilon_3^3-\epsilon_1^2 \epsilon_4)]_{c=c_h} \neq 0.$$
\end{theorem}
\vspace{0.5cm}Now, we discuss the stability and orientation of the bifurcating periodic trajectory
\begin{theorem} \label{ho-dir-thm}
The orientation of the trajectory of Hopf bifurcation around a compresent equilibrium point can be ascertained by the sign of $\vartheta$, positivity would imply supercritical Hopf bifurcation while negativity would imply subcritical Hopf bifurcation; also positive $ \beta_2$ indicates unstability while negative $\beta_2$ implies stability. The values of $\vartheta,\beta_2$ are given within the proof.
\end{theorem}
\textbf{Proof: } So as to arbitrate the stability and orientation of bifurcating periodic solution occurring because of Hopf bifurcation, we use a procedure homologous to the methodology given by Hassard et al \cite{hassa}. 
At Hopf bifurcation point, we must have a pair of conjugate imaginary eigen values, say $\pm \alpha i$ and the other two eigen values are $v_1,$ and $v_2$.\vspace{0.1cm}\\
For the purpose of reducing the given equations of bio-system into normal form, new variables $(\bar{x},\bar{y_1},\bar{y_2},\bar{y_3})$ are introduced as:
$x=\bar{x}+X^*$, $y_1=\bar{y_1}+Y_1^*$, $y_2=\bar{y_2}+Y_2^*$, $y_3=\bar{y_3}+Y_3^*$.
For sake of simplicity while writing we take, $x \to \bar{x},y_1 \to \bar{y_1},y_2 \to \bar{y_2},y_3 \to \bar{y_3}$. Also $(X^*,Y_1^*,Y_2^*,Y_3^*)$ is used in lieu of either $E_3(x^*,y_1^*,y_2^*,y_3^*)$ or $E_4(x^{**},y_1^{**},y_2^{**},y_3^{**})$.\\
After the introduction of variables, the system of equations (\ref{ma-eq}) are:
\begin{align*}
    \dot{x}=&(x+X^*)(1-x-X^*)-a_1(x+X^*)(y_2+Y_2^*)-a2(x+X^*)(y_3+Y_3^*), \\
\dot{y_1}=& u a_2(x+X^*)(y_3+Y_3^*)-b(y_1+Y_1^*)-d_1(y_1+Y_1^*) ,\\
\dot{y_2}=&b(y_1+Y_1^*)-(c-a_3(x+X^*))(y_2+Y_2^*)- d_2(y_2+Y_2^*),\\
\text{and }\dot{y_3}=&(c-a_3(x+X^*))(y_2+Y_2^*)- d_3(y_3+Y_3^*).
\end{align*}
The above system of equations can be written as:
\begin{align} \label{tran}
    \dot{X}=UX+V,
\end{align}
where, $X=\left( \begin{array}{cccc} x& y_1& y_2& y_3\\ \end{array} \right)’, \; \;$ with  $U$ being the linear part and $V$ the non-linear part.\\
Next, a non-singular matrix $A= \left( \begin{array}{cccc}
 1 & 0 & 1 & 1 \\
 a_{21} & a_{22} & a_{23} & a_{24} \\
 a_{31} & a_{32} & a_{33} & a_{34} \\
 a_{41} & a_{42} & a_{43} & a_{44} \\
\end{array}
\right)$ is procured so as to get $A^{-1}UA= \left( \begin{array}{cccc}
0&\alpha &0 & 0\\
-\alpha& 0& 0& 0\\
0& 0& v_1&0\\
0&0&0&v_2\\
\end{array} \right).$\\
For this, we would need
\begin{align*}
    a_{21}=& [-u_{13} u_{24} \left(u_{24} u_{32} \left(\alpha ^2-u_{11} u_{33}\right)+u_{14} \left(u_{21} u_{32} u_{33}+u_{22} u_{31} u_{33}+\alpha ^2 (-u_{31})\right)\right)+u_{14} \left(\alpha ^2+u_{33}^2\right) (u_{24} \left(\alpha ^2+u_{11} u_{22}\right)\\
    &-u_{14} u_{21} u_{22})-u_{13}^2 u_{24}^2 u_{31} u_{32}]/[u_{13}^2 u_{24}^2 u_{32}^2+2 u_{13} u_{14} u_{24} u_{32} \left(u_{22} u_{33}-\alpha ^2\right)+u_{14}^2 \left(\alpha ^2+u_{22}^2\right) \left(\alpha ^2+u_{33}^2\right)],\\
    a_{22}=& -\frac{\alpha  \left(u_{13} u_{24} (u_{24} u_{32} (u_{11}+u_{33})+u_{14} (-u_{21} u_{32}+u_{22} u_{31}+u_{31} u_{33}))+u_{14} \left(\alpha ^2+u_{33}^2\right) (u_{24} (u_{22}-u_{11})+u_{14} u_{21})\right)}{u_{13}^2 u_{24}^2 u_{32}^2+2 u_{13} u_{14} u_{24} u_{32} \left(u_{22} u_{33}-\alpha ^2\right)+u_{14}^2 \left(\alpha ^2+u_{22}^2\right) \left(\alpha ^2+u_{33}^2\right)},\\
    a_{23}=& -\frac{(u_{33}-v_{1}) (u_{24} (v_{1}-u_{11})+u_{14} u_{21})+u_{13} u_{24} u_{31}}{u_{13} u_{24} u_{32}+u_{14} (u_{22}-v_1) (u_{33}-v_1)},\\
    a_{24}=&-\frac{(u_{33}-v_{2}) (u_{24} (v_{2}-u_{11})+u_{14} u_{21})+u_{13} u_{24} u_{31}}{u_{13} u_{24} u_{32}+u_{14} (u_{22}-v_{2}) (u_{33}-v_{2})},\\
    a_{31}=& -[u_{13} u_{24} u_{32} (u_{11} u_{24} u_{32}-u_{14} u_{21} u_{32}+u_{14} u_{22} u_{31})+u_{14} (u_{24} u_{32} \left(u_{11} \left(u_{22} u_{33}-\alpha ^2\right)+\alpha ^2 (u_{22}+u_{33})\right)
    +u_{14} (-u_{21} u_{22} u_{32} u_{33}\\
    &+\alpha ^2 (u_{21} u_{32}+u_{31} u_{33})+u_{22}^2 u_{31} u_{33}))]/[u_{13}^2 u_{24}^2 u_{32}^2+2 u_{13} u_{14} u_{24} u_{32} \left(u_{22} u_{33}-\alpha ^2\right) +u_{14}^2 \left(\alpha ^2+u_{22}^2\right) \left(\alpha ^2+u_{33}^2\right)],\\
    a_{32}=& -\frac{\alpha  \left(u_{14} u_{24} u_{32} \left(\alpha ^2+u_{11} (u_{22}+u_{33})-u_{13} u_{31}-u_{22} u_{33}\right)-u_{13} u_{24}^2 u_{32}^2+u_{14}^2 \left(-u_{21} u_{22} u_{32}-u_{21} u_{32} u_{33}+u_{22}^2 u_{31}+\alpha ^2 u_{31}\right)\right)}{u_{13}^2 u_{24}^2 u_{32}^2+2 u_{13} u_{14} u_{24} u_{32} \left(u_{22} u_{33}-\alpha ^2\right)+u_{14}^2 \left(\alpha ^2+u_{22}^2\right) \left(\alpha ^2+u_{33}^2\right)},\\
    a_{33}=& \frac{u_{24} u_{32} (v_{1}-u_{11})+u_{14} (u_{21} u_{32}-u_{22} u_{31}+u_{31} v_{1})}{u_{13} u_{24} u_{32}+u_{14} (u_{22}-v_{1}) (u_{33}-v_{1})},\\
    a_{34}=& \frac{u_{24} u_{32} (v_{2}-u_{11})+u_{14} (u_{21} u_{32}-u_{22} u_{31}+u_{31} v_{2})}{u_{13} u_{24} u_{32}+u_{14} (u_{22}-v_{2}) (u_{33}-v_{2})},\\
    a_{41}=& [u_{13} \left(u_{24} u_{32} \left(u_{11} \left(\alpha ^2-u_{22} u_{33}\right)+\alpha ^2 (u_{22}+u_{33})\right)+u_{14} \left(-u_{21} u_{22} u_{32} u_{33}+\alpha ^2 (u_{21} u_{32}+u_{31} u_{33})+u_{22}^2 u_{31} u_{33}\right)\right)
    -u_{11} u_{14} (\alpha ^2\\   &+u_{22}^2) \left(\alpha ^2+u_{33}^2\right)+u_{13}^2 u_{24} u_{32} (u_{22} u_{31}-u_{21} u_{32})]/[u_{13}^2 u_{24}^2 u_{32}^2+2 u_{13} u_{14} u_{24} u_{32} \left(u_{22} u_{33}-\alpha ^2\right)
    +u_{14}^2 \left(\alpha ^2+u_{22}^2\right) \left(\alpha ^2+u_{33}^2\right)], \\
   a_{42}=& [\alpha  (u_{13} \left(u_{24} u_{32} \left(-\alpha ^2+u_{11} (u_{22}+u_{33})+u_{22} u_{33}\right)+u_{14} \left(-u_{21} u_{22} u_{32}-u_{21} u_{32} u_{33}+u_{22}^2 u_{31}+\alpha ^2 u_{31}\right)\right)-u_{13}^2 u_{24} u_{31} u_{32}\\
   &+u_{14} \left(\alpha ^2+u_{22}^2\right) \left(\alpha ^2+u_{33}^2\right))]/[u_{13}^2 u_{24}^2 u_{32}^2+2 u_{13} u_{14} u_{24} u_{32} \left(u_{22} u_{33}-\alpha ^2\right)+u_{14}^2 \left(\alpha ^2+u_{22}^2\right) \left(\alpha ^2+u_{33}^2\right)] ,\\
   a_{43}=& \frac{u_{13} (-u_{21} u_{32}+u_{22} u_{31}-u_{31} v_{1})-(u_{11}-v_{1}) (v_{1}-u_{22}) (v_{1}-u_{33})}{u_{13} u_{24} u_{32}+u_{14} (u_{22}-v_{1}) (u_{33}-v_{1})}, \\
   a_{44}=& \frac{u_{13} (-u_{21} u_{32}+u_{22} u_{31}-u_{31} v_{2})-(u_{11}-v_{2}) (v_{2}-u_{22}) (v_{2}-u_{33})}{u_{13} u_{24} u_{32}+u_{14} (u_{22}-v_{2}) (u_{33}-v_{2})}. \\
\end{align*}
$u_{ij}$ $(i,j=1,2,3,4)$ are the elements of the matrix U, i.e.$ U=[u_{ij}]_{4\times 4}$.\vspace{0.14cm}\\ 
To obtain the normal form, another transformation $X=AZ$ is done, where $Z=\left(\begin{array}{cccc} z_1& z_2 & z_3& z_4\\ \end{array} \right)'$. Therefore, (\ref{tran}) becomes  $\dot{X}=UX+V \; \; \; \;
\implies \dot{Z}=A^{-1} U (AZ)+A^{-1} V \; \; \; \;
\implies \dot{Z}= A^{-1} U AZ+ W $.\vspace{0.14cm}\\
The matrix $W=\left( \begin{array}{cccc} w_1& w_2& w_3& w_4\\ \end{array} \right)’$ is obtainable by transformation of V done by the following\\
$x=z_1+z_3+z_4,\;\;$ 
$y_1=a_{21} z_1+a_{22} z_2+a_{23} z_3+a_{24} z_4,\; \;$ 
$y_2=a_{31} z_1+a_{32} z_2+a_{33} z_3+a_{34} z_4,$ and 
$y_3=a_{41} z_1+a_{42} z_2+a_{43} z_3+a_{44} z_4.$\\
Now, we arrive at the stage where we can derive the stability and orientation of the bifurcating periodic solution with the help of the expressions, $g_{11},g_{02},g_{20},G_{21}, h^1_{11},h^1_{20},w_{11},w_{20},G_{110}^1,G_{101}^1,g_{21}$, all of which can be obtained as is given by Hassard et al \cite{hassa} at $(z_1,z_2,z_3,z_4)=(0,0,0,0)$. Finally, we arrive at the required expression for determining the direction of Hopf bifurcation,
$$C_1(0)=\frac{i}{2 \alpha} \left(g_{20} g_{11} -2 |g_{11}|^2-\frac{1}{3} |g_{02}|^2\right)+\frac{g_{21}}{2}.$$
$$ \therefore \; \;\vartheta=-\frac{Re C_1(0)}{Re \omega’(c_h)}, \text{ and}\hspace{0.2cm} \beta_2=2 ReC_1(0).$$
$\vartheta>0$ implies supercritical Hopf bifurcation and $\vartheta<0$ implies subcritical Hopf bifurcation. Moreover, 
$\beta_2<0$ and $\beta_2>0$ would mean stability and unstablity for the bifurcating periodic solution, respectively.
\begin{figure}[H]
    \centering
    \begin{subfigure}[b]{0.44\textwidth}
         \centering
    \includegraphics[width=\textwidth]{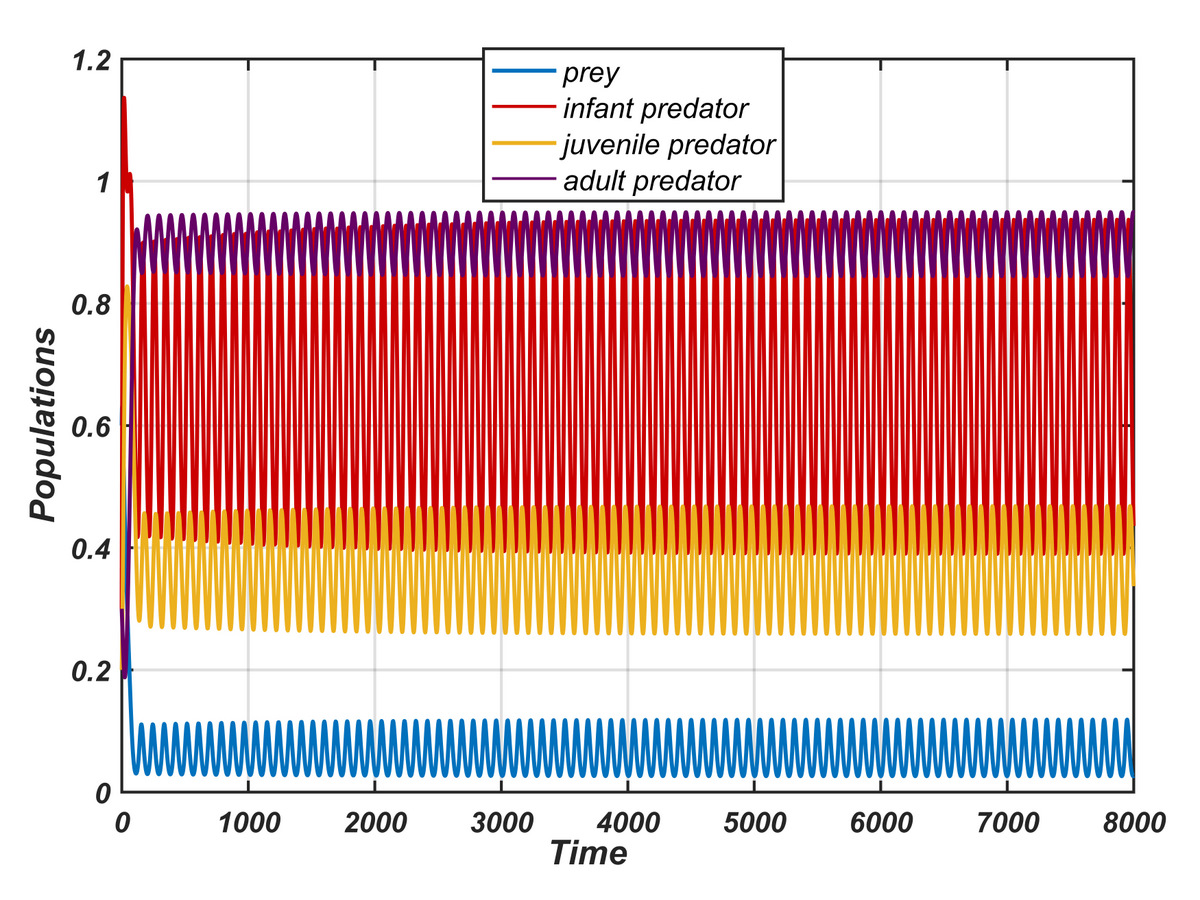}
    \caption{at $c=0.037$}
    \label{c-1}
    \end{subfigure}\hfill
    \begin{subfigure}[b]{0.44\textwidth}
         \centering
    \includegraphics[width=\textwidth]{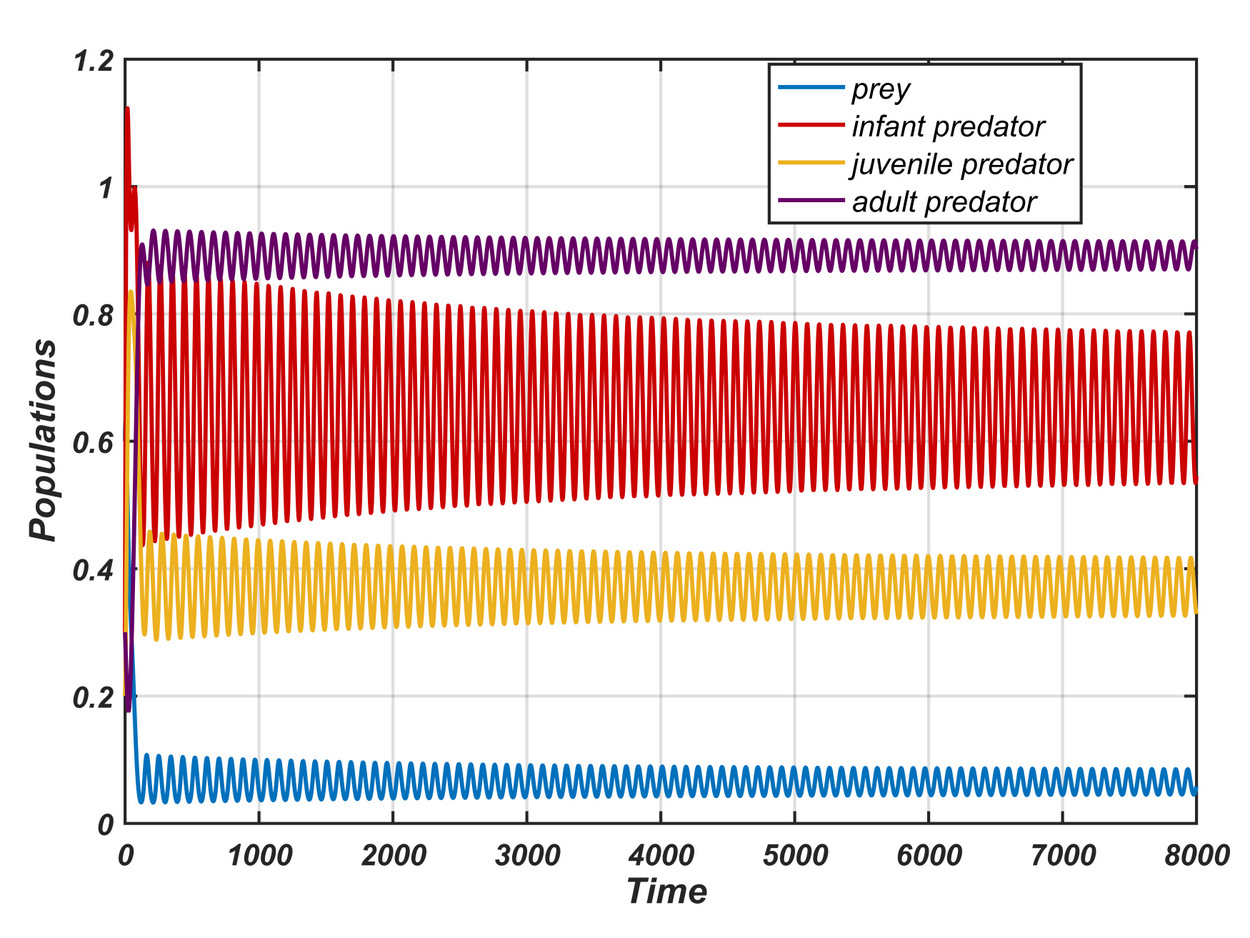}
    \caption{at $c=0.03598345$}
    \label{c-2}
    \end{subfigure}\hfill
     \begin{subfigure}[b]{0.44\textwidth}
         \centering
    \includegraphics[width=\textwidth]{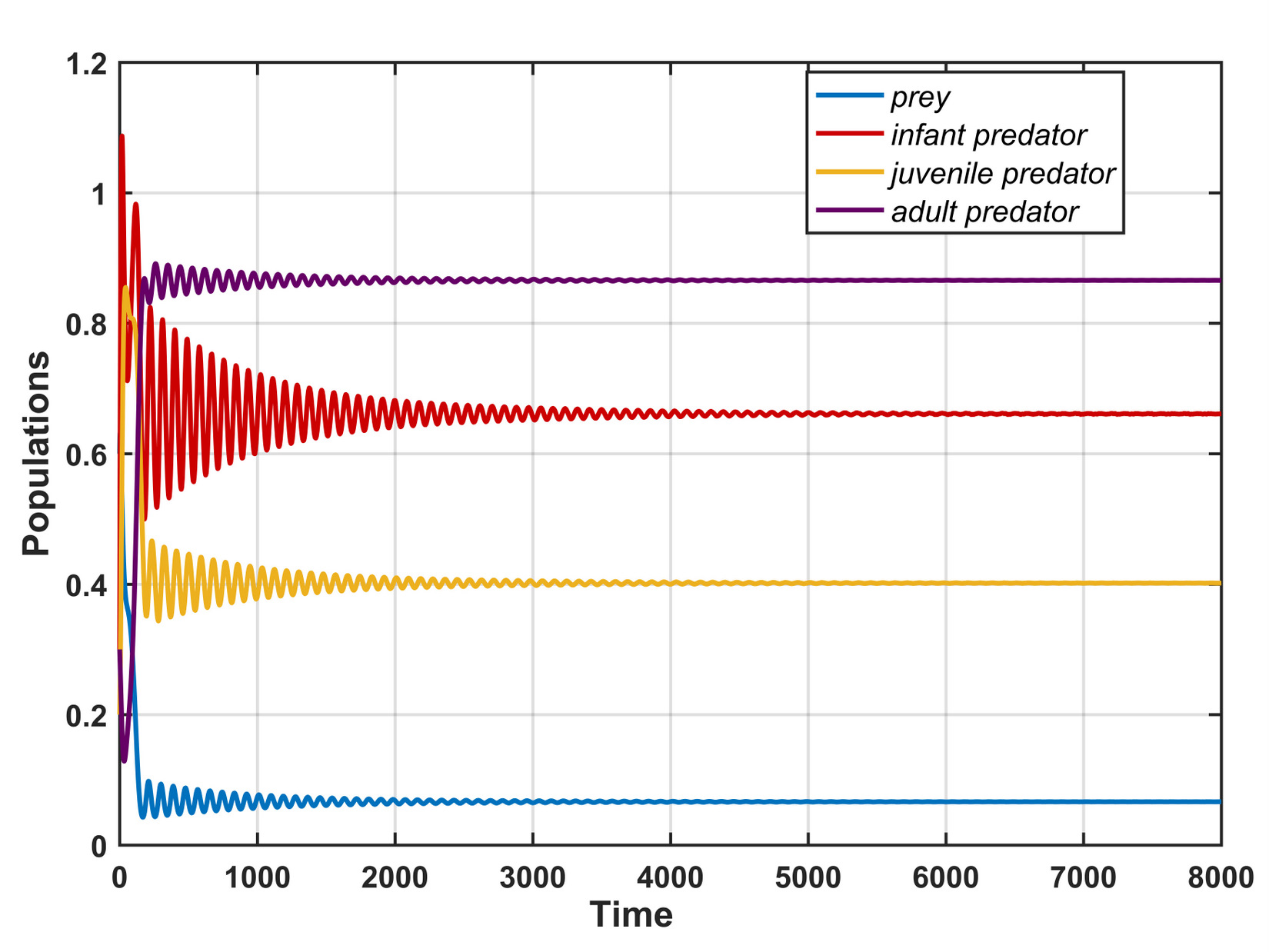}
    \caption{at $c=0.033$}
    \label{c-3}
    \end{subfigure}
    \caption{Time series depicting equilibrium states in the neighbourhood of Hopf bifurcation point of $c$, the transition rate of juvenile predators into adult predators: $(a)c>c_h$, $(b)c=c_h$, $(c)c<c_h$, $c_h$ is the Hopf bifurcation point.}
    \label{c-ho}
\end{figure}
\subsection{Transcritical bifurcation}
The transition point where two equilibria convene and the system experiences changes in stability in regards to its equilibrium points is the transcritical bifurcation point.   
\begin{theorem} \label{thm-tb} The equations of bio-system (\ref{ma-eq}) experiences interchangeability of stability between the compresent equilibrium point and prey-only equilibrium when $a_2= [(b+d¬_1) (a_3-c-d_2) d_3]/[b (a_3-c) u]$. \end{theorem}
\textbf{Proof:} The variational matrix of the equations of bio-system (\ref{ma-eq}) at $E_2(1,0,0,0)$ is: 
\begin{equation}
    \label{tra} \left(
\begin{array}{cccc}
 -1 & 0 & -a_1 & -a_2 \\
 0 & -b-d_1 & 0 & a_2 u \\
 0 & b & a_3-c-d_2 & 0 \\
 0 & 0 & c-a_3 & -d_3 \\
\end{array}
\right)
\end{equation}
For $ a_2= [(b+d_1) (a_3-c-d_2) d_3]/[b (a_3-c) u](=a_{2t})$, the determinant of the above matrix would be zero, i.e., zero would be the value of one of the eigen values.
Corresponding to this zero eigen value for the matrices (\ref{tra}) and its transpose at $a_{2t}$, the respective eigen vectors are:\\
$l=\left(-a_1 l_3 + \frac{(b + d_1) (a_3 - c – d_2) l_3}{b u}, \frac{(-a_3+c+d_2) l_3}{b}, l_3, \frac{(-a_3 + c) l_3}{d3} \right)$, and\\
$m=\left(0,m_2, ((b + d_1) m_2)/b), \frac{(b + d_1) (a_3 - c – d_2) m_2}{b (a_3 – c)} \right)$ respectively, where $l_3$ and $m_2$ are arbitrary\\
Next, the given model (\ref{ma-eq}) is rewritten in the form 
$$ \dot{Y}=BY,$$
where $Y=\left( \begin{array}{c} x\\ y_1\\ y_2\\ y_3\\ \end{array} \right) $, and $ B=\left(
\begin{array}{c}
 (1-x) x-a_1 y_2 x-a_2 y_3 x \\
 -b y_1-d_1y_1+a_2 u x y_3 \\
 b y_1-d_2 y_2-(c-a_3 x) y_2 \\
 (c-a_3 x) y_2-d_3 y_3 \\
\end{array}
\right)$. \\
According to the Sotomayor theorem \cite{soto}, the instance of transcritical bifurcation can be proven by\\ \\
$m B_{a_2}(E_2,a_{2t})=0$,\hspace{15pt} $m(D B_{a_2}(E_2,a_{2t}))\hspace{3pt}l^t \neq 0$,\hspace{15pt} $m[D^2(B_{a_2}(E_2,a_{2t}))]\hspace{3pt}(l^t,l^t) \neq 0$.\\ \\
In our case, we have
\begin{align*} &B_{a_2}(E_2,a_{2t})=\left( \begin{array}{cccc} \frac{\partial B_1}{\partial a_2} & \frac{\partial B_2}{\partial a_2} & \frac{\partial B_3}{\partial a_2} & \frac{\partial B_4}{\partial a_2} \\ \end{array} \right)_{(1,0,0,0);a_{2t}}\\
& \text{Thus, } \; m B_{a_2}(E_2,a_{2t})=0,\\
&D B_{a_2}(E_2,a_{2t})= \left( \begin{array}{cccc} \frac{\partial^2 B_1}{\partial x\, \partial a_2} & \frac{\partial^2 B_1}{\partial y_1\, \partial a_2} & \frac{\partial^2 B_1}{\partial y_2\, \partial a_2} & \frac{\partial^2 B_1}{\partial y_3\, \partial a_2} \\ 
\frac{\partial^2 B_2}{\partial x\, \partial a_2} & \frac{\partial^2 B_2}{\partial y_1\, \partial a_2} & \frac{\partial^2 B_2}{\partial y_2\, \partial a_2} & \frac{\partial^2 B_2}{\partial y_3\, \partial a_2}\\
\frac{\partial^2 B_3}{\partial x\, \partial a_2} & \frac{\partial^2 B_3}{\partial y_1\, \partial a_2} & \frac{\partial^2 B_3}{\partial y_2\, \partial a_2} & \frac{\partial^2 B_3}{\partial y_3\, \partial a_2}\\ 
\frac{\partial^2 B_4}{\partial x\, \partial a_2} & \frac{\partial^2 B_4}{\partial y_1\, \partial a_2} & \frac{\partial^2 B_4}{\partial y_2\, \partial a_2} & \frac{\partial^2 B_4}{\partial y_3\, \partial a_2}\\ \end{array} \right)_{(1,0,0,0);a_{2t}}\\
& \text{Thus, } \; 
m(D B_{a_2}(E_2,a_{2t}))l^t= m_2 u \frac{(-a_3 + c) l_3}{d3} \neq 0,\\
&D^2B(E_2,a_{2t})(l^t,l^t)=\left( \begin{array}{c} \sum_{i=1}^4\sum _{j=1}^4[\frac{\partial^2 B_1}{\partial p_i\, \partial p_j} l_{i} l_{j}] \\
\sum _{i=1}^4\sum _{j=1}^4[\frac{\partial^2 B_2}{\partial p_i\, \partial p_j} l_{i} l_{j}]\\
\sum _{i=1}^4\sum _{j=1}^4[\frac{\partial^2 B_3}{\partial p_i\, \partial p_j} l_{i} l_{j}]\\
\sum _{i=1}^4\sum _{j=1}^4[\frac{\partial^2 B_4}{\partial p_i\, \partial p_j} l_{i} l_{j}]\\
\end{array} \right)'_{(1,0,0,0);\;\text{{\small $p_i=(x,y_1,y_2,y_3)$}}}\\
\text{Thus, } \; &m D^2B_{a_2}(E_2,a_{2t})(l^t,l^t)= - \frac{2 l_3^2 m_2 ( (b + d_1) (a_3-c- d_2) - a_1 b u) (-a_3 (b + d_1) d_2 d_3 + a_2 b (a_3 - c)^2 u)}{ b^2 (a_3 - c) d_3 u}\; \;  \neq 0.
\end{align*}
Hence, the theorem is proved.\\ \\
\textbf{Note:} From direct observation, it is to be noted that for the occurance of trancritical bifurcation, we must have $a_3<c$ or $c+d_2<a_3$. Also in lieu of $a_2$ other parametric values can be used as bifurcating parameter.
\subsection{Saddle node bifurcation}
The point of collision and disappearance of two branches of equilibria is the saddle node bifurcation point (also known as limit point bifurcation). 
\begin{theorem} \label{thm-snb} The equations of bio-system (\ref{ma-eq}) experiences saddle node bifurcation when the parameters satisfy the expression $a_2bux(a_3x(-1+2x)-c(-1+2x+a_1y_2)+a_1d_3y_3)+(-b-d_1)[a_2a_3d_2xy_2-d_3(a_1a_3xy_2+(c+d_2-a_3x)(-1+2x+a_1y_2+a_2y_3))]=0$ around the equilibrium points, along with $\chi_2^T [B_b(E_4,b_s)] \neq 0, \text{ and } \chi_2^T[D^2B_b(E_4,b_s)(\chi_1,\chi_1)] \neq 0.$
\end{theorem}
\textbf{Proof:} For saddle node bifurcation to occur at an equilibrium point, one of the eigen values of the variational matrix (\ref{jaco}) of our model needs to be zero, i.e., the determinant of (\ref{jaco}) would be zero. Hence, we would need-
$$ a_2bux(a_3x(-1+2x)-c(-1+2x+a_1y_2)+a_1d_3y_3)+(-b-d_1)[a_2a_3d_2xy_2-d_3(a_1a_3xy_2+(c+d_2-a_3x)(-1+2x+a_1y_2+a_2y_3))]=0,$$
where, $x,y_1,y_2,y_3$ are to be replaced by the points of equilibrium where the bifurcation is supposed to occur. Any one parameter from the above expression can be used as the control parameter, say $b=b_s$ is the point of saddle node bifurcation.\\
Suppose $E_4$ is the equilibrium point, then we have
\begin{align*}
    &\frac{1}{2 a_2^2 a_3^2 b^2 u^2}(a_3^3 (b + d_1)^2 d_3^2 ((b + d_1) d_3 – a_2 b u) + a_2^2 b^2 c^2 u^2 (-R + a_2 b c u) + a_2 a_3 b u (2 (b + d_1) d_2 d_3 R + a_2 b c (-(b + d_1) (c + 4 d_2) d_3\\
    &+ R) u – a_2^2 b^2 c^2 u^2) – a_3^2 (b + d_1) d_3 ((b + d_1) d_3 R + a_2 b ((b + d_1) (c + 4 d_2) d_3 - R) u - 2 a_2^2 b^2 (c + 2 d_2) u^2))=0,\\
    & \hspace{1cm} \text{where, } R=\sqrt{a_2^2 b^2 c^2 u^2-2 a_2 a_3 b d_3 u (b+d_1) (c+2 d_2)+a_3^2 d_3^2 (b+d_1)^2}.
\end{align*}
To simplify the calculations, the matrix (\ref{jaco}) at compresent equilibrium point is be rewritten as 
\begin{align}
    \left( \begin{array}{cccc}
 e_1 & 0 & e_2 & e_3 \\
 e_4 & e_5 & 0 & e_6 \\
e_7 & e_8 & e_9 & 0 \\
 e_{10} & 0 & e_{11} & e_{12} \\ 
 \end{array} \right).
\label{sad-jaco}\end{align}
Hence, the determinant of the above matrix would be zero, when $e_5=\frac{ (e_{12} e_2 e_4 – e_{11} e_3 e_4 + e_1 e_{11} e_6 – e_{10} e_2 e_6) e_8}{ e_{12} e_2 e_7 – e_{11} e_3 e_7 – e_1 e_{12} e_9 + e_{10} e_3 e_9}$.\\
Therefore, the eigen vectors corresponding to zero eigen value of matrix (\ref{sad-jaco}) with the value of $e_5=\frac{ (e_{12} e_2 e_4 – e_{11} e_3 e_4 + e_1 e_{11} e_6 – e_{10} e_2 e_6) e_8}{ e_{12} e_2 e_7 – e_{11} e_3 e_7 – e_1 e_{12} e_9 + e_{10} e_3 e_9}$  and its transpose are:
\begin{align*}
    \chi_1= &\left( \phi, - \frac{ (e_{12} e_2 e_7 – e_{11} e_3 e_7 – e_1 e_{12} e_9 + e_{10} e_3 e_9) \phi}{ (e_{12} e_2 – e_{11} e_3) e_8} , -\frac{(e_1 e_{12} – e_{10} e_3) \phi}{e_{12} e_2 – e_{11} e_3}, - \frac{(e_1 e_{11} – e_{10} e_2) \phi}{-e_{12} e_2 + e_{11} e_3} \right)^t,\\
 \text{ and } \chi_2= & \left( \frac{-(e_{11} e_6 e_7 + e_{12} e_4 e_9 – e_{10} e_6 e_9) \phi_1}{-e_{12} e_2 e_7 + e_{11} e_3 e_7 + 
 e_1 e_{12} e_9 – e_{10} e_3 e_9}, \phi_1, \frac{- (e_{12} e_2 e_4 – e_{11} e_3 e_4 + e_1 e_{11} e_6 – e_{10} e_2 e_6) \phi_1}{e_{12} e_2 e_7 – e_{11} e_3 e_7 – e_1 e_{12} e_9 + e_{10} e_3 e_9}, \frac{- (e_2 e_6 e_7 + e_3 e_4 e_9 – e_1 e_6 e_9) \phi_1}{e_{12} e_2 e_7 – e_{11} e_3 e_7 – e_1 e_{12} e_{9} + e_{10} e_3 e_9}\right)^t,\\
\text{where } \phi, \phi_1& \text{ are arbitrary real numbers}.
\end{align*}
According to Sotomayor theorem, saddle node bifurcation would occur when
$\chi_2^T [B_b(E_4,b_s)] \neq 0, \text{ and } \chi_2^T[D^2B_b(E_4,b_s)(\chi_1,\chi_1)] \neq 0.$
This can be proven in a similar manner as done in Theorem \ref{thm-tb}. This is verified numerically in section \ref{numeri} for figure \ref{b-equ}.
\section{Numerical simulations}
\label{numeri}
In this section, the aim is to validate the analytical findings, and also to explicate the dynamics of the equations of bio-system using ode45, Matcont package present in Matlab software \cite{matc}. For the illustrations, biologically feasible, and hypothetical values of the parameters are considered as given in tables 1 and 2. The values may vary within a suitable range.
\begin{multicols}{2}
\begin{table}[H]
    \centering
    \begin{tabular}{|c|c|c|}
    \hline
   Biological interpretation& Parameter&Value\\
   \hline
       Predation rate by juvenile predators& $a_1$ &0.46 \\
        \hline
        Predation rate by matured predators& $a_2$&0.625\\
         \hline
         Transformation rate of consumed prey into infant predator &$u$& 0.8\\
         \hline
         Transition rate of infant predators into juvenile stage &$b$& 0.112\\
         \hline
         Transition rate of juvenile predators into adult stage &$c$&0.09\\
         \hline
        Rate of maturation delay in juvenile predators& $a_3$&0.06\\
         \hline
                Natural death rate of infant predators& $d_1$&0.15\\
         \hline
         Natural death rate of juvenile predators &$d_2$&0.1\\
         \hline
         Natural death rate of adult predators &$d_3$&0.05\\
         \hline
    \end{tabular}
    \caption{Parametric values.}
    \label{table1}
\end{table}

\begin{table}[H]
    \centering
    \begin{tabular}{|c|c|c|}
    \hline
    Parameter&Value\\
   \hline
        $a_1$ &0.6 \\
        \hline
        $a_2$&0.8\\
         \hline
        $u$& 0.82\\
         \hline
        $b$& 0.031\\
         \hline
        $c$&0.035\\
         \hline
        $a_3$&0.075\\
         \hline
          $d_1$&0.026\\
         \hline
        $d_2$&0.023\\
         \hline
        $d_3$&0.013\\
         \hline
    \end{tabular}
    \caption{Parametric values.}
    \label{table2}
\end{table}
\end{multicols}
\par For the data set as given in table 1 except for $u$, which is taken $u=0.7$, it is observed that the equations of bio-system admits a locally asymptotically stable prey-only equilibrium $E_2(1,0,0,0)$.  Analytically, as given in the proof of Theorem \ref{thm-lo-pf}, either one of the three cases requisite for the stability of $E_2(1,0,0,0)$ should be satisfied. And since we have $(a_3=0.06)<(c=0.09)$ and $(a_2=0.625)<( \frac{(b+d_1)(a_3-c-d_2)d_3}{bu(a_3-c)}=0.6336) $, the equilibrium point is locally stable, which is illustrated in figure \ref{lo-axi}. Also, the trajectories of widely varying co-existing initial populations can be seen leading to the prey-only equilibrium, even if not all of them satisfy the global  stability condition \ref{thm-glo-pf}, because of its sufficient nature.
\par For the parametric values from table 1, both of the compresent equilibrium points exist, but only one of them is espied to be stable. For the equilibrium point $E_3(0.97,0.04,0.03,0.02)$, from Theorem \ref{sta-co} we have $\epsilon_1=1.41, \, \epsilon_4=-0.000019, \, \epsilon_1 \epsilon_2-\epsilon_3=0.6329, \, \epsilon_1 \epsilon_2 \epsilon_3 - \epsilon_3^2- \epsilon_1^2 \epsilon_4=0.0335$; and for equilibrium point $E_4(0.762, 0.30159, 0.23419, 0.207)$, $\epsilon_1=1.219, \, \epsilon_4=0.000126, \, \epsilon_1 \epsilon_2-\epsilon_3=0. 44539, \, \epsilon_1 \epsilon_2 \epsilon_3 - \epsilon_3^2- \epsilon_1^2 \epsilon_4=0. 02096$. Hence, the system is locally asymptotically stable only around the equilibrium point $E_4$, due to satisfying the Routh-Hurwitz criteria.   Figure \ref{lo-xyz} illustrates this scenario, also at the same time, the widely varying initial populations approaching equilibrium $E_4$ depicts the global stability of the system.\\
\par
\begin{figure}[H]
    \centering
    \includegraphics[width=10cm]{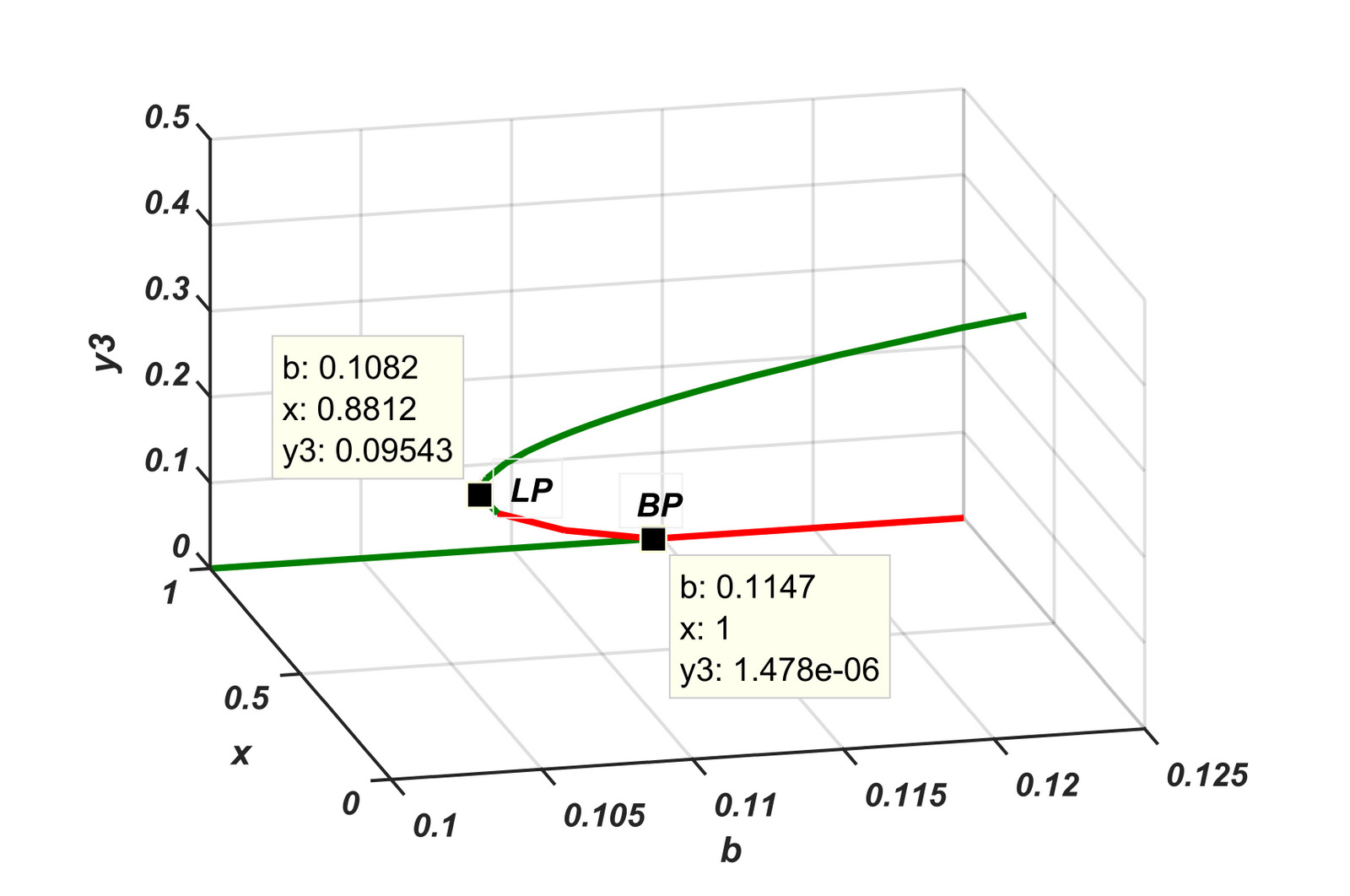}
    \caption{Equilibrium curve in context of $b$, the transition rate of infant predators into the juvenile stage. Here, \textbf{$LP$} denotes the saddle node bifurcation point, and \textbf{$BP$} denotes the transcritical bifurcation point.}
    \label{b-equ}
\end{figure}

  To investigate the bifurcation behavior of the system, we begin by varying the transition rate of infant predators into the juvenile stage $b$ (the rest of the parametric values are taken from table 1). Graphically, a saddle node bifurcation (also known as fold bifurcation) is espied at $b=0.108186$, the point where a stable compresent equilibrium $E_4$ and an unstable compresent equilibrium $E_3$ encounter and annihilate each other. In addition, transcritical bifurcation is espied at $b=0.114706$, where the unstable compresent equilibrium encounters the prey-only equilibrium. Figure \ref{b-equ}  illustrates all this, where the red line depicts unstable equilibrium and green line depicts stable equilibrium. Numerically, for $b=0.114706$, we get $\frac{(b+d_1)(a_3- c-d_2)d_3}{b(a_3- c)u}=0.625\; (=a_2)$, which is the condition of transcritical bifurcation of Theorem \ref{thm-tb}, and for $b=0.108186$, we have $a_2b u x(a_3 x(-1 + 2x) -c(-1+ 2x + a_1 y_2) +a_1d_3y_3)+ (-b -d_1)[a_2 a_3 d_2 x y_2 - d_3(a_1 a_3 x y_2 + (c + d_2- a_3x)(-1 +2x +a_1y_2+ a_2y_3))] \approx 0$, which is the condition of saddle node bifurcation of Theorem \ref{thm-snb} along with $\chi_2^T[B_b(E_4,b_s)]=0.02745 \neq 0$ and $\chi_2^T[D^2B_b(E_4,b_s)(\chi_1,\chi_2)]=-0.00421 \neq 0$. On that account, it can be said that the theoretical part agrees with the numerical simulations.\\
\par
\begin{figure}[H]
   \begin{subfigure}[b]{0.46\textwidth}
         \centering
    \includegraphics[width=\textwidth]{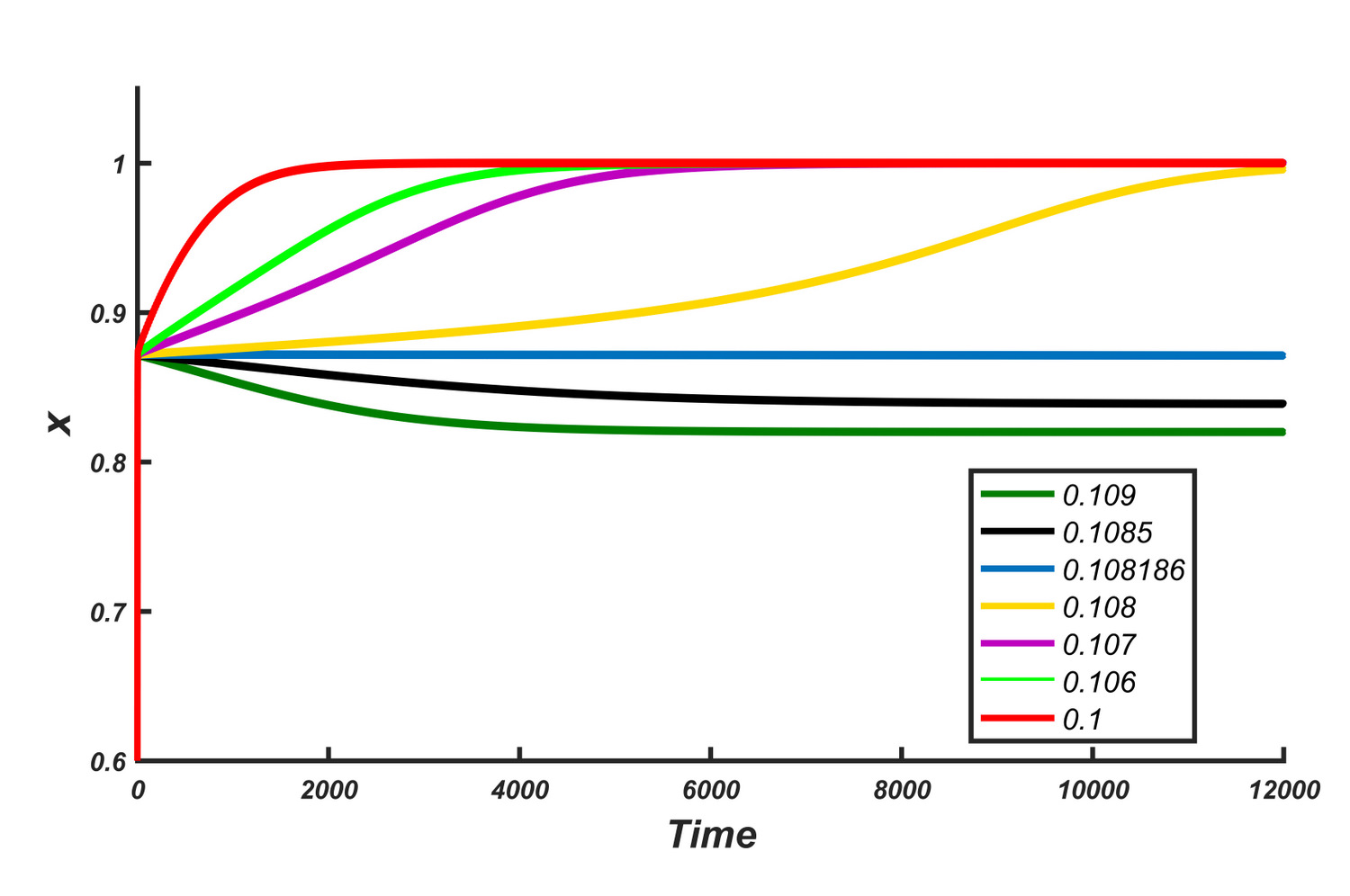}
    \caption{}
    \label{bx}
    \end{subfigure}\hfill
    \begin{subfigure}[b]{0.46\textwidth}
         \centering
    \includegraphics[width=\textwidth]{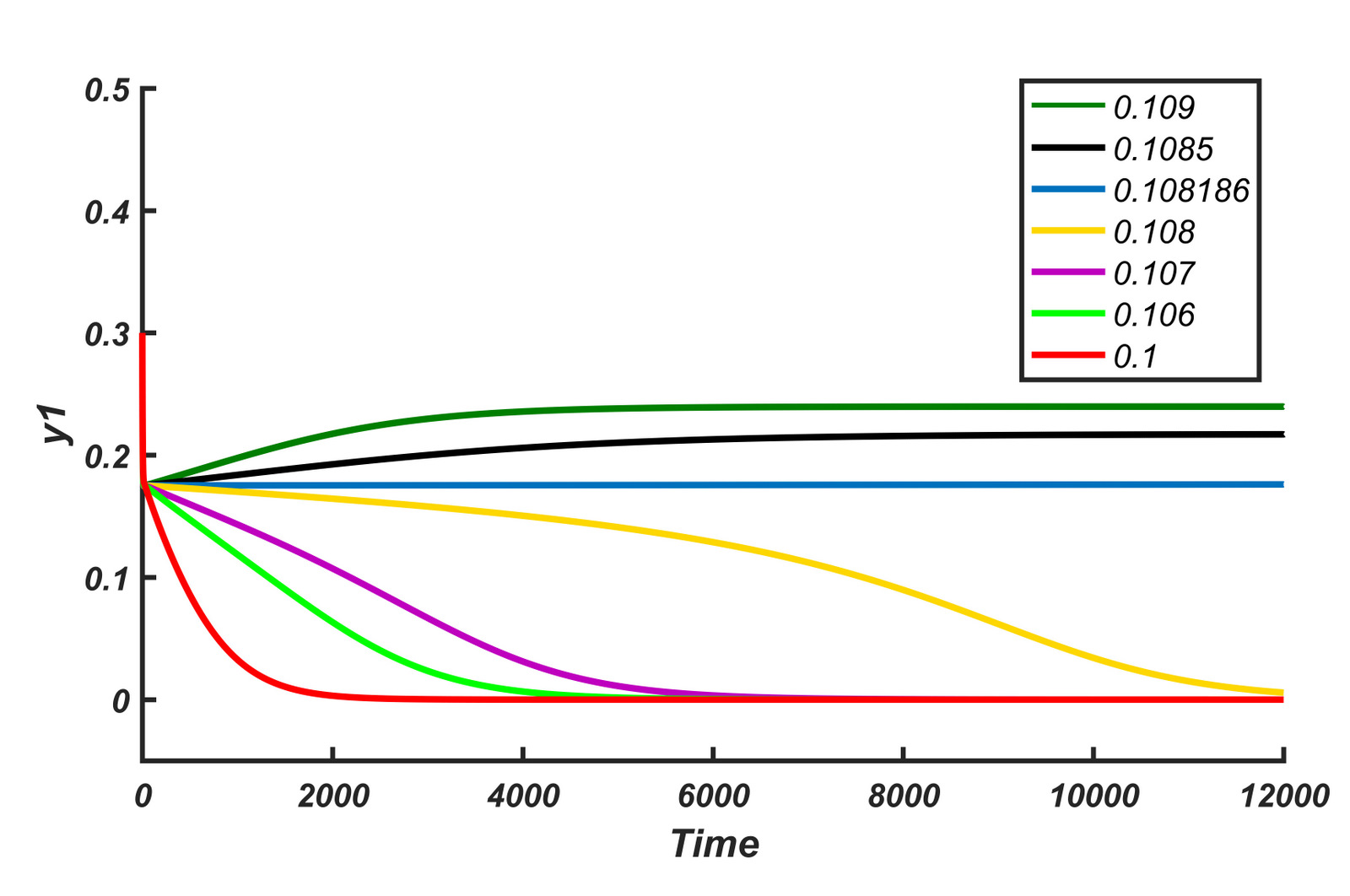}
    \caption{}
    \label{by1}
    \end{subfigure}\hfill
    \begin{subfigure}[b]{0.46\textwidth}
         \centering
    \includegraphics[width=\textwidth]{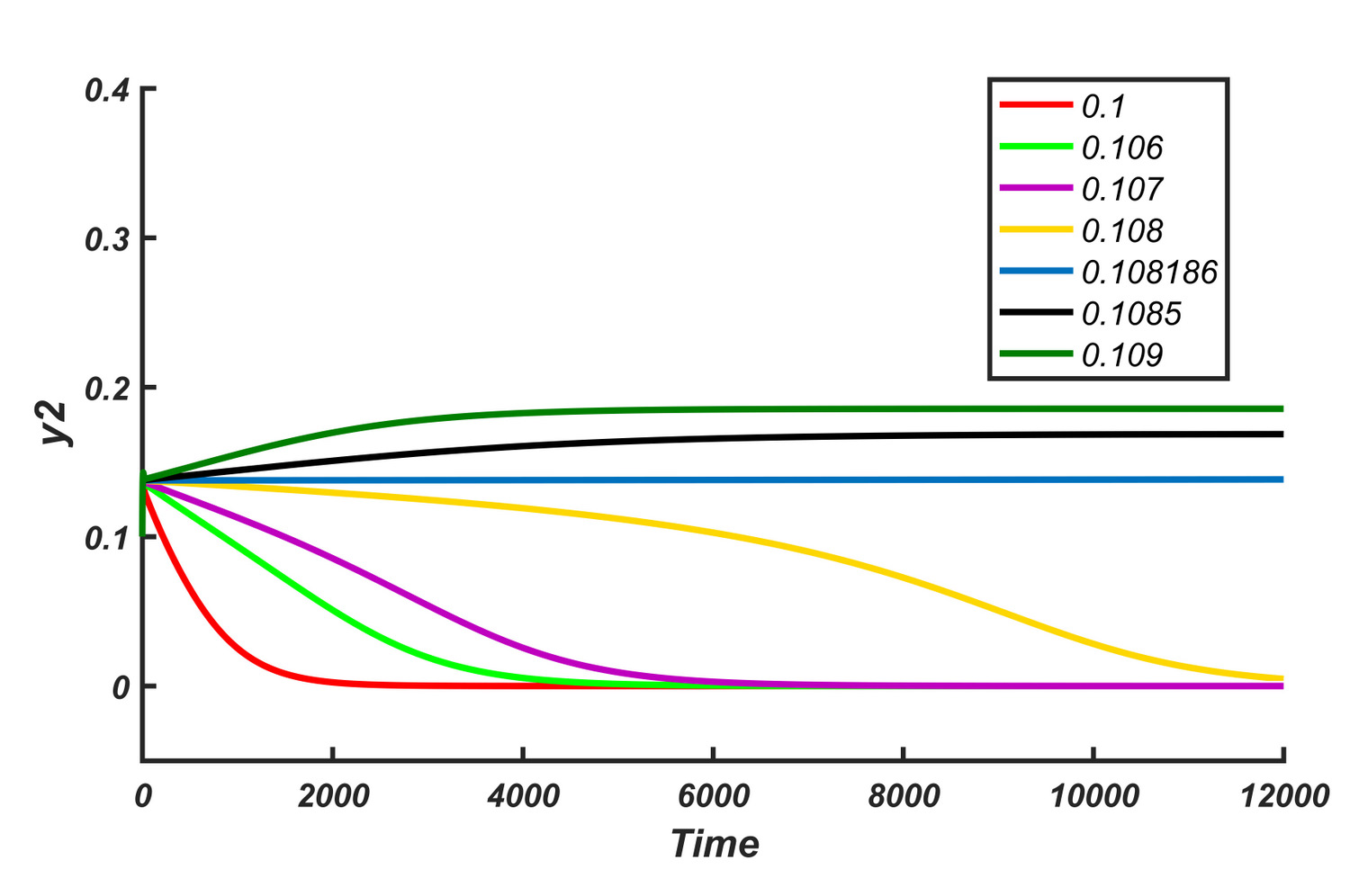}
    \caption{}
    \label{by2}
    \end{subfigure}\hfill
    \begin{subfigure}[b]{0.46\textwidth}
         \centering
    \includegraphics[width=\textwidth]{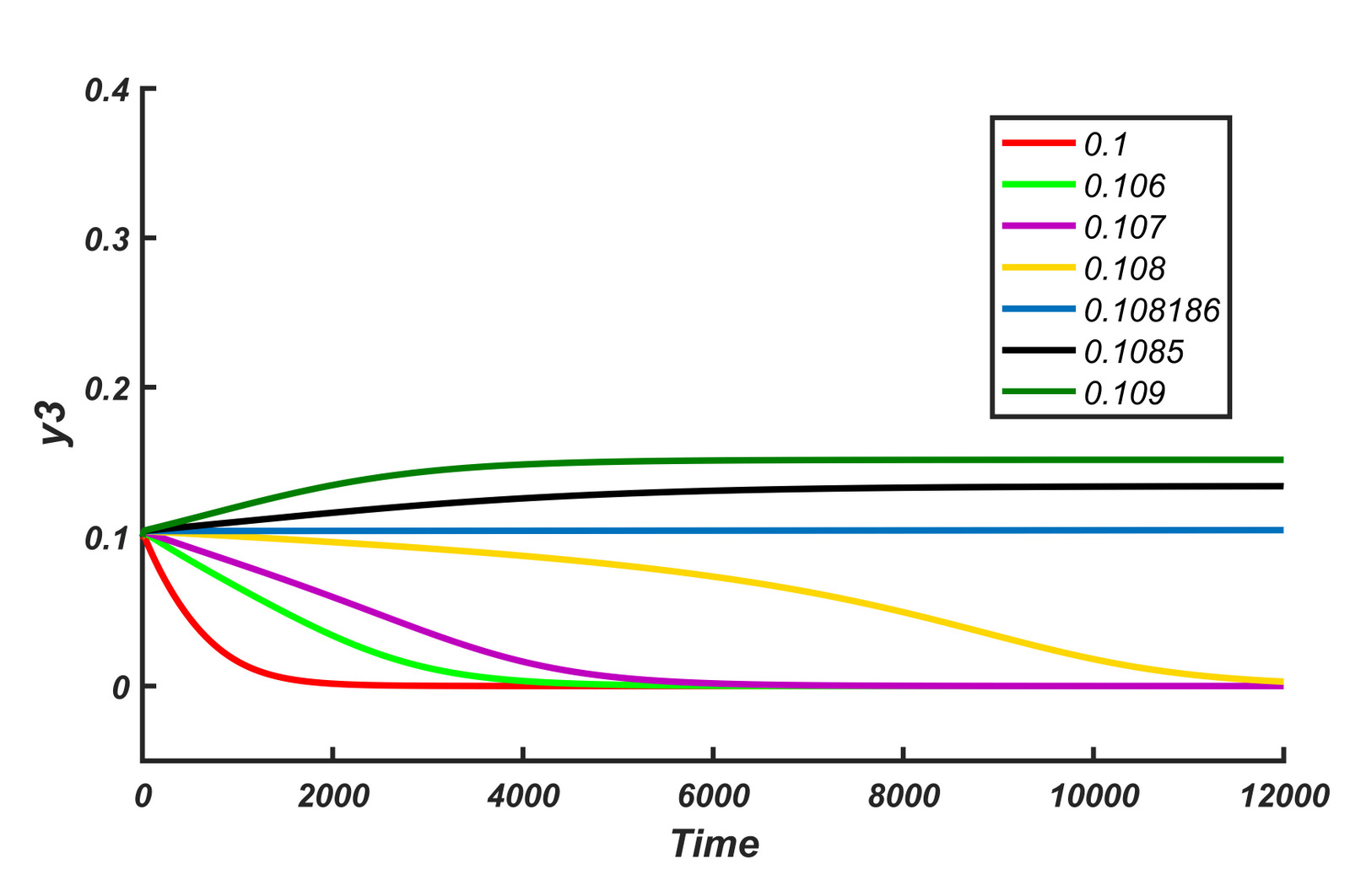}
    \caption{}
    \label{by3}
    \end{subfigure}
    \caption{Equilibrium states of the system \ref{ma-eq} at different parametric values of $b$ in the neighbourhood of its saddle node bifurcation point in context of (a)prey population, (b)infant predator population, (c)juvenile predator population, (d)adult predator population.}
    \label{b-lim}
\end{figure}

 Figure \ref{b-lim} depicts the way the equations of bio-system (\ref{ma-eq}) behaves in the vicinity of the saddle node bifurcation point. Varying the parametric values of the transition rate of infant predators, $b$ in the neighborhood of $0.108186$, it is seen that for $b<0.108186$, equilibrium $E_2$ is stable, but the closer it is to the bifurcation point, the more time would the compresent initial population take to reach the prey-only equilibrium point. Contrariwise, for $b>0.108186$, the compresent equilibrium is espied to be stable.
\begin{figure}[H]
    \centering
    \begin{subfigure}[b]{0.46\textwidth}
         \centering
    \includegraphics[width=\textwidth]{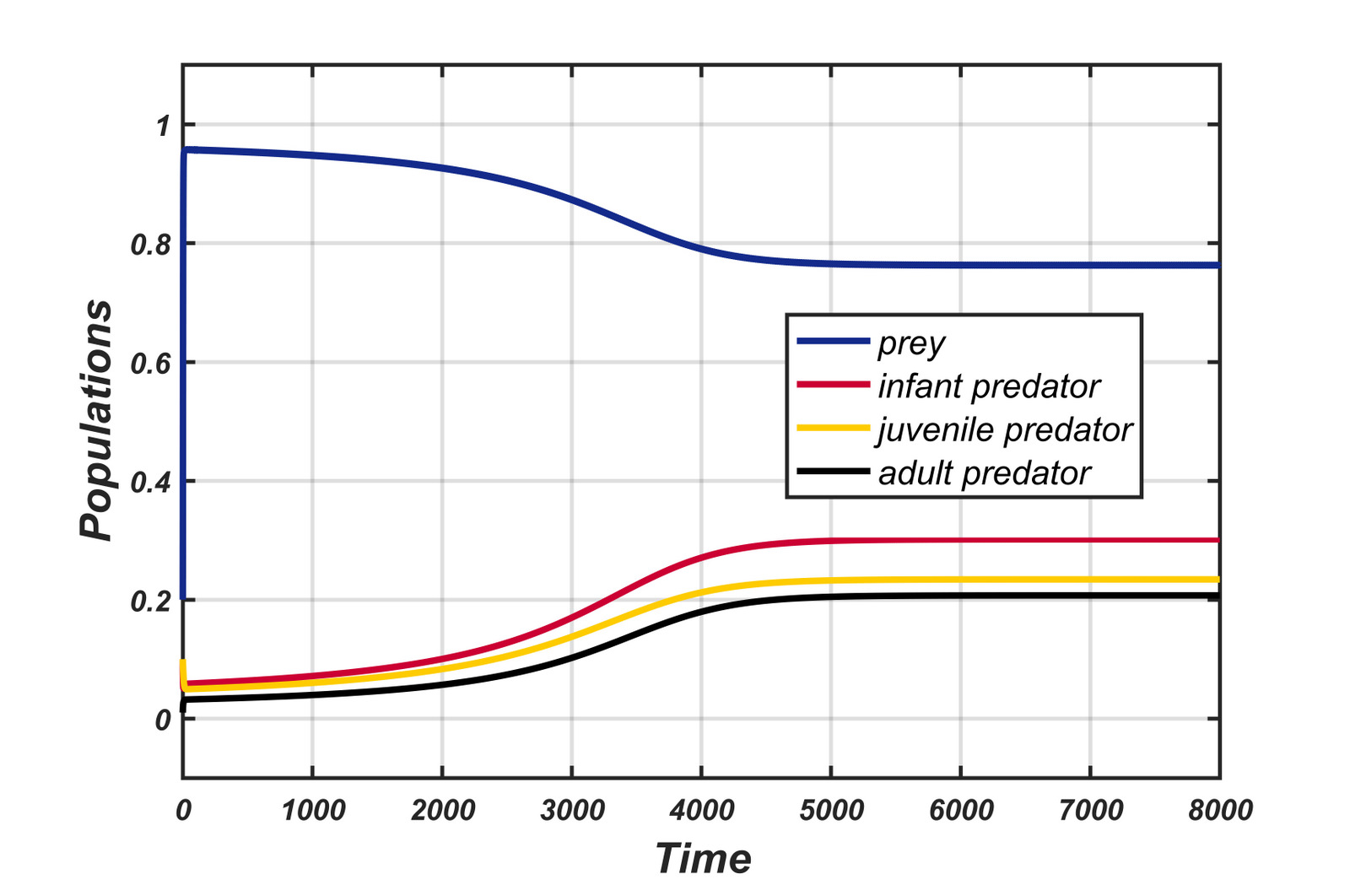}
    \caption{}
    \label{bi1}
    \end{subfigure}\hfill
    \begin{subfigure}[b]{0.46\textwidth}
         \centering
    \includegraphics[width=\textwidth]{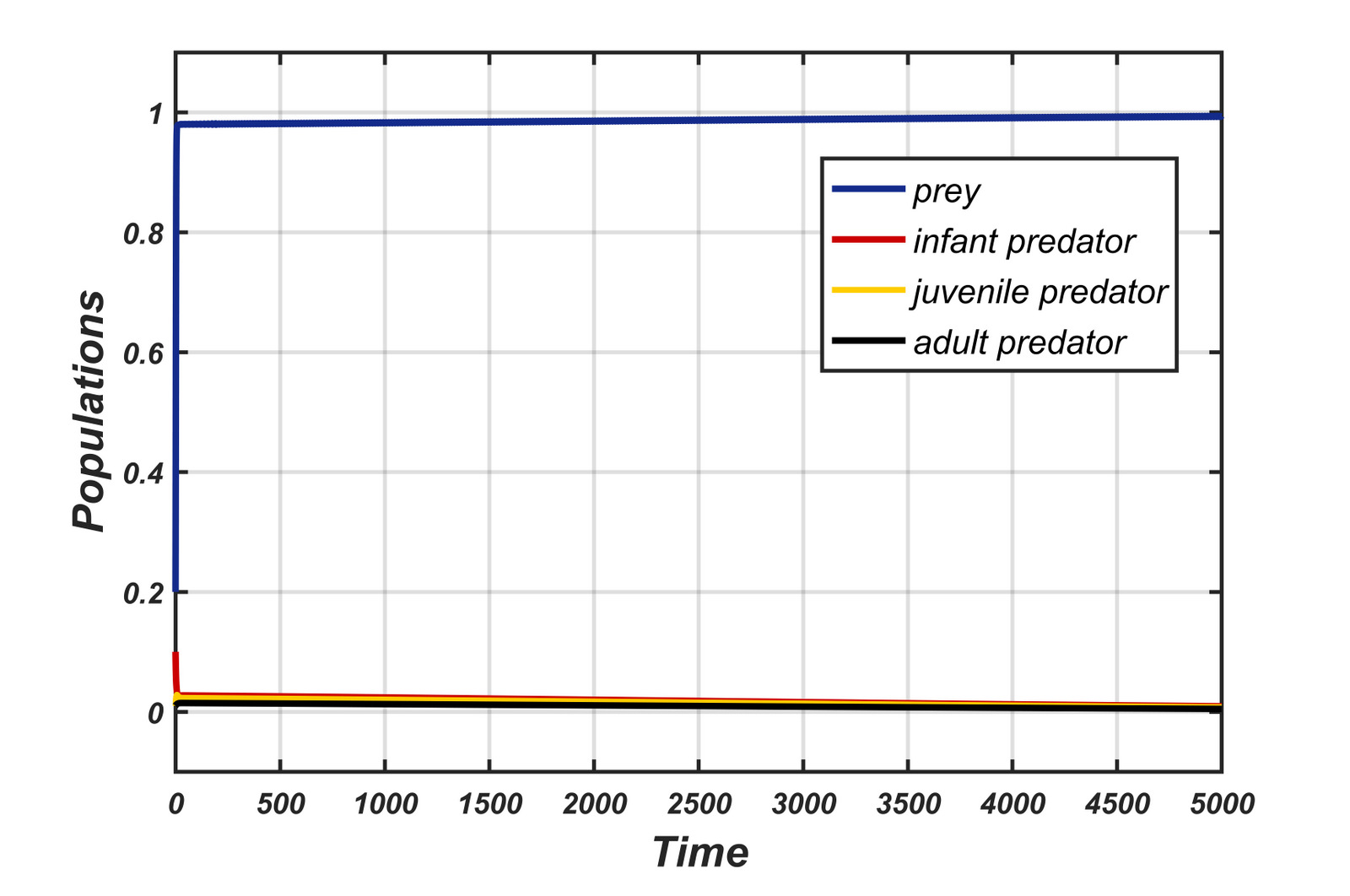}
    \caption{}
    \label{bi2}
    \end{subfigure}\hfill
     \begin{subfigure}[b]{0.5\textwidth}
         \centering
    \includegraphics[width=\textwidth]{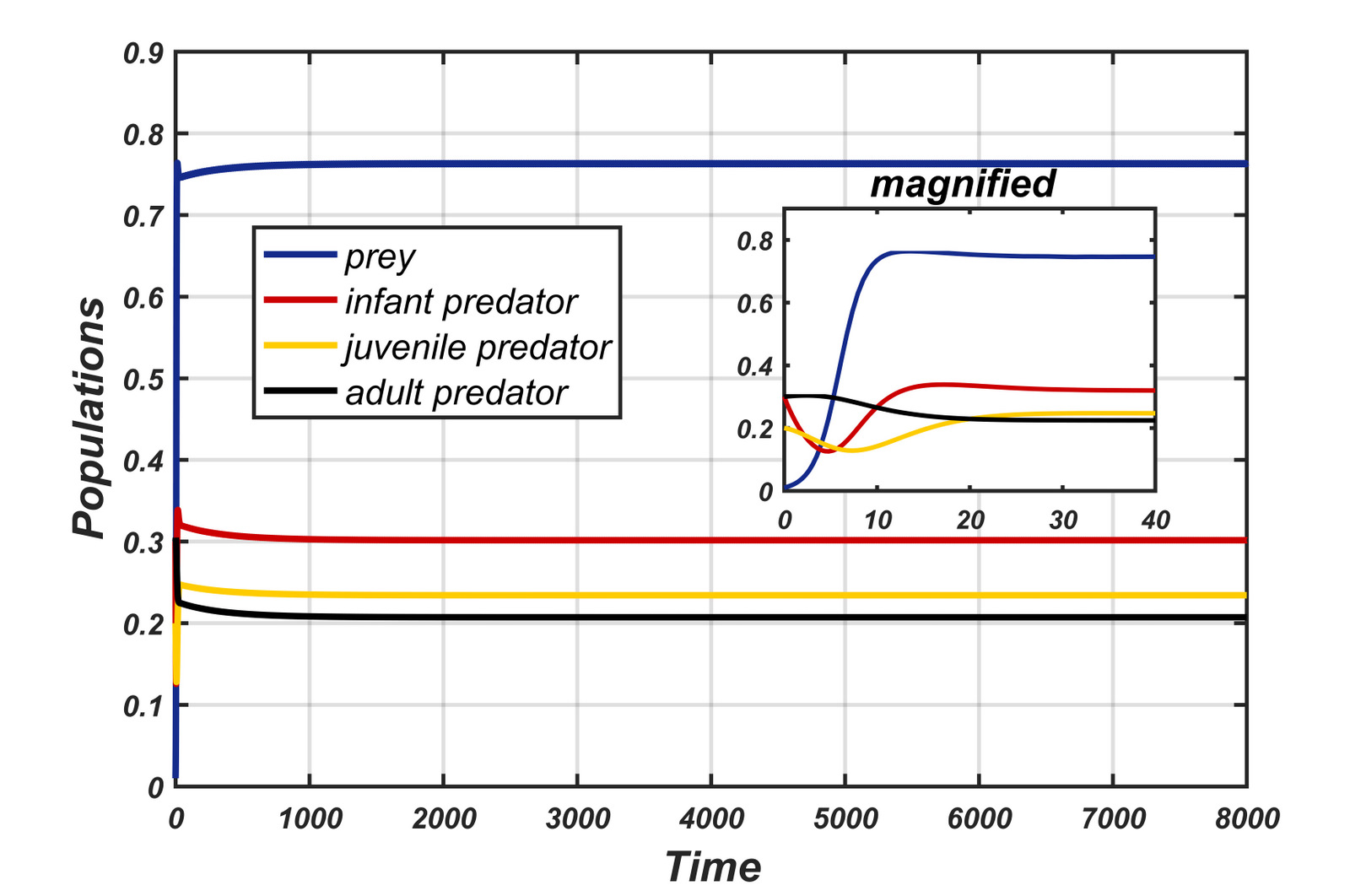}
    \caption{}
    \label{bi3}
    \end{subfigure}
    \caption{Depicting bi-stability of the equations of bio-system \ref{ma-eq} for the same parametric values with different initial population values, i.e., (a) at initial value $(0.2,0.1,0.1,0.01)$, (b) at initial value $(0.2,0.1,0.01,0.01)$, (c) at initial value $(0.01,0.3,0.2,0.3)$.}
    \label{bi-stab}
\end{figure}
\begin{figure}[H]
   \begin{subfigure}[b]{0.5\textwidth}
         \centering
    \includegraphics[width=\textwidth]{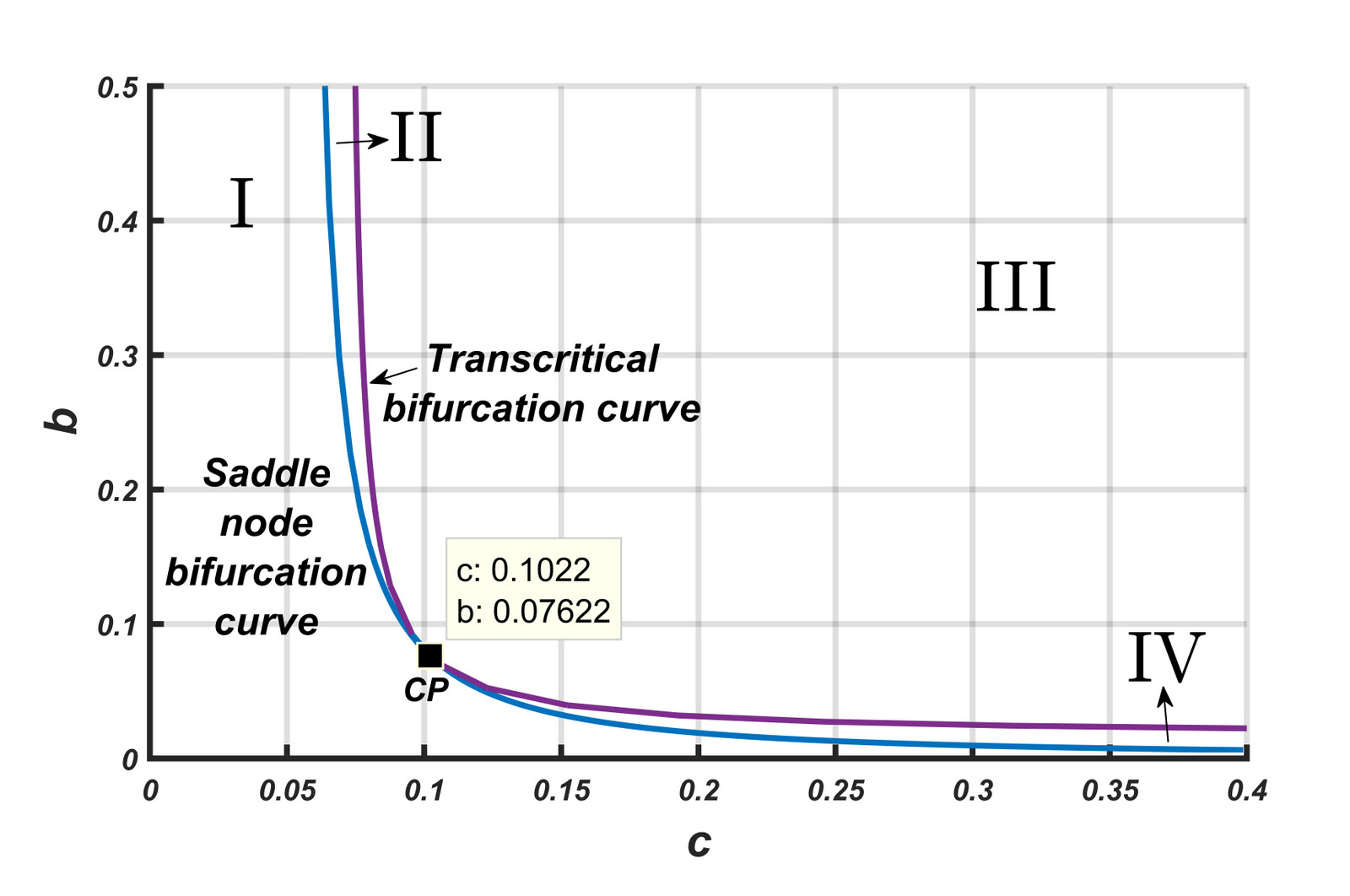}
    \caption{}
    \label{bipara1}
    \end{subfigure}\hfill
    \begin{subfigure}[b]{0.5\textwidth}
         \centering
    \includegraphics[width=\textwidth]{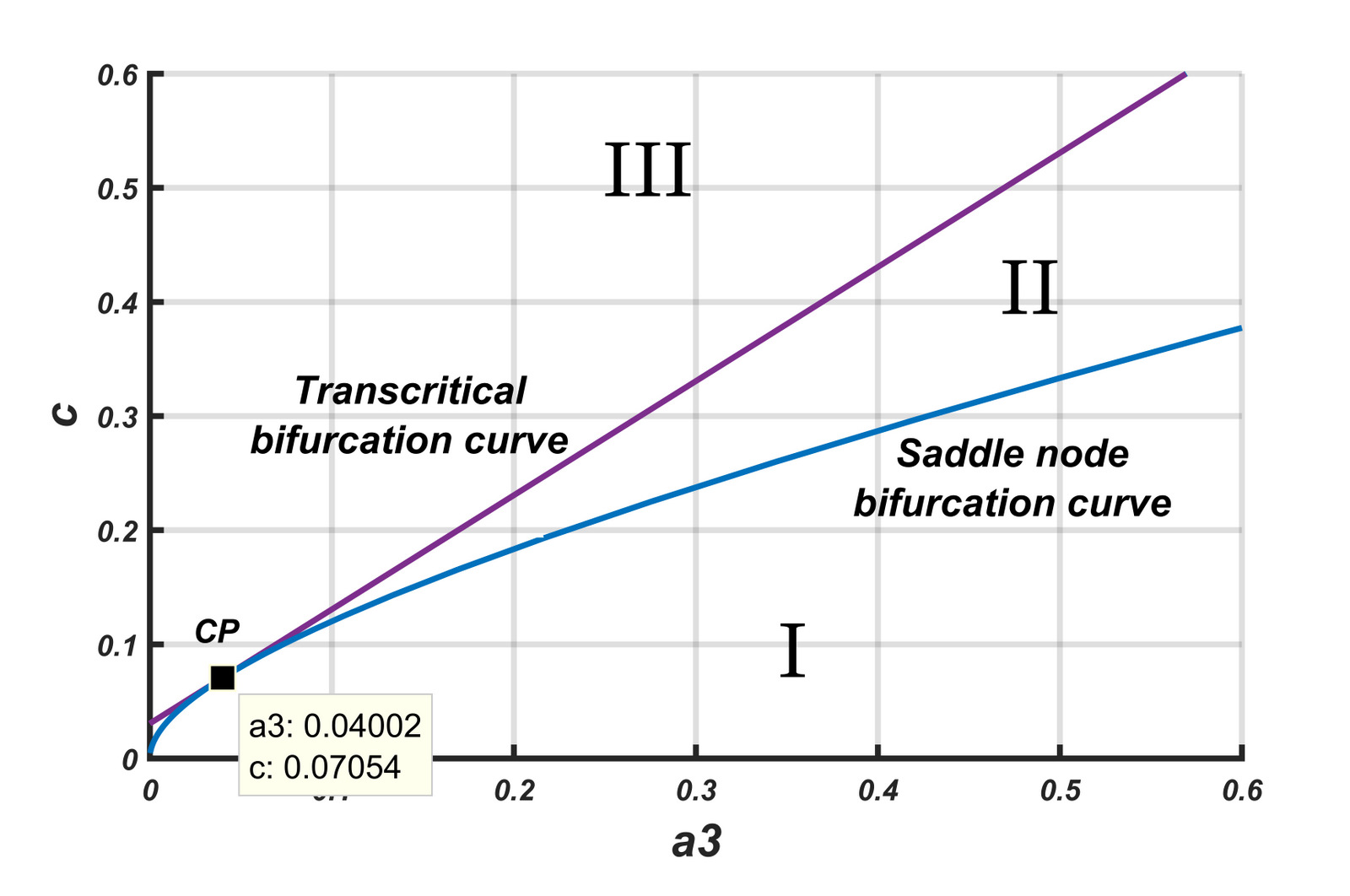}
    \caption{}
    \label{bipara2}
    \end{subfigure}
    \caption{Bi-parametric bifurcations: (a) transition rate of juvenile predator i.e., $c$, and transition rate of infant predator i.e., $b$  (b) rate of maturation delay in juvenile predators i.e, $a_3$ and transition rate of juvenile predator i.e., $c$. Here, $CP$ denotes cusp point}
    \label{bipara}
\end{figure}
\par Bi-stability is a fascinating circumstance whereby there is an existence of multiple attractors for the same set of parametric values, i.e., the solution trajectories converge to different attractors, which is extensively conditional to the initial populations. In our system (\ref{ma-eq}), there is a bi-stability between the compresent equilibrium point $E_4$ and the prey-only equilibrium point $E_2$. As can be evidenced from figure \ref{b-equ}, for $b<0.114706$, the prey-only equilibrium point is stable, and for $b>0.108186$, the compresent equilibrium is stable; which is to say that for $b \in (0.108186,0.114706)$, both the equilibrium points $E_2$ and $E_4$ are stable. 
Using the set of parametric values from table 1 (i.e., for $b=0.112\; \; \in (0.108186,0.114706)$ ) at different initial populations, figure \ref{bi-stab} is visualized. In figure \ref{bi1}, the initial population of the predator species is sparse, the time evolution of the prey population can be espied to rapidly increase towards the environment’s carrying capacity and then gradually decrease to the optimum population size for the co-existence of both species, while at the same time, all the stages of the predator population steadily increase to their respective optimum population sizes. But when the dimensions of the population of predator species is incommensurate, as is in figure \ref{bi2}, then unlike the previous scenario, the predator species are not in the position to revive again, and thus the equations of bio-system settles into the prey-only equilibrium point. Furthermore, when the initial biomass of prey population is low, as is in figure \ref{bi3}, the total biomass of the predator species decreases initially due to scarcity of food resources, giving the prey an opportunity to increase its biomass, after which the predator population size too increases, and thereby settling into the co-existing optimum population sizes, i.e., $E_4(0.762, 0.30159, 0.23419, 0.207)$. This can be espied in the magnified window of the figure \ref{bi3}. Interestingly, this same scenario, identified as ‘bloom’ phenomenon \cite{khajan1, khajan2}, is witnessed even when the initial prey population is taken incredibly low (for example the initial prey population is $0.1*10^{-9}$).  Therefore, it can be concluded that the prey-only equilibrium is attainable in this region of bi-stability only when the total predator population size is sufficiently low.\\
At this particular set of values as given in table 1, the similar scenario as witnessed in figures \ref{b-equ},\ref{b-lim},\ref{bi-stab} for the transition rate of infant predators, also occurs for transition rate of juvenile predators($c$), transformation rate($u$) and, predation rate of adult predators($a_2$) , and the scenario corresponding to the rate of maturity is also akin to this scenario, albeit in the reverse  orientation, i.e., increasing the parametric value of rate of maturation delay of juvenile predators ($a_3$), spawns the stability of equilibrium $E_2(1,0,0,0)$. While the variations in the parametric values of predation rate by juvenile predator ($a_1$) causes no changes in the system of equations.\\
\par So as to investigate the joint impact of the transition rate of infant predator $b$, and the transition rate of juvenile predator $c$, the bi-parametric bifurcation illustration is given in figure \ref{bipara1}. The entire bi-parametric region is segregated into regions of stability and bi-stability for different equilibrium points by the curves of transcritical bifurcation (solid violet colored lines) and saddle node bifurcation (solid blue colored lines). Region $I$ and region $III$ are the areas of mono-stability for the prey-only equilibrium $E_2$ and compresent equilibrium $E_4$ respectively, while region $II$ is the area of bi-stability for both the equilibrium points $E_2$ and $E_4$. For the parametric values from the left and lower side of transcritical bifurcation curve, equilibrium point $E_2(1,0,0,0)$ is stable else it is unstable. Similarly, for the values from only the right and upper side of the saddle node bifurcation curve, the equilibrium point $E_4$ is stable, region $II$ which is common to both the cases becomes a bi-stable area. The same bi-stable scenario should have been for region $IV$, but the presence of cusp bifurcation point at $(c,b)=(0.1022,0.07622)$ re-configures that. The point where transcritical bifurcation curve meets saddle node bifurcation curve is known as cusp bifurcation point. Because of cusp bifurcation, there is an existence of a phenomenon called ‘hysterics’. ‘Hysterics’, spawned because of the disappearance of the traced stable compresent equilibrium $E_4$ via the saddle node bifurcation, is a catastrophic ‘jump’ to a different stable equilibrium, which in our case is the prey-only equilibrium $E_2$. Thus, in region $IV$ only $E_2$ is stable. \\
In similar manner, figure \ref{bipara2} illustrates the dynamical relation between the parameters, rate of maturation delay $a_3$ and the transition rate $c$ of the juvenile predator. The cusp bifurcation point is at $(a_3,c)=(0.04002,0.7054)$. (Stability regions are same as that of the previous figure \ref{bipara1}).\\
\par The parametric values as given in table 2 is considered, except $c$, the transition rate of juvenile predators into adult stage, which is varied, to investigate the occurrence of Hopf bifurcation (also known as Andronov-Hopf bifurcation) point. Analytically, at $c=0.03598345(=c_h)$, we have for the equilibrium point $E_3(0.452951, 0.609794, 0.755778, 0.116978)$, $\epsilon_1=0.547,\, \epsilon_4=-0.0001, \, \epsilon_1 \epsilon_2 - \epsilon_3=0.02, \, \epsilon_1\epsilon_2 \epsilon_3-\epsilon_3^2-\epsilon_1^2 \epsilon_4=0.00007$; and for the equilibrium point $E_4(0.063266, 0.649761,0.37137,0.89239)$, we have $\epsilon_1=0.187, \, \epsilon_4=0.00003, \, \epsilon_1 \epsilon_2 - \epsilon_3=0.0013, \, \epsilon_1\epsilon_2 \epsilon_3-\epsilon_3^2-\epsilon_1^2 \epsilon_4=0$, therefore the conditions of Hopf Bifurcation along with the transversality condition ($\frac{\epsilon_1\epsilon_2 \epsilon_3-\epsilon_3^2-\epsilon_1^2 \epsilon_4}{a_3}=-0.000025 \neq 0$) is satisfied for $E_4$. Therefore, there is an emergence of limit cycle around the equilibrium point $E_4$, at $c_h$, according to Theorem \ref{ho-thm}. Next for the examination of direction and stability of the limit cycle, as specified in Theorem \ref{ho-dir-thm}, we have, 
 $g_{11}=-0.018595+0.0328119 i, \; g_{02}= 0.0751727+0.397665 i, \; g_{20}= -0.112363-0.332041 i, \; g_{21}= -0.00766054+0.438565 i $. Hence we get, $C_1(0)= -0.0214226+0.134726i $, $\gamma_2= 0.0580115 $, and $\beta_2= -0.0428451$. Here, $\gamma_2>0$ implies a supercritical Hopf bifurcation, and $\beta_2<0$ implies the limit cycle being stable. 
Figure \ref{c-2} depicts the emergence of limit cycles at $c_h$. For the values of $c$ greater than $c_h$, the system is unstable, and the Hopf bifurcation being a supercritical one, stable oscillation is witnessed around the equilibrium $E_4$ (figure \ref{c-1}), the greater the distance of the value of $c$ from $c_h$ the greater would be the amplitude of the oscillation. On the contrary, for $c<c_h$, the equations of bio-system evidences a stable focus point $E_4$ (figure \ref{c-3}). \\ 
\par
\begin{figure}[H]
    \centering
    \begin{subfigure}[b]{0.8\textwidth}
         \centering
    \includegraphics[width=\textwidth]{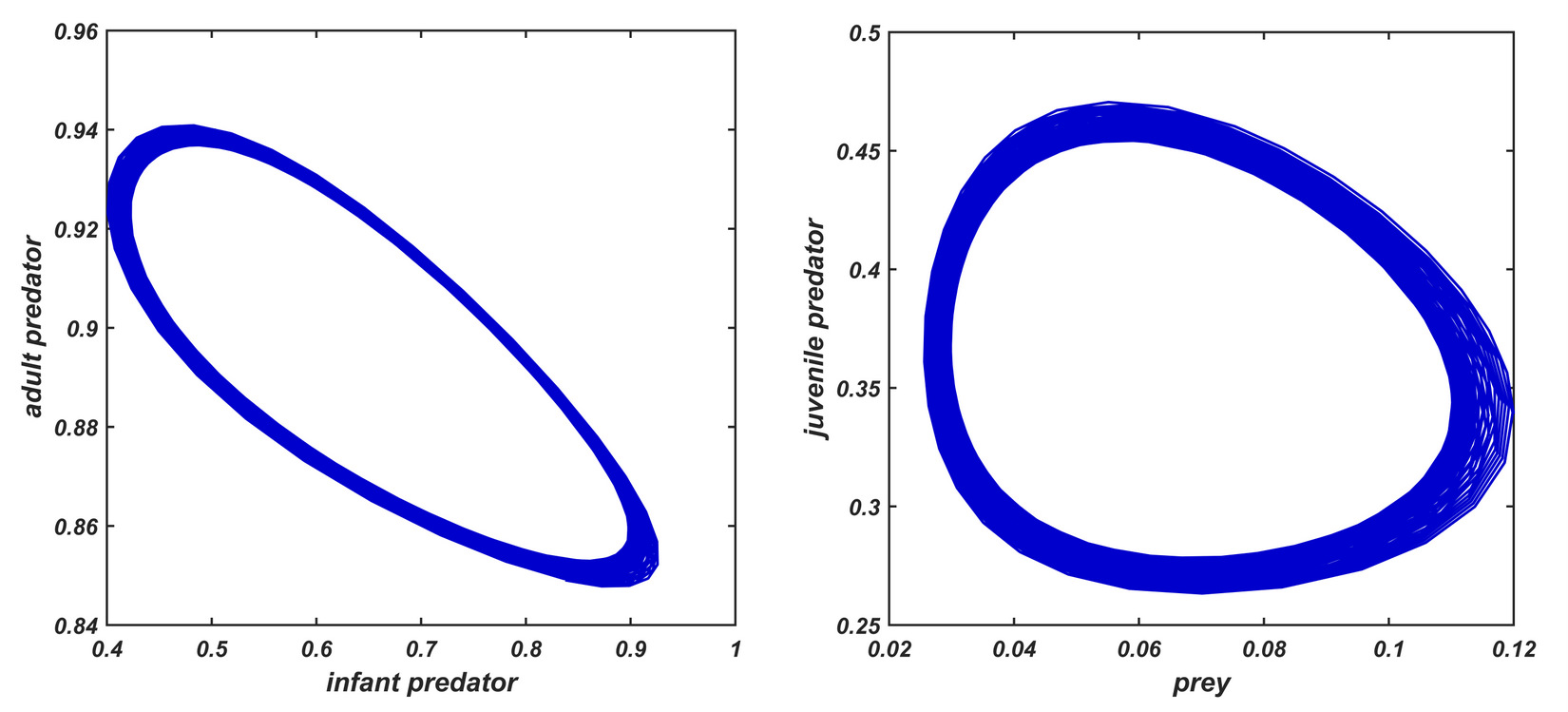}
    \caption{at $a_3=0.055$}
    \label{a3-1}
    \end{subfigure}\hfill
    \begin{subfigure}[b]{0.8\textwidth}
         \centering
    \includegraphics[width=\textwidth]{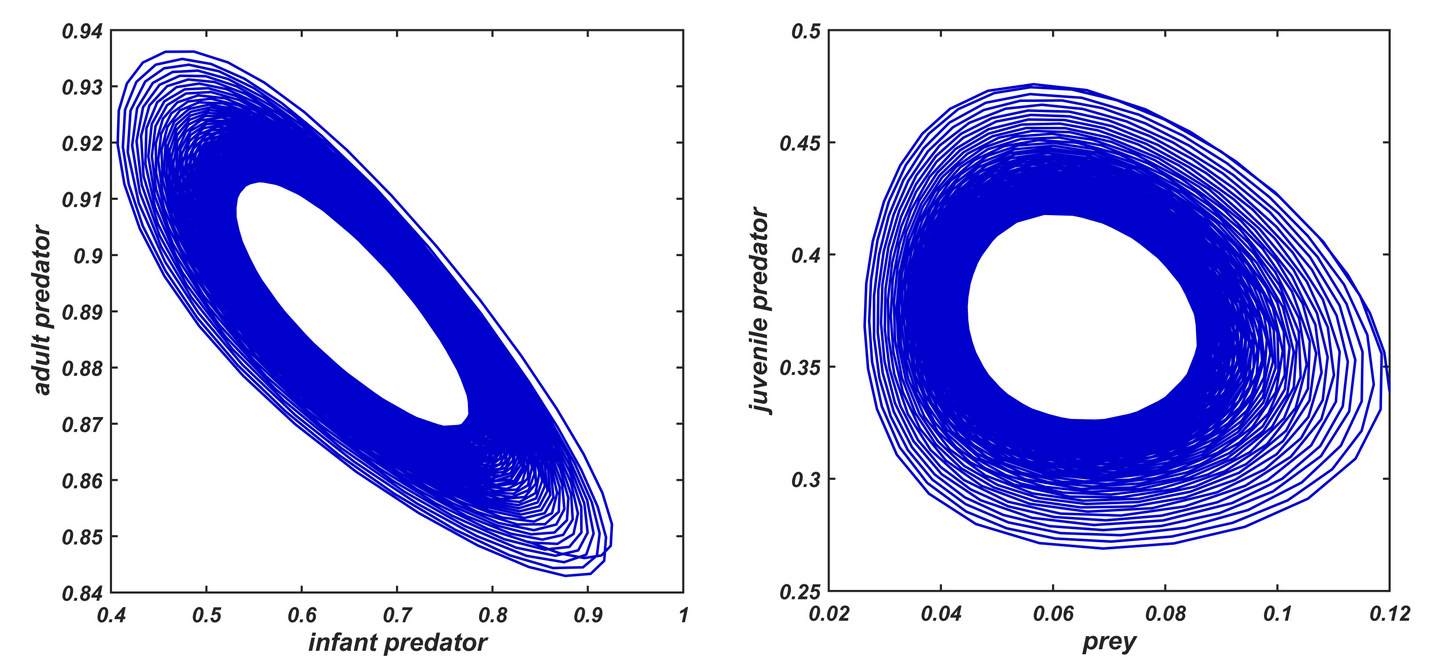}
    \caption{at $a_3=0.0604877$}
    \label{a3-2}
    \end{subfigure}\hfill
    \begin{subfigure}[b]{0.8\textwidth}
         \centering
    \includegraphics[width=\textwidth]{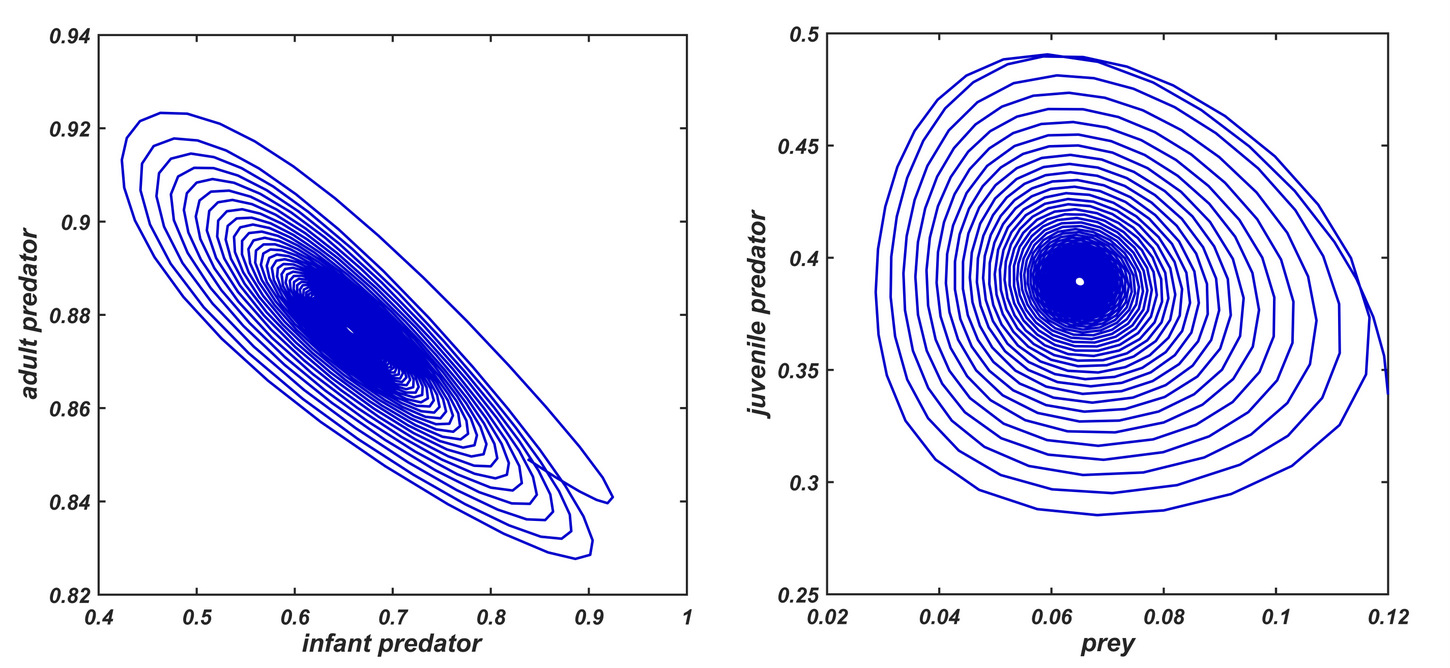}
    \caption{at $a_3=0.084$}
    \label{a3-3}
    \end{subfigure}
    \caption{Phase portrait in a bi-dimensional space depicting predator-prey population dynamics in the neighborhood of Hopf bifurcation point of the rate of maturation delay in juvenile predators, $a_3$: $(a)a_3<a_{3h}$, $(b) a_3=a_{3h}$, $(c)a_3>a_{3h}$, $a_{3h}$ is the Hopf bifurcation point.}
    \label{a3-2dim}
\end{figure}

 The Hopf bifurcation plot with respect to $a_3$, the rate of maturity delay in juvenile predators is illustrated in figure \ref{a3-2dim}, the bifurcating parametric value is $a_3=0.0604877(=a_{3h})$, when the values of the other parameters are taken from table 2. The first Lyapunov coefficient at $a_{3h}$ is found to be $-1.461335e^{-02}$, its negativity signifies the supercritical nature of Hopf bifurcation. For $a_3>a_{3h}$ (figure \ref{a3-3}), the trajectory with the initial co-habitance population spirals into the compresent equilibrium point $E_4$, while for $a_3<a_{3h}$, it spirals around the formerly stable equilibrium point $E_4$ (figure \ref{a3-1}). The stable limit cycle being born with its focus at the unstable compresent equilibrium point $E_4$ is seen in figure \ref{a3-2}  \\
\par
\begin{figure}[H]
   \begin{subfigure}[b]{0.5\textwidth}
         \centering
    \includegraphics[width=\textwidth]{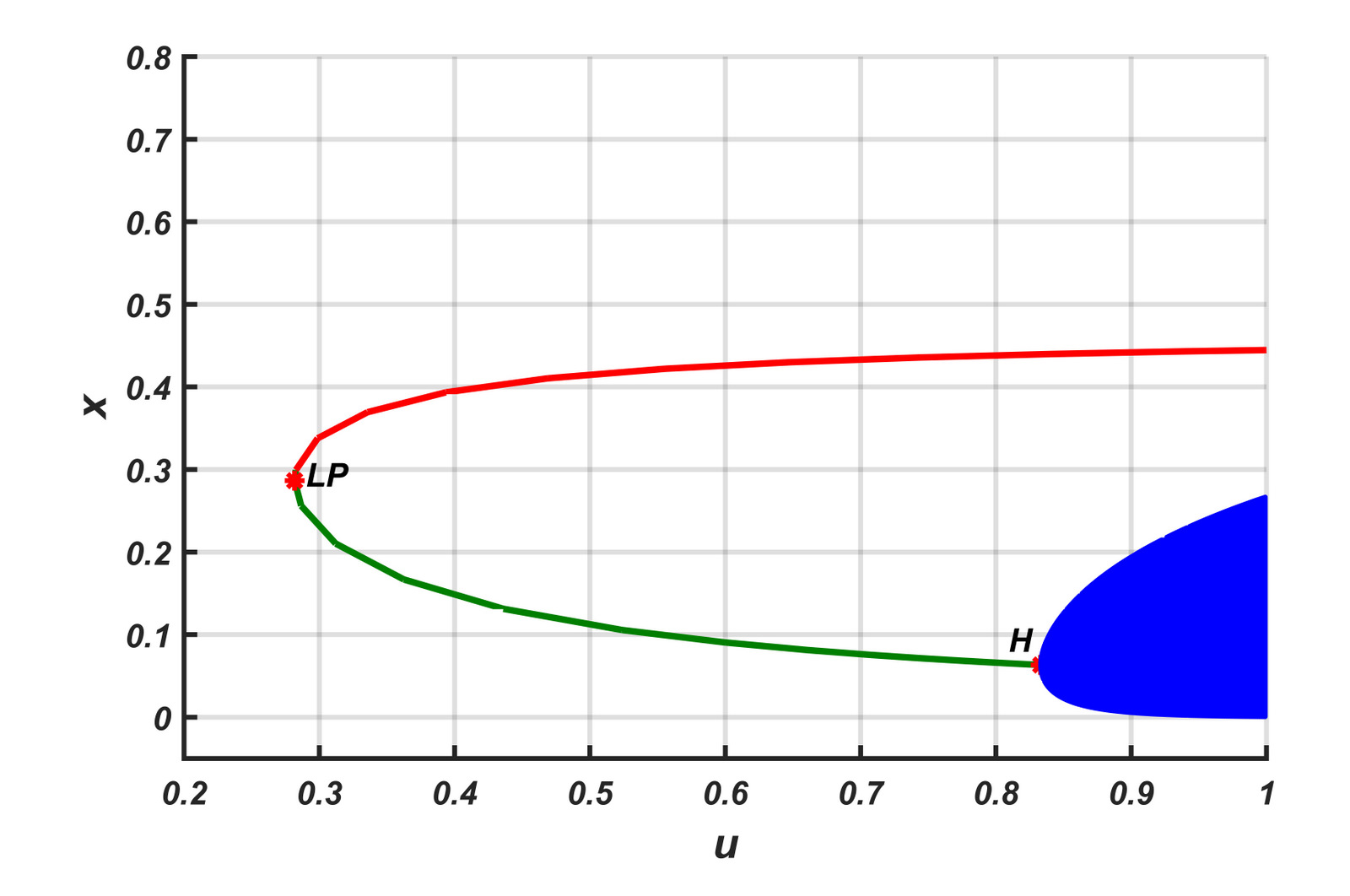}
    \caption{}
    \label{ux}
    \end{subfigure}\hfill
    \begin{subfigure}[b]{0.5\textwidth}
         \centering
    \includegraphics[width=\textwidth]{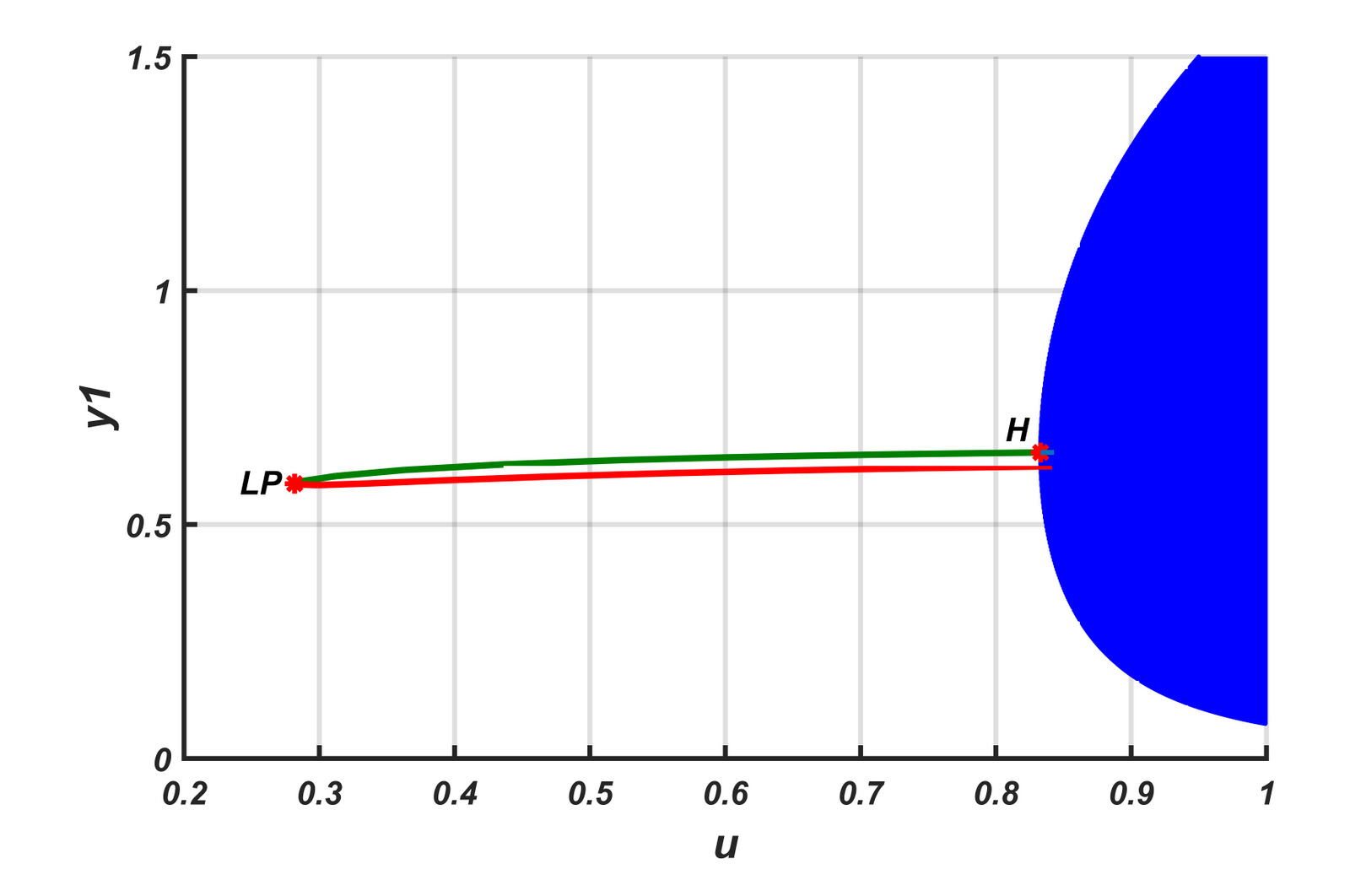}
    \caption{}
    \label{uy1}
    \end{subfigure}\hfill
    \begin{subfigure}[b]{0.5\textwidth}
         \centering
    \includegraphics[width=\textwidth]{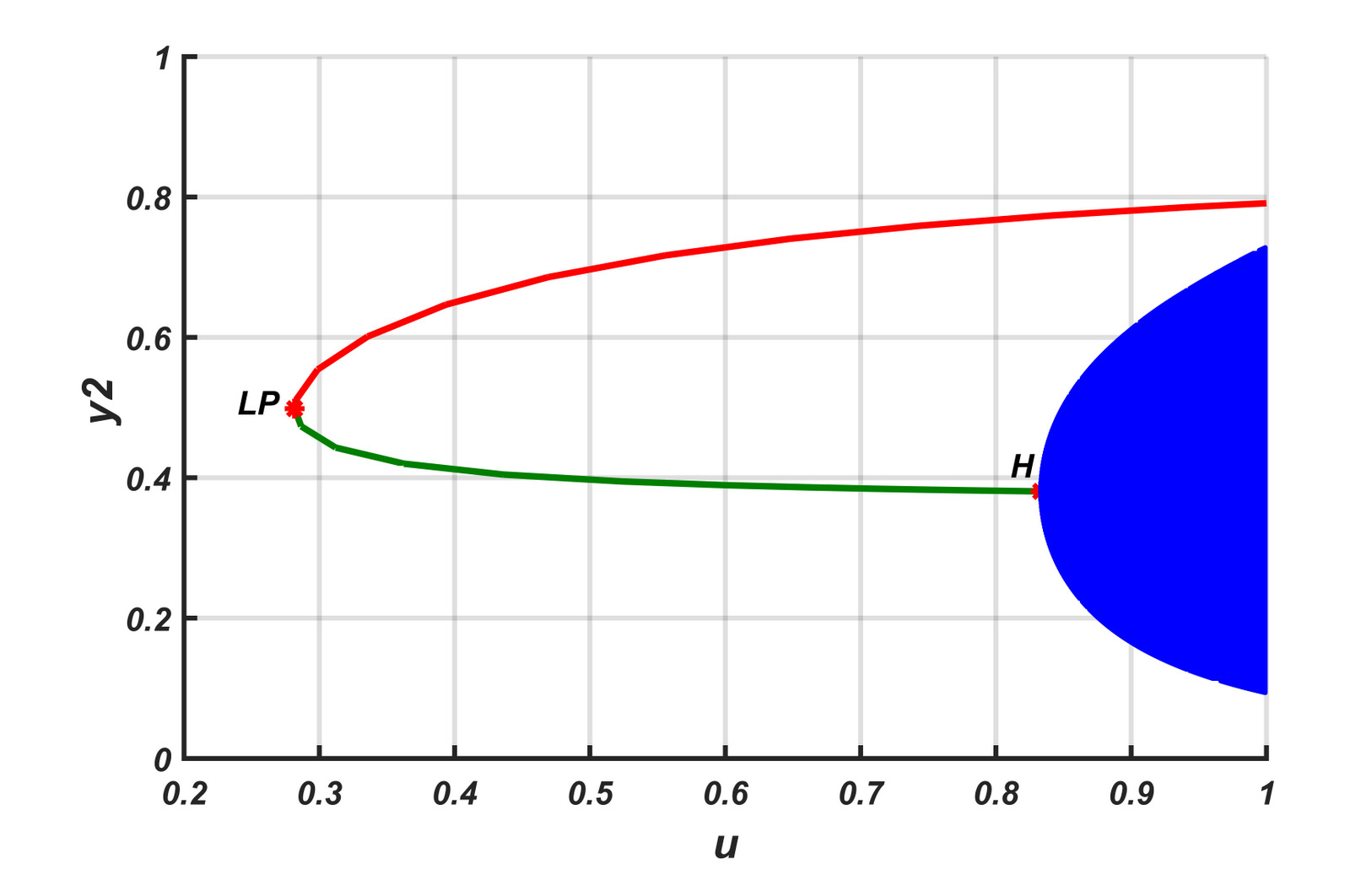}
    \caption{}
    \label{uy2}
    \end{subfigure}\hfill
    \begin{subfigure}[b]{0.5\textwidth}
         \centering
    \includegraphics[width=\textwidth]{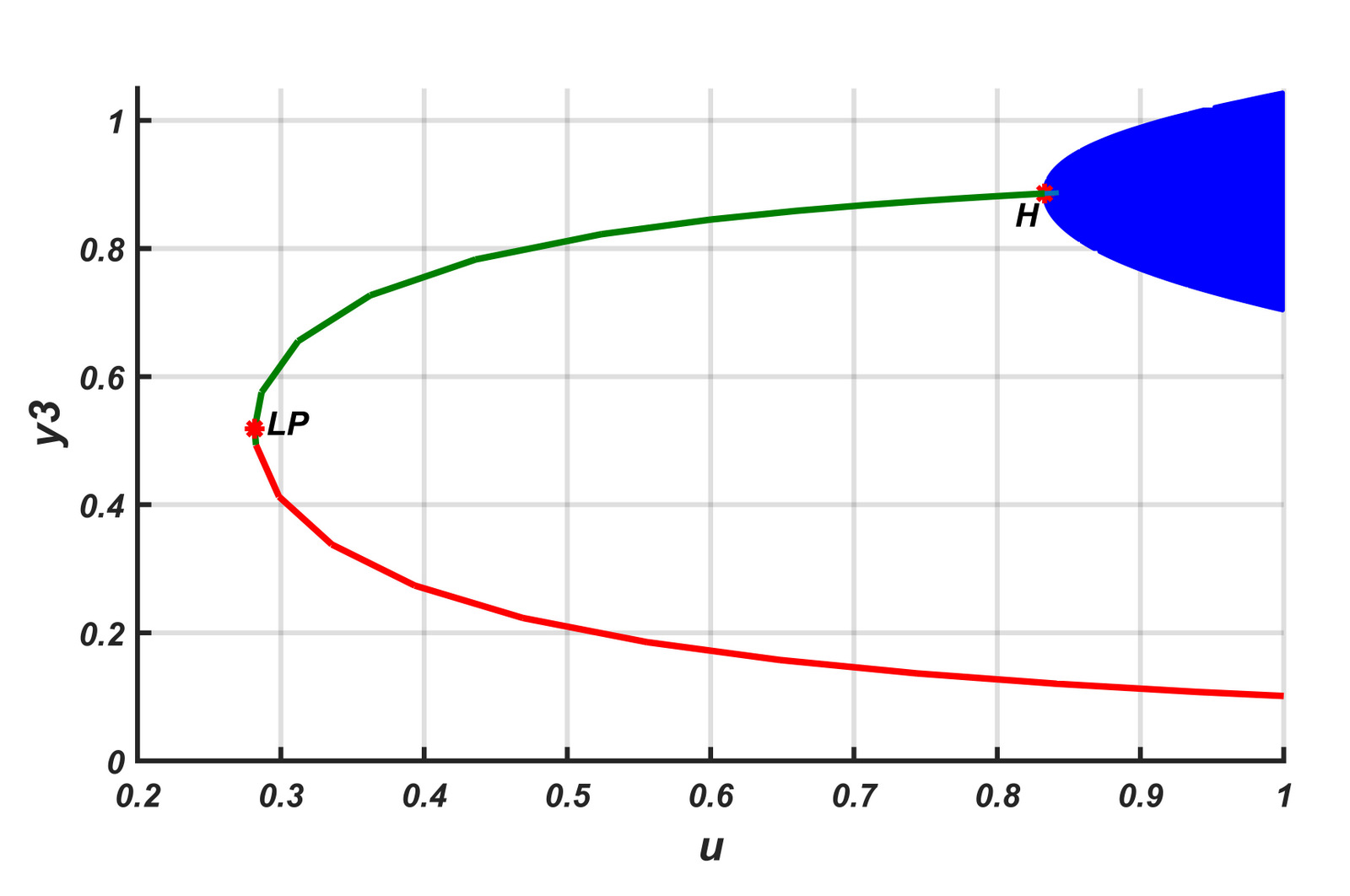}
    \caption{}
    \label{uy3}
    \end{subfigure}
    \caption{Equilibrium states of the prey, infant predator, juvenile predator and, adult predator populations with respect to the transformation rate of consumed prey into infant predator, $u$ for (a) prey population, (b) infant predator population, (c) juvenile predator population, (d) adult predator population. \textbf{$LP$} denotes the saddle node bifurcation point and, \textbf{$H$} denotes the Hopf bifurcation point.}
    \label{u-lh}
\end{figure}

 Next, the bifurcations of the system with respect to a vital parameter, $u$, the transformation rate of consumed prey into infant predators, is examined in figure \ref{u-lh} (other parametric values are as is in table 2).  At $u=0.281804$, the system evinces a saddle node bifurcation where a stable compresent equilibrium $E_4$ and an unstable compresent equilibrium $E_3$ collide, while at $u=0.833189$ the system experiences the Hopf bifurcation, which is supercritical as the first Lyapunov coefficient is $- 1.502700e^{-02}$. Therefore, for $0.281804<u<0.833189$, the equations of bio-system (\ref{ma-eq}) harbours one stable compresent equilibrium $E_4$, and when $u$ crosses the threshold value at $0.833189$, it exchanges its stability with stable limit cycles around the now unstable equilibrium $E_4$. An augmentation of the modulations in the populations of infant predators, juvenile predators, adult predators and prey can be seen with further increment of the value of $u$. 
 
\begin{figure}[H]
   \begin{subfigure}[b]{0.5\textwidth}
         \centering
    \includegraphics[width=\textwidth]{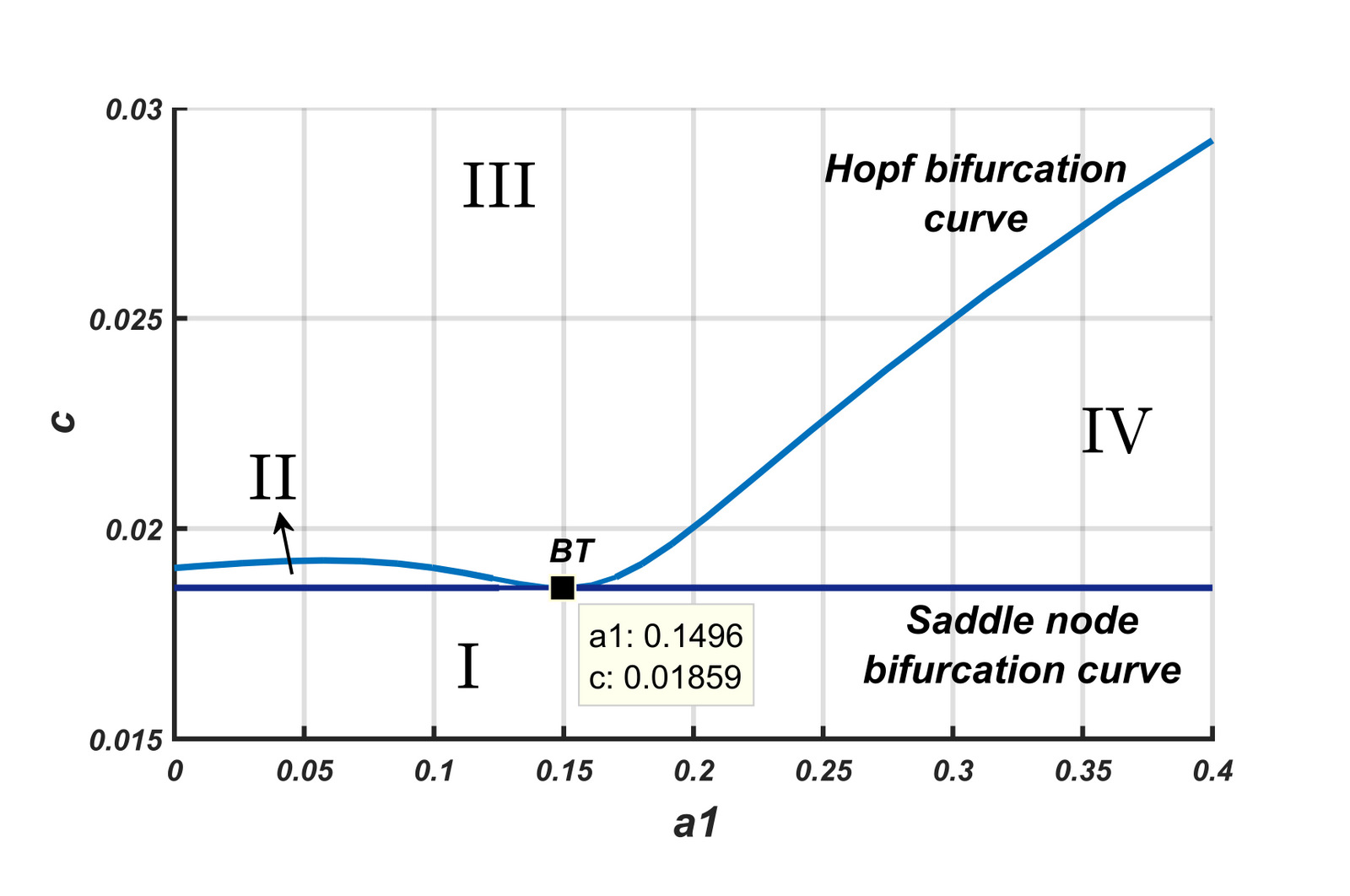}
    \caption{}
    \label{a12}
    \end{subfigure}\hfill
    \begin{subfigure}[b]{0.5\textwidth}
         \centering
    \includegraphics[width=\textwidth]{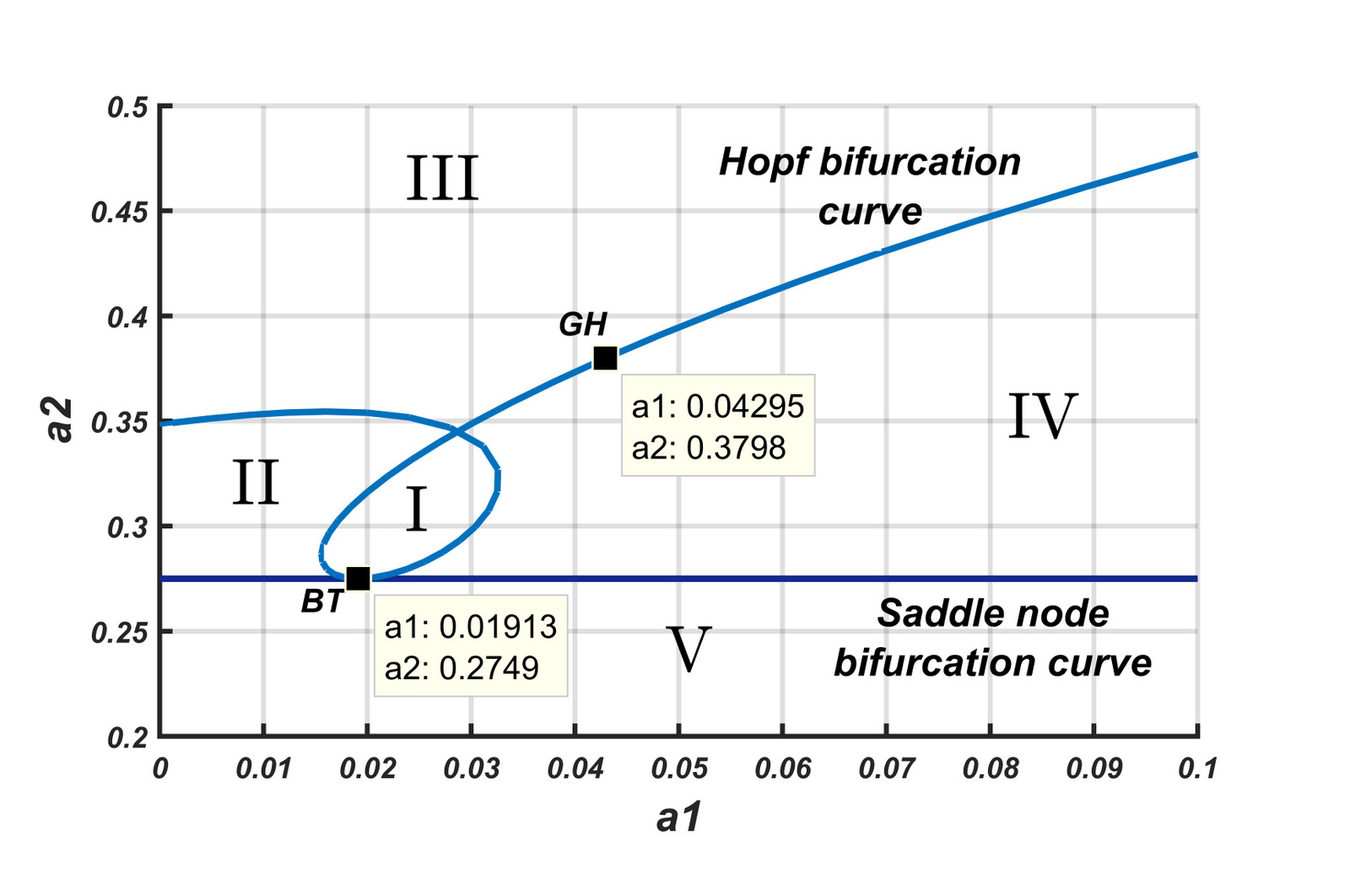}
    \caption{}
    \label{a11}
    \end{subfigure}
    \caption{Bi-parametric bifurcations: (a) predation rate of juveniles, $a_1$, and transition rate of juvenile predator, $c$  (b) predation rate of juveniles, $a_1$, and predation rate of adult predators, $a_2$.}
    \label{fig-a1}
\end{figure}
\begin{figure}[H]
   \begin{subfigure}[b]{0.5\textwidth}
         \centering
    \includegraphics[width=\textwidth]{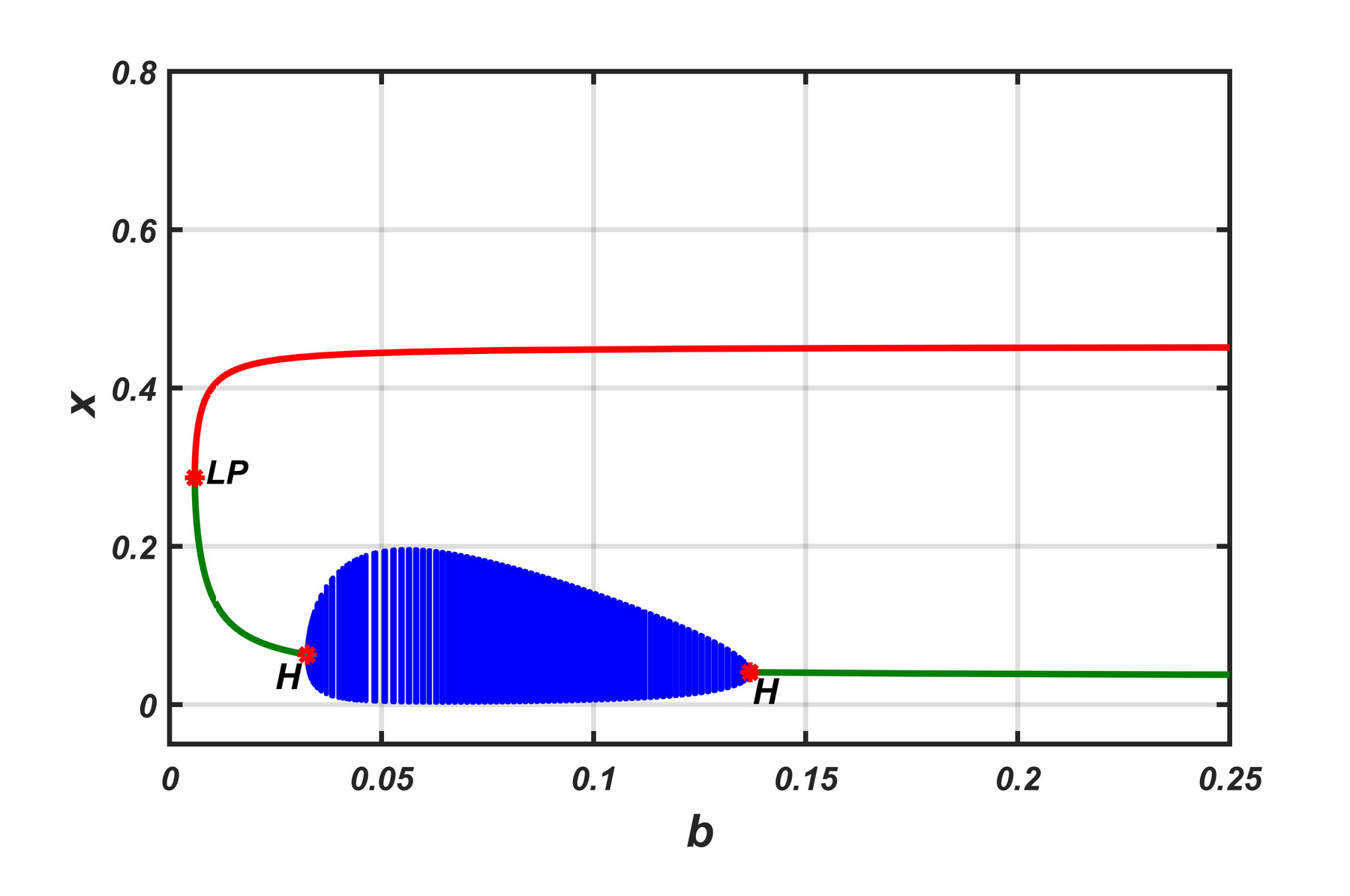}
    \caption{}
    \label{bhox}
    \end{subfigure}\hfill
    \begin{subfigure}[b]{0.5\textwidth}
         \centering
    \includegraphics[width=\textwidth]{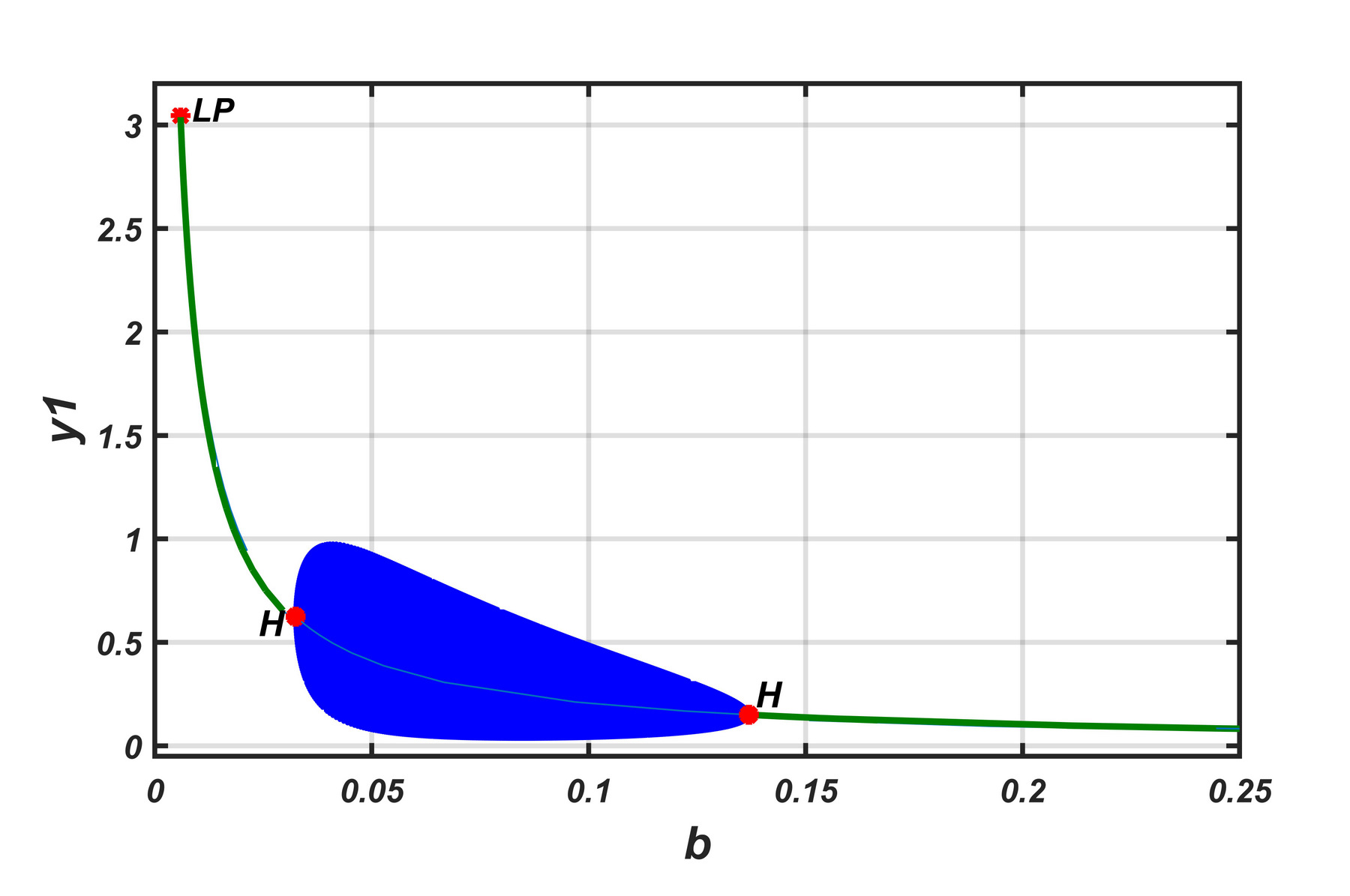}
    \caption{}
    \label{bhoy2}
    \end{subfigure}\hfill
    \begin{subfigure}[b]{0.5\textwidth}
    \centering
         \includegraphics[width=\textwidth]{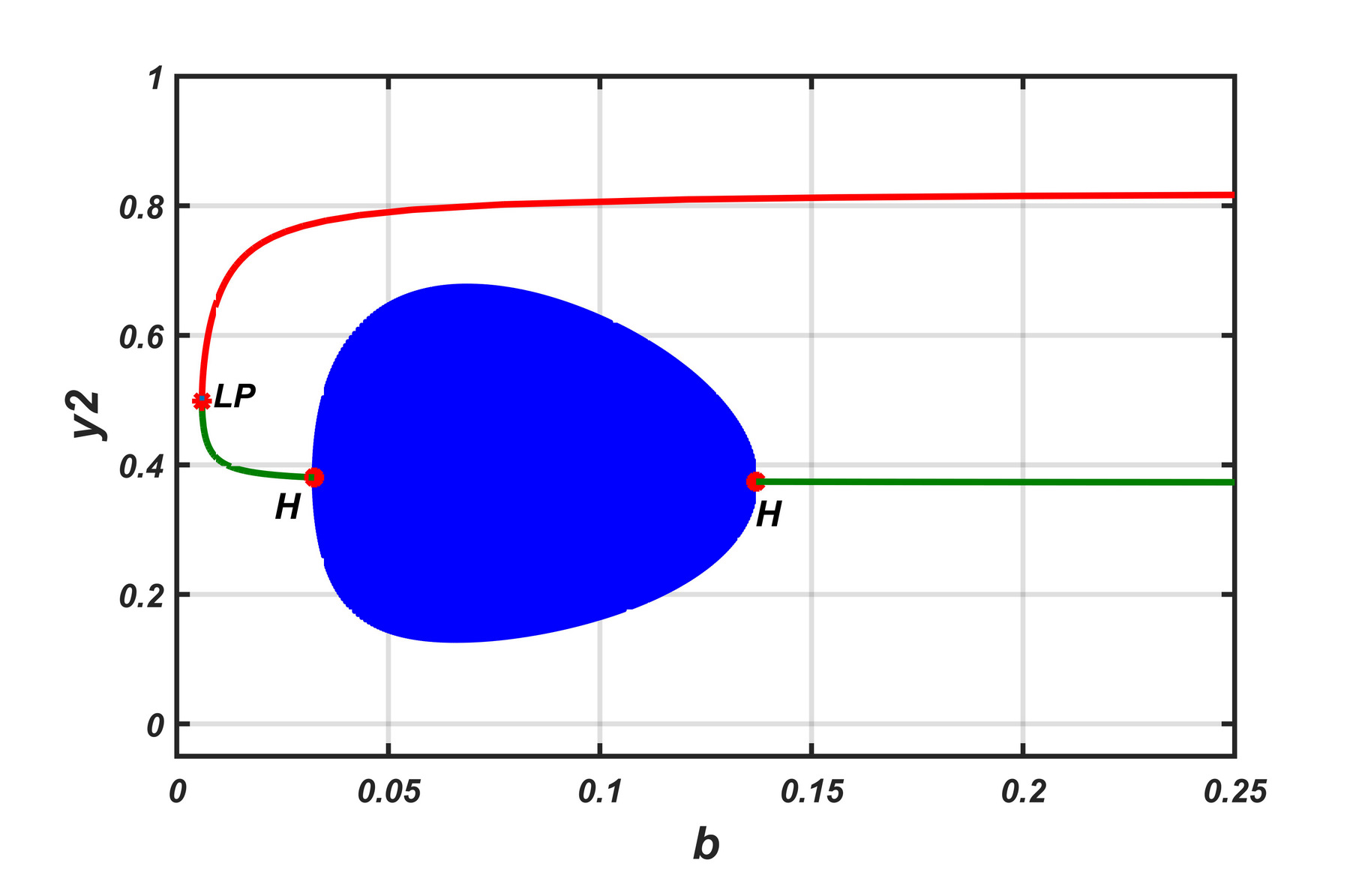}
         \caption{}
    \end{subfigure}\hfill
    \begin{subfigure}[b]{0.5\textwidth}
    \centering
         \includegraphics[width=\textwidth]{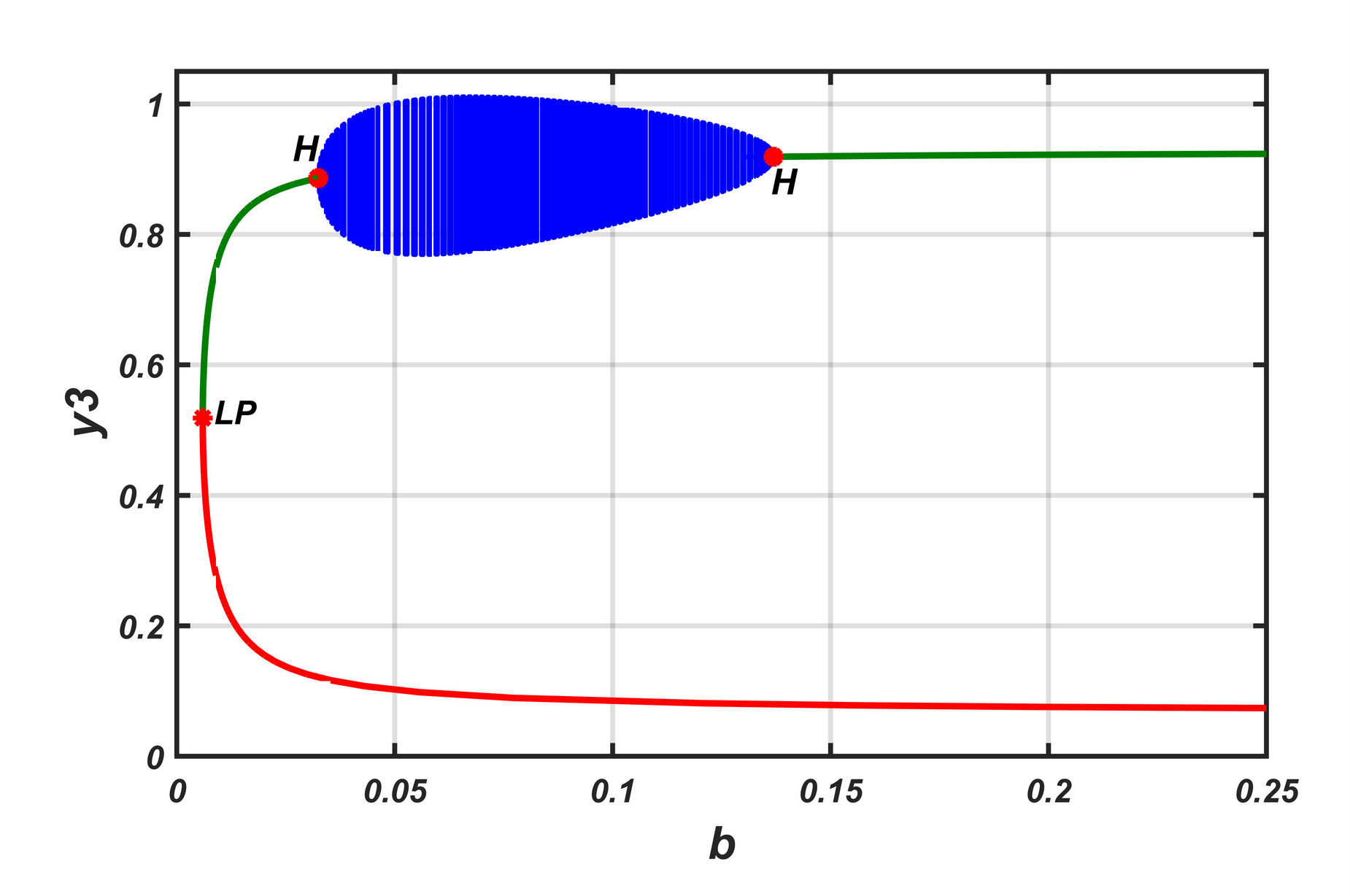}
         \caption{}
    \end{subfigure}
    \caption{Equilibrium states of (a) prey population, (b) infant predator population, (c) juvenile predator population, (d) adult predator population depicting bifurcations with respect to the transition rate of infant predator into juvenile stage, $b$.}
    \label{fig-bho}
\end{figure}
In figure \ref{u-lh}, the red line denotes unstability, the green line denotes stability, and the blue line denotes the amplitude of fluctuation at that particular parametric value of $u$. An illustration alike figure \ref{u-lh} can be espied for the parameter $a_2$, the predation rate of adult predators. While varying $a_2$, the saddle node bifurcation is experienced by the system at $a_2=0.274931$, the Hopf bifurcation point is detected at $a_2=0.810103$.\\
\par Figure \ref{a12} illustrates the change in dynamical behavior of the equations of bio-system (\ref{ma-eq}) with variations in the predation rate of juvenile predators ($a_1$) and the transition rate of juvenile predators ($c$), by continuing  the curves of saddle node bifurcation and Hopf bifurcation. These two curves partition the whole bi-parametric region into sub-regions.  The Hopf bifurcation curve separates the unstable compresent equilibrium region ($region\,  III$) from the stable one. The Hopf curve being a supercritical one, there exists stable limit cycle in the neighborhood of the curve. The region $IV$ bounded between the Hopf curve and the saddle node curve is the area of stability of compresent equilibrium point $E_4$, after crossing the saddle node curve and stepping into region $I$, the stability of $E_4$ vanishes because of collision with the unstable equilibrium $E_3$. Furthermore, the two bifurcation curves meet tangentially at $(a_1,c)=( 0.149588, 0.018589)$, called Bogdanov-Takens ($BT$) point. Because of $BT$ point, the non-saddle equilibrium $E_4$ which undergoes Hopf bifurcation, loses the occurrence of limit cycle via the saddle node bifurcation, in addition to that $BT$ point causes the separation of the saddle node curve into two parts: stable curve and unstable curve. Hence, in region $II$, neither the occurrence of limit cycle nor the stability of $E_4$ is found. Similarly, in figure \ref{a11}, the Hopf curve and saddle node curve meet tangentially at Bogdanov-Takens point, $(a_1,a_2)=( 0.019133, 0.274931)$, thereby making region $II$ synonymous with region $II$ of figure \ref{a12}. Also, in figure \ref{a11}, there exists a generalized Hopf bifurcation point (also known as Bautin bifurcation point) at  $(a_1,a_2)=( 0.042954, 0.379816)$, where the first Lyapunov coefficient vanishes. A division of the Hopf curve into supercritical part and subcritical part, is done by Bautin bifurcation point. For higher values of $a_1$ and $a_2$ (i.e., at the right hand side of Bautin bifurcation point) the Hopf curve is supercritical. Colossally, like the previous scenario, region $III$ is the area of unstability of compresent equilibrium, and regions $I$ and $IV$ are the regions of stability of compresent equilibrium $E_4$, while $E_4$ loses its stability via saddle node bifurcation curve in region $V$. The horizontal line corresponding to saddle node bifurcation curve in figures signifies that the values of both transition rates of juvenile predator ($c$) or predation rate of adult predator ($a_2$) remain the same in correspondence to all values of $a_1$, the predation rate of juvenile predators, i.e., $a_1$ is not responsible for the occurrence of saddle node bifurcation.  However, it should be noted that for a completely different set of values, $a_1$ may cause the said bifurcation, by reason of the presence of the parameter in the condition given in theorem \ref{thm-snb} for saddle node bifurcation.\\
\par The local bifurcations while varying the transition rate of infant predators, $b$ are visualized in figure \ref{fig-bho}. The equations of bio-system changes its stability through Hopf bifurcation not once but twice, the bifurcating points being $b=0.032488$ and $b=0.136940$ whose Lyapunov coefficients are $-1.619062e^{-02} $ and $ -7.177191e^{-02}$ respectively. Hence, both the points are supercritical Hopf bifurcation points. For $b \in (0.032488, 0.136940)$, there is the existence of a stable limit cycle whose amplitude increases and then gradually decreases as the value of $b$ travels from one bifurcation point to another, forming a ‘balloon’ like structure as can be seen in the figure \ref{fig-bho}. For $b>0.136940$ and $b_{s1}<b<0.032488$, compresent equilibrium $E_4$ is stable, where $b_{s1}=0.005977$ is the saddle node bifurcation point.  The solid red and green lines denote unstable and stable equilibriums respectively, while the blue lines depict the amplitudes.
\section{Discussion and conclusion}
Juveniles of the predator have been chronicled to hunt their prey. But predation without any prior experience, i.e., without knowing the prey specific approach would be a hellacious venture. Concomitantly, the prey oft-times counterattacks as a desperate remedy for intercepting their demise, thereby effectuating injuries, which are induced particularly in the inexperienced predator.  Also, as evidenced from all the documentations, there is always a period of time in the beginning of the predator’s life where they are completely reliant on their parents for survival. Thereupon, while considering predation by immature predators, it would be perspicuous for the predator species to be structured into three stages (infant, juvenile, and adult), while simultaneously incorporating injury, as the baleful consequence of predation by juveniles, ensuing a delay in the maturation rate for the juvenile predator into the adult stage. In this paper, equations of bio-system representing this very scenario has been assembled, thus providing a much more detailed and richer dynamic of the prey-predator interaction. \\
\par The boundedness of the equations of bio-system (\ref{ma-eq}) is authenticated, proving its feasibility in the real-life scenario, cause none of the species should be able to grow exponentially owing to the constraints of natural resources.  The system of equations (\ref{ma-eq}) harbours four biologically viable equilibrium points: one each of extinction equilibrium and prey-only equilibrium, and two of compresent equilibrium points. With the intent of studying the behavior of the system in the proximity of the equilibrium points, the conditions for local and global stability are examined for the equilibrium points in Theorems \ref{thm-lo-pf}, \ref{thm-glo-pf}, \ref{sta-co} and, \ref{thm-glo-co}, which are then authenticated by visualizing graphically in figures \ref{lo-axi}, and \ref{lo-xyz}. \\
\par An interesting state of affairs known as `bi-stability' (a situation with multiple attractors for the system of equations) is observed between compresent equilibrium and prey-only equilibrium points in figures \ref{b-equ}, \ref{bi-stab}. The latter equilibrium point is obtained only with a considerably low initial population of one or more stages of the predator species. It is noteworthy that even with negligible initial prey population, the system of equations instead of settling into extinction equilibrium, goes towards compresent equilibrium, which figure \ref{bi3} shows happening, by virtue of rapidly decreasing predator population (ascribable to inanition), in parallelism the prey population rapidly increases, before both of the species converges to the compresent equilibrium point. This circumstance is known as ‘bloom’ phenomenon.\\
\par Analysis of bifurcations is a prerequisite for uncovering the rich dynamics of the equations of bio-system (\ref{ma-eq}),  because it enables the determination of the transition points while tracing the paths of steady as well as unstable equilibriums. The equations of bio-system exhibits the occurrence of Hopf bifurcation around compresent equilibrium, i.e., the transition point where the said equilibrium point exchanges its stability with either a stable or unstable limit cycle (Theorem \ref{ho-thm}). The analytical methodology used to determine the stability and direction of the thus occurred limit cycle is shown Theorem \ref{ho-dir-thm} and then verified numerically through figure \ref{c-ho}. Likewise, the manifestation of saddle node bifurcation is corroborated both analytically (Theorem \ref{thm-snb}) and graphically (figure \ref{b-equ}).  The expression $a_2 b u (a_3-c)= (b + d_1)(a_3 - c - d_2)d_3$ denotes the surface on which transcritical bifurcation occurs (Theorem \ref{thm-tb}), the bifurcation responsible for the interchangeability of stability between the equilibrium points. The transcritical bifurcation is found to be existent between the prey-only equilibrium and one of the compresent equilibrium. As can be witnessed from the expression above, all the parameters are responsible for causing this bifurcation except for $a_1$, the predation rate of juveniles.\\
\par When the influence of transition rate of juvenile predators into the adult stage ($c$) is investigated, an exchange of stability from co-existing population to the stability of prey-only population is witnessed after decreasing the value of $c$ after a threshold point because of transcritical bifurcation for a certain parametric set of values (table 1). While for another set of values for the parameters i.e., table 2, the compresent equilibrium $E_4$ is stable, but increasing the parametric value of $c$ after the threshold value of Hopf bifurcation point, each and every stage of the predator population and the prey population, incapable of abiding by a specific population size due to the unstability of $E_4$, keeps on oscillating (figure \ref{c-ho}). A similar sequence of events can be espied for the parameter $a_3$, the rate of maturity delay of juvenile predators, but in the reverse orientation, for example, a decrease in parametric value after the Hopf bifurcation point leads to the birth of the limit cycle around the compresent equilibrium (figure \ref{a3-2dim}). Also, the existence of the fold bifurcation point where the unstable compresent equilibrium $E_3$  and stable compresent equilibrium $E_4$ collide, can be found for both $c$ and $a_3$. Both of these parameters combined together as $(c-a_3x)$ gives the comprehensive maturation rate of juvenile predators into the adult stage. Hence, biologically we can surmise that when this rate increases, the juvenile population decreases leading to lesser juvenile predator maturing in the next instance, correspondingly the adult predator would increase at first but their number too would gradually decrease both due to lesser recruitment and lesser availability of prey (decreased number of prey is attributable to more of them being killed because of increased adult predator population). Then when the biomass of the predator population would be less, the prey species would increase their number, resulting in the predator population following suit because of the now availability of food, and thus creating the oscillations. On the contrary, with low comprehensive maturation rate, the predator population not being able to survive goes towards extinction, i.e.  the equilibrium point $E_2(1,0,0,0)$. 
\par Predation by juveniles can be appraised as ancillary, owing to the fact that it doesn’t equate to proliferation of the species, but rather to amplification of expertise of the juveniles while hunting for their later stage.  Varying this parameter can cause Hopf-bifurcation in the system, since in figures \ref{bipara1}, \ref{bipara2}, for different values of $c$ and $a_2$, the Hopf bifurcating point of $a_1$ too differs (when the values of other parameters are taken from table 2). For values lesser than the threshold value of $a_1$, the equations of bio-system has unstable compresent equilibrium. That is for lesser value of $a_1$, more prey survives, which corresponds to a decrease in comprehensive maturation rate of juvenile predators ($c-a_3x$), the biological implication of which can be comprehended as the juveniles are not able to kill the prey with the reason of the anti-predator behavior of the prey being too strong, which also results into more juveniles getting injured. Also, since $a_1$ doesn’t possess the ascendancy in causing transcritical bifurcation, the prey-only equilibrium remains unstable.  
\par From a new vantage point, a stage in the predator species is introduced, identified in the present study as infant predators, which is similar to that of the immature predators as is considered in the conventional prey- predator model with stage structure in predator \cite{khajan1,khajan2}. Their transition rate into the juvenile stage is extensively investigated through graphical visualization. A decrease in $b$ is espied to cause the extinction of the predator species, that is the attainment of prey-only equilibrium (figure \ref{b-equ}). Biologically, this can be perceived as- with lesser predator individuals attaining the juvenile stage leads to a lesser number of adult predators, which corresponds to frivolous birth of infants, thereby the eventual extinction.  Also, as is espied in figure \ref{fig-bho}, an increment  of this transition rate beyond a threshold value results in destabilisation of compresent equilibrium, similar to the scenario of the birth of oscillations in all the populations with the  increment of transition rate of juvenile predators,  the solution trajectory keeps on encircling the equilibrium point beyond a threshold parametric value of the transition rate of infant predators. But the circumference of the solution trajectory at first increases with the increment of the parametric value, and then again decreases due to the manifestation of another supercritical Hopf bifurcation, ensuing the stability of compresent equilibrium point all over again. Therefore, there exists a neighbourhood in the parametric value of the infants' transition rate where none of the populations can conform to its optimal number, that is, after crossing the second threshold parametric value, both of the species can obtain their respective optimal populations 
\par The existence of the predators obviously relies a lot on the predation rate of adult predators ($a_2$), and the transformation of consumed prey into infant predators ($u$), which can be clearly seen with the existence of supercritical Hopf bifurcation and saddle node bifurcation with respect to these two parameters.  Furthermore, diminishing either of them beyond a threshold value insinuates non-viability of the subsistence of the predator species; contrarily, an increment beyond a certain value interdicts the existence of both the species. 
\par Occurrence of bifurcation of equilibria in a two-parameter region such as cusp bifurcation point, Bogdanov Takens point, and Bautin bifurcation point implies the richer dynamics of the equations of bio-system. The ascendancy of these bifurcation points over the system is shown in figures \ref{bipara}, and \ref{fig-a1} and explained in section \ref{numeri}.\\
\par With no biological data being accessible, the parametric values are based on hypothetical data, but the analysis is hoped to be useful for researchers doing their work on related fields with experimental data.  For further analysis, other instances such as refuge, fear effect, group hunting, and others can be incorporated in the current equations of bio-system (\ref{ma-eq}).\\ \\
\textbf{Competing interests: } The authors declare no competing interests.


\begin{thebibliography}{99}
\bibitem{p1}  Lotka, A. J. (1910). Contribution to the theory of periodic reactions. The Journal of Physical Chemistry, 14(3), 271-274.
\bibitem{p2} Volterra, V. (1927). Fluctuations in the abundance of a species considered mathematically. Nature, 119(2983), 12-13.
\bibitem{p3}  Gurney, W. (1986). The systematic formulation of models of stage-structured populations.
 \bibitem{p4} Nisbet, R. M., \& Gurney, W. S. C. (1984). “Stage-structure” models of uniform larval competition. In Mathematical ecology (pp. 97-113). Springer, Berlin, Heidelberg.
\bibitem{p5} Nisbet, R. M., Gurney, W. S. C., \& Metz, J. A. J. (1989). Stage structure models applied in evolutionary ecology. In Applied Mathematical Ecology (pp. 428-449). Springer, Berlin, Heidelberg.
\bibitem{p6}  Aiello, W. G., \& Freedman, H. I. (1990). A time-delay model of single-species growth with stage structure. Mathematical biosciences, 101(2), 139-153.
\bibitem{p7} Aiello, W. G., Freedman, H. I., \& Wu, J. (1992). Analysis of a model representing stage-structured population growth with state-dependent time delay. SIAM Journal on Applied Mathematics, 52(3), 855-869.
\bibitem{puma} Elbroch, L. M., Feltner, J., \& Quigley, H. B. (2017). Stage‐dependent puma predation on dangerous prey. Journal of Zoology, 302(3), 164-170.
\bibitem{snake} Farrell, T. M., Smiley-Walters, S. A., \& McColl, D. E. (2018). Prey species influences foraging behaviors: Rattlesnake (Sistrurus miliarius) predation on Little Brown Skinks (Scincella lateralis) and Giant Centipedes (Scolopendra viridis). Journal of Herpetology, 52(2), 156-161.
\bibitem{caiman} Figueiredo, A., Alves-Martins, N., \& Nogueira-Costa, P. (2021). Predation attempt by the Spectacled Caiman, Caiman crocodilus (Linnaeus, 1758), on the microhylid Elachistocleis carvalhoi Caramaschi, 2010 in the southeastern Amazon of Brazil. Herpetology Notes, 14, 1227-1229.
\bibitem{spider}Gajski, D., Petráková, L., \& Pekár, S. (2020). Ant-eating spider maintains specialist diet throughout its ontogeny. Journal of Zoology, 311(3), 155-163.
\bibitem{snail} Gosselin, L. A., \& Chia, F. S. (1996). Prey selection by inexperienced predators: do early juvenile snails maximize net energy gains on their first attack?. Journal of experimental marine biology and ecology, 199(1), 45-58.
\bibitem{clinus} Gibbons, M. J. (1988). Impact of predation by juvenile Clinus superciliosus on phytal meiofauna: are fish important as predators?. Marine ecology progress series, 45, 13-22.
\bibitem{raptor}Zuluaga, S., Vargas, F. H., Aráoz, R., \& Grande, J. M. (2022). Main aerial top predator of the Andean Montane Forest copes with fragmentation, but may be paying a high cost. Global Ecology and Conservation, 37, e02174.
\bibitem{shark} Dhellemmes, F., Smukall, M. J., Guttridge, T. L., Krause, J., \& Hussey, N. E. (2021). Predator abundance drives the association between exploratory personality and foraging habitat risk in a wild marine meso‐predator. Functional Ecology, 35(9), 1972-1984.
\bibitem{a1} Huang, R. K., Webber, Q. M., Laforge, M. P., Robitaille, A. L., Bonar, M., Balluffi-Fry, J., Zabihi-Seissan, S., \& Vander Wal, E. (2021). Coyote (Canis latrans) diet and spatial co-occurrence with woodland caribou (Rangifer tarandus caribou). Canadian Journal of Zoology, 99(5), 391-399.
\bibitem{a2} MacNulty, D. R. (2002). The predatory sequence and the influence of injury risk on hunting behavior in the wolf (Doctoral dissertation, University of Minnesota).
\bibitem{a3} Brown, J. S., Embar, K., Hancock, E., \& Kotler, B. P. (2016). Predators risk injury too: the evolution of derring-do in a predator–prey foraging game. Israel Journal of Ecology and Evolution, 62(3-4), 196-204.
\bibitem{a4} Miranda, E. B., Menezes, J. F. D., \& Rheingantz, M. L. (2016). Reptiles as principal prey? Adaptations for durophagy and prey selection by jaguar (Panthera onca). Journal of Natural History, 50(31-32), 2021-2035.
\bibitem{a5}Pokharel, P., Sippel, M., Vilcinskas, A., \& Petschenka, G. (2020). Defense of milkweed bugs (Heteroptera: Lygaeinae) against predatory lacewing larvae depends on structural differences of sequestered cardenolides. Insects, 11(8), 485.
\bibitem{muk} Mukherjee, S., \& Heithaus, M. R. (2013). Dangerous prey and daring predators: a review. Biological reviews, 88(3), 550-563.
\bibitem{b1} Fernández Moya, S., Iglesias Pastrana, C., Marín Navas, C., Ruíz Aguilera, M. J., Delgado Bermejo, J. V., \& Navas González, F. J. (2021). The Winner Takes it All: Risk Factors and Bayesian Modelling of the Probability of Success in Escaping from Big Cat Predation. Animals, 12(1), 51.
\bibitem{ex1}  Escalante, I. (2015). Predatory behaviour is plastic according to prey difficulty in naïve spiderlings. Journal of insect behavior, 28(6), 635-650.
\bibitem{ex2}  García, L. F., Franco, V., Robledo-Ospina, L. E., Viera, C., Lacava, M., \& Willemart, R. H. (2016). The predation strategy of the recluse spider Loxosceles rufipes (Lucas, 1834) against four prey species. Journal of Insect Behavior, 29(5), 515-526.
\bibitem{ex3} Perlman, Y., \& Tsurim, I. (2008). Daring, risk assessment and body condition interactions in steppe buzzards Buteo buteo vulpinus. Journal of Avian Biology, 39(2), 226-228.
\bibitem{ppr1} Kaushik, R., \& Banerjee, S. (2022). Predator–prey system with multiple delays: prey’s countermeasures against juvenile predators in the predator–prey conflict. Journal of Applied Mathematics and Computing, 68(4), 2235-2265.
\bibitem{ppr2} Kaushik, R., \& Banerjee, S. (2021). Predator-prey system: Prey’s counter-attack on juvenile predators shows opposite side of the same ecological coin. Applied Mathematics and Computation, 388, 125530.
\bibitem{ppr3} Li, J., Liu, X., \& Wei, C. (2022). Stationary distribution of a stage-structure predator–prey model with prey’s counter-attack and higher-order perturbations. Applied Mathematics Letters, 129, 107921.
\bibitem{ppr4} Bhattacharjee, D., Roy, T., Acharjee, S., \& Dutta, T. K. (2022). Stage structured prey-predator model incorporating mortal peril consequential to inefficiency and habitat complexity in juvenile hunting. Heliyon, 8, e11365. 
\bibitem{ppr5} Yao, P., Wang, Z., \& Wang, L. (2022). Stability Analysis of a Ratio-Dependent Predator-Prey Model. Journal of Mathematics, 2022, 46052.
\bibitem{ppr6} Jawad, S., \& Naji, R. K. (2022). The Influence of Stage Structure and Prey Refuge on the Stability of the Predator-Prey Model. I. J. Engineering and Manufacturing, 3, 51-59.
\bibitem{matc} Dhooge, A., Govaerts, W., Kuznetsov, Y. A., Meijer, H. G. E., \& Sautois, B. (2008). New features of the software MatCont for bifurcation analysis of dynamical systems. Mathematical and Computer Modelling of Dynamical Systems, 14(2), 147-175.
\bibitem{soto}Perko, L. (2013). Differential equations and dynamical systems (Vol. 7). Springer Science \& Business Media.
\bibitem{hassa} Hassard, B. D., Hassard, B. D., Kazarinoff, N. D., Wan, Y. H., \& Wan, Y. W. (1981). Theory and applications of Hopf bifurcation (Vol. 41). Cambridge, New York, Cambridge University Press.
\bibitem{khajan1} Khajanchi, S. (2017). Modeling the dynamics of stage-structure predator-prey system with Monod–Haldane type response function. Applied Mathematics and Computation, 302, 122-143.
\bibitem{khajan2} Khajanchi, S. (2014). Dynamic behavior of a Beddington–DeAngelis type stage structured predator–prey model. Applied Mathematics and Computation, 244, 344-360.
\bibitem{p-ex1} Georgescu, P., \& Hsieh, Y. H. (2007). Global dynamics of a predator-prey model with stage structure for the predator. SIAM Journal on Applied Mathematics, 67(5), 1379-1395.
\bibitem{p-ex2} Liu, Q., Jiang, D., Hayat, T., \& Alsaedi, A. (2018). Dynamics of a stochastic predator–prey model with stage structure for predator and Holling type II functional response. Journal of Nonlinear Science, 28, 1151-1187.
\end{thebibliography}
\end{document}